\input ./style/arxiv-general.cfg
\documentclass[aop,MSNbibl,dvips]{arximspdf}
\makeatletter
   \@ifpackageloaded{graphicx}{}{\usepackage{graphicx}}
\makeatother
\doi{10.1214/15-AOP1022}% Updated by VTEXPTS2LaTeX.exe, 21.05.2015
\volume{44}
\issue{3}
\pubyear{2016}
\firstpage{2264}
\lastpage{2348}
\docsubty{FLA}

\makeatletter
\newcommand{\lleft}{\left}
\newcommand{\rrvert}{\vert}
\newcommand{\rright}{\right}
\newcommand{\rrVert}{\Vert}
\newcommand{\llvert}{\vert}
\newcommand{\llVert}{\Vert}

\def\vafrac#1#2{(#1)/(#2)}
\def\sklfrac#1#2{(#1/#2)}
\def\sklvfrac#1#2{((#1)/#2)}
\def\sklafrac#1#2{(#1/(#2))}
\def\sklvafrac#1#2{((#1)/(#2))}
\renewcommand{\mid}{|}
\newtheorem{theorem}{Theorem}
\newtheorem{lemma}{Lemma}[section]

\newtheorem{corollary}[lemma]{Corollary}
\newtheorem{proposition}[lemma]{Proposition}
\newproclaim{definition}[lemma]{Definition}
\newproclaim{assumption}{Assumption}
\newproclaim{Remark}[lemma]{Remark}
\newcommand{\R}{\mathbb{R}}
\newcommand{\C}{\mathbb{C}}
\newcommand{\N}{\mathbb{N}}
\newcommand{\p}{\mathbb{P}}
\newcommand{\K}{\mathrm{K}}
\newcommand{\A}{\mathrm{A}}
\newcommand{\B}{\mathrm{B}}
\newcommand{\E}{\mathrm{E}}
\newcommand{\Ures}{\Upsilon_{\mathrm{res}}}
\newcommand{\Ttilde}{\widetilde\Theta}
\newcommand{\supp}{\operatorname{Supp}}
\newcommand{\dist}{\operatorname{dist}}
\def\Tr{\operatorname{Tr}}
\def\re{\operatorname{Re}}
\newcommand{\Ai}{\mathrm{Ai}}
\newcommand{\Aii}{\mathrm{Ai}}
\renewcommand{\d}{\mathrm{d}}
\makeatother

\begin{document}
\begin{frontmatter}

%\dochead{}
\title{Large complex correlated Wishart matrices: Fluctuations and
asymptotic independence at~the~edges}
\runtitle{Complex correlated Wishart matrices}

\begin{aug}
% Corresponding author: Adrien Hardy - ahardy@kth.se% Updated by
%VTEXPTS2LaTeX.exe, 21.05.2015 09:49
\author[A]{\fnms{Walid}~\snm{Hachem}\ead[label=e1]{walid.hachem@telecom-paristech.fr}\thanksref{T1}},
\author[B]{\fnms{Adrien}~\snm{Hardy}\corref{}\ead[label=e2]{ahardy@kth.se}\thanksref{T2}}
\and
\author[C]{\fnms{Jamal}~\snm{Najim}\ead[label=e3]{najim@univ-mlv.fr}\thanksref{T1}}
\runauthor{W. Hachem, A. Hardy and J. Najim}
\affiliation{T\'{e}l\'{e}com ParisTech, KTH Royal Institute of
Technology and Universit\'e Paris-Est}
%\dedicated{}
\address[A]{W. Hachem\\
CNRS LTCI\\
T\'{e}l\'{e}com ParisTech\\
46 rue Barrault\\
75013 Paris Cedex 13\\
France\\
\printead{e1}}
\address[B]{A. Hardy\\
Department of Mathematics\\
KTH Royal Institute of Technology\\
Lindstedtsv\"agen 25\\
10044 Stockholm\\
Sweden\\
\printead{e2}}
\address[C]{J. Najim\\
CNRS LIGM\\
Universit\'e Paris-Est\\
Cit\'e Descartes\\
5 Boulevard Descartes\\
Champs sur Marne\\
77 454 Marne la Vall\'ee Cedex 2\\
France\\
\printead{e3}}
\end{aug}
\thankstext{T1}{Supported in part by the program ``mod\`eles num\'eriques'' of the French Agence Nationale de la Recherche
under the Grant ANR-12-MONU-0003 (project DIONISOS).}
\thankstext{T2}{Supported by the KU Leuven research Grant OT/12/73 and the Grant KAW 2010.0063 from the
Knut and Alice Wallenberg Foundation.}

% HISTORY:
%
\received{\smonth{9} \syear{2014}}% Updated by VTEXPTS2LaTeX.exe,
%21.05.2015 09:49
%
\revised{\smonth{2} \syear{2015}}% Updated by VTEXPTS2LaTeX.exe,
%21.05.2015 09:49

% ABSTRACT
%
\begin{abstract}
We study the asymptotic behavior of eigenvalues of large complex
correlated Wishart matrices
at the edges of the limiting spectrum. In this setting, the support of
the limiting
eigenvalue distribution may have several connected components. Under
mild conditions for the population matrices, we show that for every
generic positive edge of that support, there exists an extremal
eigenvalue which converges almost surely toward that edge and
fluctuates according to the Tracy--Widom law at the scale $N^{2/3}$.
Moreover, given several generic positive edges, we establish that the
associated extremal eigenvalue fluctuations are asymptotically
independent. Finally, when the leftmost edge is the origin (hard edge),
the fluctuations of the smallest eigenvalue are described by mean of
the Bessel kernel at the scale $N^2$.
\end{abstract}

% KEYWORDS
% Pirmas kwd is didziosios raides
%
\end{frontmatter}

\tableofcontents[level=2]
%\maxcsec{10}
%\maxcsubsec{10}
%\setattribute{tocline}{skip}{\space}

\setcounter{footnote}{2}
%s1 #&#
\section{Introduction}\label{sec1}
Correlated Wishart matrices and more generally empirical covariance
matrices are ubiquitous models in applied mathematics. After Mar\v
cenko and Pastur's seminal contribution \cite{mar-and-pas-67}, a
systematic study of their large dimension properties has been
undertaken (see, e.g., \cite
{book-bai-silverstein,book-pastur-shcherbina} and the many references
therein), which found many applications, for example, in multivariate
statistics \cite{book-bai-chen-liang-2009}, electrical engineering
\cite{book-couillet-debbah}, mathematical finance \cite
{bouchaud-noise-dressing-99,guhr-credit-risk-2014}, etc.

Now that many global properties of their spectrum are well understood
(cf. \cite
{bai-silverstein-98,silverstein-exact-separation-99,bai-silverstein-aop-2004,najim-yao-2013-preprint,silverstein-95}),
attention has shifted to local properties (cf. \cite
{BBP-2005,elkaroui-07-TW,bloemendal-virag-2013-PTRF}, etc.) and their
underlying universal phenomenas; cf. \cite{kuijlaars-handbook-2011}
and references therein.

The main contribution of this article is to provide a local analysis of
the spectrum of large complex correlated Wishart matrices near the
edges of the limiting support: it is well known that such random
Hermitian matrices have a real spectrum whose limiting support may
display several disjoint intervals. Beside the behavior of the largest
and smallest random eigenvalues, we investigate here the fluctuations
of the eigenvalues that converge to any endpoint of the limiting
support. These eigenvalues are referred to as \textit{extremal
eigenvalues}, for which we shall provide a precise definition later.

\begin{longlist}
\item[\textit{The model}.] Let ${\mathbf X}_N$ be a $N\times n$ matrix with
independent and identically distributed (i.i.d.) standard complex
Gaussian entries ${\mathcal N}_\C(0,1)$, and let ${\bolds\Sigma}_N$ be
a $n\times n$ deterministic positive definite Hermitian matrix. The
random matrix of interest here is the $N\times N$ matrix
%
%e1 #&#
\begin{equation}
\label{mainmatrixmodel} {\mathbf M}_N = \frac{1}N {\mathbf X}_N {
\bolds\Sigma}_N {\mathbf X}_N^*
\end{equation}
which has $N$ nonnegative eigenvalues $0\leq x_1\leq\cdots\leq x_N$,
but which may be of different nature: $\min(n,N)$ of them are
nonnegative random (i.e., nondeterministic) eigenvalues, while the
other $N-\min(n,N)$ eigenvalues are deterministic and equal to zero. A
companion matrix of interest is the $n\times n$ sample covariance matrix
%
%e2 #&#
\begin{equation}
\widetilde{\mathbf M}_N = \frac{1}N {\bolds\Sigma}_N^{1/2}
{\mathbf X}_N^* {\mathbf X}_N {\bolds\Sigma}_N^{1/2},
\end{equation}
which models the empirical covariance of a sample of $N$ independent
observations
\[
\bigl\{{\bolds\Sigma}_N^{1/2} {\bigl[{\mathbf X}_N^*
\bigr]}_{k}, 1\leq k\leq N \bigr\},
\]
where $[{\mathbf X}_N^*]_{k}$ stands for the $k$th column of ${\mathbf X}_N^*$,
with population covariance matrix~${\bolds\Sigma}_N$. Indeed, matrices
$\mathbf M_N$ and $\widetilde{\mathbf M}_N$ share the same nonnull
eigenvalues with the associated multiplicities.

We shall consider the asymptotic regime where $n=n(N)$, $N\to\infty$ and
%
%e3 #&#
\begin{equation}
\label{evdistrMN} \lim_{N\rightarrow\infty} \frac{n}{N}=\gamma\in(0,
\infty).
\end{equation}
This regime will be simply referred to as $N\to\infty$ in the sequel.

The random matrix $\mathbf M_N$ can also be interpreted as a
multiplicative deformation of the Laguerre unitary ensemble (LUE) and
is related to multiple Laguerre polynomials. A close matrix model is
the additive deformation of the Gaussian unitary ensemble (GUE), also
known as GUE with an external source; it involves multiple Hermite
polynomials instead. For further information, see \cite
{bleher-kuijlaars-05} and references therein. Capitaine and P\'ech\'e
\cite{capitaine-peche-2014-preprint} recently studied the fluctuations
of extremal eigenvalues for this model.

We now briefly review the literature and present our contribution.
\end{longlist}

\begin{longlist}
\item[\textit{Global regime}.]
Denote by $\mu_N$ the empirical distribution of the eigenvalues of~${\mathbf M}_N$, also called spectral measure (or distribution) of ${\mathbf
M}_N$ in the sequel. Namely,
\[
\mu_N =\frac{1}N \sum_{i=1}^N
\delta_{x_i},
\]
where $\delta_{x}$ is the Dirac measure at point $x$.
In the uncorrelated case where ${\bolds\Sigma}_N=I_n$,
it is well known \cite{mar-and-pas-67} that $\mu_N$ almost surely
(a.s.) converges weakly toward the Mar\v cenko--Pastur (\v MP)
distribution of parameter $\gamma$,
%
%e4 #&#
\begin{equation}
\label{MPlaw} \mu_{\mathrm{MP}}^\gamma(\d x)=(1-\gamma)^+
\delta_0 + \frac{1}{2\pi
x}\sqrt{(\frak b-x) (x-\frak a)}
\mathbf 1_{[\frak a,\frak
b]}(x)\,\d x,
\end{equation}
where $x^+=\max(x,0)$ and the endpoints of its support read $\frak
a=(1-\sqrt{\gamma})^2$ and $\frak b= (1+\sqrt{\gamma})^2$.

In the general case where $\bolds\Sigma_N$ is not the identity, say
with eigenvalues $0<\lambda_1\leq\cdots\leq\lambda_n$, a similar
result holds true \cite{silverstein-95} under the additional
assumption that the spectral measure
%
%e5 #&#
\begin{equation}
\nu_N=\frac{1}n \sum_{j=1}^n
\delta_{\lambda_j}
\end{equation}
of $\bolds\Sigma_N$ converges weakly toward a limiting distribution
$\nu$. In the latter case, the limit $\mu$ of $\mu_N$ only depends
on the limiting parameters $\gamma$ and $\nu$ but is no longer
explicit; this dependence $\mu=\mu(\gamma,\nu)$ will be indicated
when needed. However, its Cauchy--Stieltjes transform satisfies an
explicit fixed-point equation from which many properties of $\mu$ can
be inferred. For example, it is known that if $\nu(\{0\})=0$, then
%
%e6 #&#
\begin{equation}
\mu(\d x)=(1-\gamma)^+ \delta_0 + \rho(x)\,\d x,
\end{equation}
where $\rho(x)$ is a nonnegative and continuous function on
$(0,+\infty)$. Depending on the properties of $\gamma$ and $\nu$,
the support of $\rho(x)\,\d x$ may have several connected components;
see Section~\ref{secbackground} for more precise informations.
Alternatively, one can describe $\mu(\gamma,\nu)$ in terms of the
free multiplicative convolution of \v MP distribution~(\ref{MPlaw})
with $\nu$; see \cite{voiculescu-1991}. From now we shall refer to
the support of $\rho(x)\,\d x$ as the \textit{bulk} and to the endpoints of
its connected components as the \textit{edges}. Also, a~positive edge is
called \textit{soft edge} and the terminology \textit{hard edge} is here used
when the edge is the origin.
\end{longlist}

\begin{longlist}
\item[\textit{Left and right edges}.] We say that an edge $\frak a$ is a
left edge, respectively, $\frak b$ is a right edge, if for every
$\delta>0$ small enough,
\[
\int_{\frak a}^{\frak a+\delta}\rho(x)\,\d x>0,\quad\mbox{respectively,}\quad
\int_{\frak b-\delta}^\frak b\rho(x)\,\d x>0.
\]
The leftmost edge can be a soft edge or a hard edge depending on the
value of $\gamma$, as explained in Section~\ref{secbackground}. Of
course, any other left edge and any right edge are soft edges.
\end{longlist}

\begin{longlist}
\item[\textit{Local regime}: \textit{Behavior at the rightmost edge}.]
If ${\bolds\Sigma}_N$ is the identity, Geman \cite{geman-1980} proved
the a.s. convergence of the largest eigenvalue $x_{\max}$ of $\mathbf
M_N$ to the right edge of \v MP's bulk $\frak b=(1+\sqrt{\gamma})^2$,
for independent, not necessarily Gaussian, real entries of $\mathbf
X_N$. Johansson \cite{johansson-2000-shape}
established Tracy--Widom fluctuations for $x_{\max}$ at the scale
$N^{2/3}$ for complex Gaussian entries; Johnstone \cite
{johnstone-2001} established a similar result for real Gaussian
entries. Subsequent works \cite
{peche-2009-PTRF,pillai-yin-2014,soshnikov-2002,wang-2012-RMTA} then
relaxed the Gaussian assumption, illustrating a phenomenon of universality.

If $\bolds\Sigma_N$ is a finite-rank perturbation of the identity,
the limiting eigenvalue distribution is still given by \v MP
distribution (\ref{MPlaw}). Baik and Silverstein \cite
{baik-silverstein-2006} studied the limiting behavior of $x_{\max}$
for general entries. In the complex Gaussian case, Baik et al. \cite
{BBP-2005} thoroughly described the fluctuations of the largest
eigenvalues at the right edge and unveiled a remarkable phase
transition phenomenon (referred to as BBP phase transition in the
sequel). They
established that the convergence and fluctuations of $x_{\max}$ are
actually highly sensitive to the way $\nu_N$ converges to $\delta_1$.
More precisely, depending on the strength of the perturbation, they
established that deformed Tracy--Widom fluctuations near the right edge
$\frak b$ at the scale $N^{2/3}$ can arise, and that $x_{\max}$ may
also\vspace*{1pt} converge outside the bulk with Gaussian-like\footnote{By
Gaussian-like, we mean that the largest eigenvalue of ${\mathbf M}_N$, when
correctly centered and rescaled and when associated to a large
perturbation of the identity $\bolds\Sigma_N$ of finite multiplicity
$k$, asymptotically converges to the distribution of the largest
eigenvalue of a fixed $k\times k$ GUE.}
fluctuations at the scale $N^{1/2}$; in the latter case $x_{\max}$ is
referred to as an outlier. Thus, depending on the way $\nu_N$
converges toward its limit, the universality phenomenon may break down.
Finally, Bloemendal and Vir\'ag \cite
{bloemendal-virag-2013-PTRF,bloemendal-virag-2011-preprint} and Mo
\cite{Mo-2012} extended the results in \cite{BBP-2005} for real
Gaussian entries; see also \cite{bloemendal-et-al-preprint} for
further extensions.

For general $\bolds\Sigma_N$'s and complex Gaussian matrices, El
Karoui \cite{elkaroui-07-TW} ($n/N\leq1$) and then Onatski \cite
{onatski-2008-TW} ($n/N> 1$) followed the approach developed in \cite
{BBP-2005} to establish Tracy--Widom fluctuations for $x_{\max}$,
under mild conditions concerning $\bolds\Sigma_N$'s spectral measure
$\nu_N$ provided that the rightmost edge satisfies some regularity
condition. The Gaussian assumption has recently been relaxed by Bao et
al. \cite{pan-et-al-preprint} (the random variables remaining complex)
and the complex one, by Lee and Schnelli \cite{lee-schnelli-preprint}
who handle the real Gaussian case and also the real non-Gaussian case
for diagonal $\bolds\Sigma_N$'s. Knowles and Yin \cite
{knowles-yin-2014-preprint} extend \cite{lee-schnelli-preprint} to
general $\bolds\Sigma_N$'s; see also the comment on universality in
Section~\ref{secAIresults}.
\end{longlist}

\begin{longlist}
\item[\textit{Local regime}: \textit{Behavior at the leftmost edge}.] When $\bolds
\Sigma_N$ is the identity,
Bai and Yin \cite{bai-yin-AoP-1993} established the a.s. convergence
of the smallest eigenvalue $x_{\min}$ of $\mathbf M_N$ to \v MP's left
edge $\frak a=(1-\sqrt\gamma)^2$; see also \cite{book-bai-silverstein}, Chapter~5. The nature of the fluctuations of
$x_{\min}$ dramatically changes whether $\gamma=1$ (hard edge) or
$\gamma\neq1$ (soft edge). In the soft edge case, the fluctuations
remain of a Tracy--Widom nature; see Borodin and Forrester \cite
{borodin-forrester-2003} and further extensions by Feldheim and Sodin
\cite{feldheim-sodin-2010}. In the hard edge case, the fluctuations of
$x_{\min}$ arise at the scale $N^2$; if $n=N$, then the limiting
distribution follows the exponential law as shown by Edelman \cite
{edelman-1988-SIAM} (cf. \cite{tao-vu-2010-GAFA} for further
extensions), while if $n=N+\alpha$ with $\alpha$ independent of $N$,
then the limiting distribution has been described by Forrester \cite
{forrester-1993-bessel} with the help of Bessel kernels;
%\textcolor{red}{, the so-called hard-edge Tracy-Widom law}
see Section~\ref{secmain-results} for a precise definition. The
Gaussian assumption has been relaxed by Ben Arous and P\'ech\'e \cite{benarous-peche-2005}.
%[REF Universality of local eigenvalue statistics for sample covariance
%matrices, G. Ben Arous, S. P\'ech\'e, Communications of Pure and
%Applied Mathematics, vol 58, 10, pp 1316-1357, 2005]

To the best knowledge of the authors, no result for the fluctuations at
the leftmost edge in the general $\bolds\Sigma_N$ case is available
in the literature.
\end{longlist}

\begin{longlist}
\item[\textit{Local regime}: \textit{When $\nu$ is the weighted sum of two Dirac measures}.] When $\bolds\Sigma_N$ has exactly two fixed eigenvalues,
each with multiplicity of order $N$, a full asymptotic analysis is
known for the correlation kernel $\K_N(x,y)$ associated with the
eigenvalues of $\mathbf M_N$; see Sections~\ref{determinantalprocess}~and~\ref{kernelofWishart}. More precisely, around each edge a local
uniform convergence for $\K_N(x,y)$ has been obtained, using the
connection to multiple Laguerre polynomials, by Lysov and Wielonsky
\cite{lysov-2008}
%[V. Lysov, F. Wielonsky, Strong asymptotics for multiple Laguerre
%polynomials, Constr. Approx., 28 (2008), 61--111.]
and Mo \cite{mo-jmva-2010}.
%[M. Y. Mo, Universality in Complex Wishart ensembles: The 2 cut case,
%Preprint 41p http://arxiv.org/abs/0809.3750].
This provides a first step toward Tracy--Widom fluctuations.
\end{longlist}

\begin{longlist}
\item[\textit{Local regime}: \textit{Asymptotic independence}.] When $\bolds
\Sigma_N$ is the identity and $\gamma>1$ (and also in the case of the
GUE), Basor, Chen and Zhang \cite{basor-et-al-2012} proved that
$x_{\min}$ and $x_{\max}$, properly rescaled, are asymptotically
independent as $N\rightarrow\infty$. Their approach heavily relies on
orthogonal polynomials techniques, which are not available for complex
correlated Wishart matrices.
Using different techniques, the asymptotic independence for the GUE's
smallest and largest eigenvalues was also obtained by Bianchi et al.
\cite{bianchi-et-al-2010} and Bornemann \cite{bornemann-2010}.

Again, it seems there is no result concerning the asymptotic
independence for the extremal eigenvalues, even for the smallest and
largest eigenvalues, in the general $\bolds\Sigma_N$ case.
\end{longlist}

\begin{longlist}
\item[\textit{Main results}.] Recall the asymptotic regime (\ref
{evdistrMN}) of interest. We first state the main assumptions related
to matrix ${\mathbf M}_N$ [cf. (\ref{mainmatrixmodel})]
and then informally state the main results of the paper; pointers to
the precise definitions and statements are provided in the next paragraph.
\end{longlist}

%The assumptions we make on the matrix model \eqref{mainmatrixmodel}
%are the following ones.

%
%as1 #&#
\begin{assumption}
\label{assgauss}
The entries of $\mathbf X_N$ are i.i.d. standard complex Gaussian
random variables.
\end{assumption}

%
%as2 #&#
\begin{assumption} %\phantom{Let ${\bolds\Sigma}_N$ be a }
\label{assnu}
The following properties hold true:
\begin{longlist}[(2)]
\item[(1)]
The spectral measure $\nu_N$ of $\bolds\Sigma_N$ weakly converges
toward a
limiting probability distribution $\nu$ as $N\rightarrow\infty$.
\item[(2)]
The eigenvalues $0<\lambda_1\leq\cdots\leq\lambda_n$ of $\bolds
\Sigma_N$ stay in a compact subset of $(0,+\infty)$ which is
independent of $N$, namely,
%
%e7 #&#
\begin{equation}
\liminf_{N\rightarrow\infty}\lambda_1>0,\qquad \limsup
_{N\rightarrow\infty}\lambda_n<+\infty.
\end{equation}
In particular, $\nu(\{0\})=0$.
\end{longlist}
\end{assumption}

% {\red La version d'Adrien}
%
% Another important assumption is the fact that the considered edges
%need to be {\mathbf regular}. By this, we mean an edge which satisfies the
%regularity condition
% of Definition~\ref{defreg}. This condition essentially rules out
% pathological behaviors at the edges, such as deformed
% Tracy-Widom fluctuations (cf. \cite{BBP-2005}), but enables the
%appearance of outliers.\\
%
% {\red version de Jamal + commentaires}

Another important assumption is the fact that the considered edges need
to be \textit{regular}. By this, we mean an edge which satisfies the
regularity condition
of Definition~\ref{defreg}. This condition essentially rules out
pathological behaviors at edges, for example, when the limiting
eigenvalue density does not vanish like a square root. It does,
however, enable the appearance of outliers.

% {\red comments: il est vrai que techniquement notre m{\chr"C3}{
%\chr"A9}thode ne permet pas de traiter les deformed TW, mais c'est
%juste par limitation en page et en temps,
% pas par limitation technique. Si on {\chr"C3}{\chr"A9}crit
% "this condition rules out pathological behaviors at the edges such as
%deformed
% Tracy-Widom fluctuations" il me semble qu'on rentre plut{\chr"C3}{
%\chr"B4}t dans la 2{\chr"C3}{\chr"A8}me cat{\chr"C3}{\chr"A9}gorie.
%Autre chose, les gens {\chr"C3}{\chr"A0} qui j'ai parl{\chr"C3}{
%\chr"A9} de notre r{\chr"C3}{\chr"A9}sultat (Antti Knowles, Gregory
%Shcher) m'ont spontan{\chr"C3}{\chr"A9}ment dit "il faut $\sqrt{\left\vert x\right\vert }$
%at the edge" ce qui est impr{\chr"C3}{\chr"A9}cis et un peu faux mais
%quand m{\chr"C3}{\chr"AA}me qui est signifiant pour les gens. Je pense
%qu'on aurait int{\chr"C3}{\chr"A9}r{\chr"C3}{\chr"AA}t {\chr"C3}{
%\chr"A0} le mentionner.

% [ADRIEN: je te laisse arbitrer et d{\chr"C3}{\chr"A9}cider une fois
%pour toutes lors de ta relecture]

% }

%
%th1 #&#
\begin{theorem}\label{thebigtheorem} Let Assumptions~\ref{assgauss}~and~\ref{assnu} hold true. Then:
\begin{longlist}[(a)]
\item[(a)] Extremal eigenvalues: Given a regular right (resp.,
left) edge, there are perfectly located maximal (resp., minimal)
eigenvalues which converge a.s. toward this edge as $N\to\infty$;
these eigenvalues are called extremal eigenvalues.
\item[(b)] Tracy--Widom fluctuations: Given a regular right
(resp., left) soft edge, the associated extremal eigenvalue, properly
rescaled,\vspace*{1pt} converges in law to the Tracy--Widom distribution (resp.,
reversed Tracy--Widom distribution) at the scale $N^{2/3}$.
\item[(c)] Asymptotic independence: Given a finite family of
regular soft edges, the associated extremal eigenvalues, properly
rescaled, are asymptotically independent as $N\to\infty$.
\item[(d)] Hard edge fluctuations: In the case where $\gamma
=1$, the bulk displays a hard edge at 0. If $n=N+\alpha$ with $\alpha
\in\mathbb{Z}$ independent of $N$, then the fluctuations of the
smallest eigenvalue, properly rescaled, are described by mean of the
Bessel kernel with parameter $\alpha\in\mathbb{N}$ at the scale $N^2$.
%\textcolor{red}{
%the smallest eigenvalue, properly rescaled, is shown to converge to
%the hard-edge Tracy-Widom distribution with parameter $\alpha\in
%\mathbb{N}$ at the scale $N^2$.}
\end{longlist}
\end{theorem}

Close to our work is the recent paper by Capitaine and P\'ech\'e \cite
{capitaine-peche-2014-preprint} where the fluctuations of the extremal
eigenvalues for the additive deformation of the GUE are established,
that is the counterpart of part~(b) of Theorem~\ref{thebigtheorem},
together with Gaussian-like fluctuations for outliers and fluctuations
of the eigenvalue process at cusp points (i.e., when two bulks merge
together) with the appearance of the Pearcey process.
As the involved techniques are extremely model-dependent, the technical
difficulties are substantially different for the model under study.
The study of the fluctuations of the eigenvalue process at a cusp point
for large complex correlated Wishart matrices will appear elsewhere
\cite{HHN-cusp-preprint}.

Let us now briefly comment on Theorem~\ref{thebigtheorem}.

In part (a), we rely on results by Silverstein et al. \cite
{bai-silverstein-98,silverstein-exact-separation-99,silverstein-choi-1995}
on the support of limiting spectral distributions and on fine
asymptotic properties of the empirical spectrum to define regular edges
and to properly express the convergence of extremal eigenvalues.

In part (b), we first obtain an asymptotic
Fredholm determinantal representation of the extremal eigenvalues'
distribution and then
perform an asymptotic analysis of the associated kernels to prove convergence
toward the Airy kernel. The latter analysis is based on a steepest
descent analysis involving contours deformations.
Contrary to the analysis performed by Baik, Ben Arous and P\'ech\'e
\cite{BBP-2005}, El Karoui
\cite{elkaroui-07-TW} and Capitaine and P\'ech\'e \cite
{capitaine-peche-2014-preprint}, who work out explicit
deformed contours, our analysis relies on a more abstract argument
where the
existence of appropriate contours is obtained by mean of the maximum principle
for subharmonic functions. This argument has the advantage to work for
every regular right or left edge (and also for cusp points, cf. \cite
{HHN-cusp-preprint}) up to minor modifications. Let us also stress that
we do not follow the same strategy as in \cite
{BBP-2005,elkaroui-07-TW}, concerning the involved operators convergence.

%Part (b) is the most involved part of the paper. We first obtain an
%asymptotic
%determinantal representation of the extremal eigenvalues' distribution
%and then
%perform an asymptotic analysis of the associated kernels to prove
%convergence
%toward the Airy kernel. The later analysis is based on contours
%deformations.
%Contrary to the analysis performed by Baik et al. \cite{BBP-2005} and
%El Karoui
%\cite{elkaroui-07-TW} {\red[jam: et les filles, elles font quoi?
%faudrait
%peut-{\chr"C3}{\chr"AA}tre nous positionner par rapport {\chr"C3}{
%\chr"A0} elles ici]} and based on explicit
%deformed contours, our analysis relies on a more abstract argument
%where the
%existence of appropriate contours is obtained by mean of the maximum
%principle
%for subharmonic functions. The advantage of this strategy is to avoid
%the
%need for an explicit contour (potentially very complicated to obtain
%in our
%setting).
% Attention: l'argument ne se duplique pas tel quel. Par exemple le cas
% c < 0 reclame un developpement particulier.
% The second advantage is that the argument can be duplicated without
% almost any change to every regular edge.

In part (c), our proof of the asymptotic independence builds upon the
operator-theoretic approach developed by Bornemann \cite
{bornemann-2010} in the
context of the GUE. We actually show that a weaker mode of convergence
for the involved operators than the one required in \cite
{bornemann-2010} is sufficient to establish the asymptotic
independence; it has the advantage to be compatible with the previous
asymptotic analysis.

Part (d) also relies on an asymptotic analysis of the rescaled kernel.
It is
based on an appropriate representation of the Bessel kernel as a double complex
integral.

\begin{longlist}
\item[\textit{Organization of the paper}.] In Section~\ref
{secbackground}, we provide a precise description for the bulk and the
extremal eigenvalues and introduce the notion of regular edge. The
precise statement of part (a) of Theorem~\ref{thebigtheorem} is
provided in Theorem~\ref{thdef-phi-N} and proved. In Section~\ref
{secmain-results}, we state our results
concerning the fluctuations of the extremal eigenvalues and their asymptotic
independence. Parts (b), (c) and (d) of Theorem~\ref{thebigtheorem}
are, respectively, stated in Theorems~\ref{thfluctuations-TW},~\ref{thasymptindep}~and~\ref{thBessel}. We also recall there the definition of the
Tracy--Widom distribution and the hard edge
distribution described by mean of the Bessel kernel (Sections~\ref
{secdescription-TW}~and~\ref{secdescription-BE}). We close this
section with an asymptotic study of the condition number of large
correlated Wishart matrices, a
discussion on nonregular edges and spikes phenomena and provide some graphical
illustrations. Section~\ref{SectionTW} is devoted to the proof for
Theorem~\ref{thfluctuations-TW} (Tracy--Widom fluctuations).
Section~\ref{secasymptotic-independence} is devoted to the proof of
Theorem~\ref{thasymptindep} (asymptotic independence for
extremal eigenvalues). Finally, Section~\ref{Besselsection} is devoted
to the proof of Theorem~\ref{thBessel} (hard edge fluctuations).
\end{longlist}

%s2 #&#
\section{Bulk description, regularity and extremal eigenvalues}\label{secbackground}

In this section, we introduce the notion of \textit{regular soft edges}
(cf. Definition~\ref{defreg}) and \textit{extremal eigenvalues} (cf.
Theorem~\ref{thdef-phi-N}), the main properties of which are gathered
in Propositions~\ref{propproperties-regular-left} and~\ref
{propproperties-regular-right}. Theorem~\ref{thdef-phi-N} provides a
precise statement for Theorem~\ref{thebigtheorem}(a). %We finally
%gather their main properties in
%Their main properties are gathered in Proposition~\ref{gN->g} and
%Theorem~\ref{thdef-phi-N}, which provides a precise statement for
%Theorem~\ref{thebigtheorem}-(a).
Before this, we provide a precise description of the bulk, mainly based
on \cite{silverstein-choi-1995}.

% In this section we provide a precise description for the bulk, based
%on the work by Silverstein and Choi (see also
%\cite{mestre-2008a,mestre-2008} by Mestre); we introduce the notion of
%{\mathbf regular edge} and describe its consequences. Finally, Theorem
%\ref{thdef-phi-N} provides a precise statement for part (a) of the
%main theorem; its proof ends the section.

%s2.1 #&#
\subsection{Description of the limiting bulk}

In \cite{mar-and-pas-67}, Mar\v cenko and Pastur characterized the
Cauchy--Stieltjes transform\footnote{Note that our definition of the
Cauchy--Stieltjes transform differs
by a sign from the one in \cite{mar-and-pas-67} but will turn out to
be more convenient in the sequel.} of the limiting distribution
$\mu= \mu(\gamma,\nu)$ of the eigenvalues of $\mathbf M_N$ as $N\to
\infty$,
\[
\label{cauchy} m(z)=\int\frac{1}{z-\lambda} \mu(\d\lambda),\qquad z\in\C_+=\bigl\{z
\in\C\dvtx  \operatorname{Im}(z)>0\bigr\},
\]
as the unique solution $m \in\C_-=\{z\in\C\dvtx    \operatorname{Im}(z)<0\}$ of the
fixed-point equation
%
%e8 #&#
\begin{equation}
\label{pt-fixe} m= \biggl(z-\gamma\int\frac{\lambda}{1- m\lambda} \nu(\d\lambda )
\biggr)^{-1} \qquad\mbox{for any } z \in\C_+.
\end{equation}

Recall that by Assumption~\ref{assgauss}, $\gamma=\lim n/N\in
(0,+\infty)$, the probability measure $\nu$ is the limiting
eigenvalue distribution of $\bolds
\Sigma_N$, and its compact support is included in $(0,+\infty)$. In
particular, $\nu(\{0\})=0$.

In \cite{silverstein-choi-1995}, Silverstein and Choi showed that
%
%e9 #&#
\begin{equation}
\mu(\d x)=(1-\gamma)^+\delta_0+\rho(x)\,\d x,
\end{equation}
where $\rho$ is a nonnegative and continuous function on $(0,\infty
)$ which is analytic wherever it is positive.
Moreover, following a procedure already described by Mar\v cenko and Pastur,
they showed rigorously how to extract from the fixed point equation
above a
characterization of the support of $\mu$, and thus of $\rho(x)\,\d x$.
Specifically, the function $m(z)$ has
an explicit inverse on $m(\C_+)$ given by
%
%e10 #&#
\begin{equation}
\label{g} g(z)=\frac{1}{z}+\gamma\int\frac{\lambda}{1-z \lambda} \nu(\d
\lambda),
\end{equation}
and this inverse extends analytically to a neighborhood of $\C_- \cup
D$ where $D$
is the open subset of the real line
%
%e11 #&#
\begin{equation}
\label{def-D} D= \bigl\{x\in\R\dvtx  x\neq0, x^{-1}\notin\supp(\nu) \bigr\}.
\end{equation}
Except in the proof of Proposition~\ref{gN->g} below, we shall confine the
notation $g$ to the restriction of this function to $D$.
On any interval $I$ of $\R\setminus\supp(\mu)$, the function $m$
exists, is real
and is decreasing (as a Cauchy--Stieltjes transform).
Consequently, its inverse also exists and is
decreasing on $m(I)$. Silverstein and Choi showed
that $g$ is this inverse, and that $\R\setminus\supp(\mu)$
coincides with
the values of $g(x)$ where this function is decreasing on $D$:
% The following two propositions hold true for any probability measure $
%\nu$
% supported by $\R_+$ and any $\gamma> 0$:

%pr2.1 #&#
\begin{proposition}[(Silverstein and Choi \cite{silverstein-choi-1995})]\label{SC-support}
%Let $\mu$ be the probability measure whose Cauchy-Stieltjes transform
%is
%determined by~\eqref{pt-fixe}.
For any $x \in\R\setminus\supp(\mu)$, let $p = m(x)$. Then $p \in D$,
$x = g(p)$ and $g'(p) < 0$. Conversely, let $p \in D$ such that $g'(p)
< 0$.
Then $x = g(p) \in\R\setminus\supp(\mu)$ and $p = m(x)$.
\end{proposition}

\begin{figure}

\includegraphics{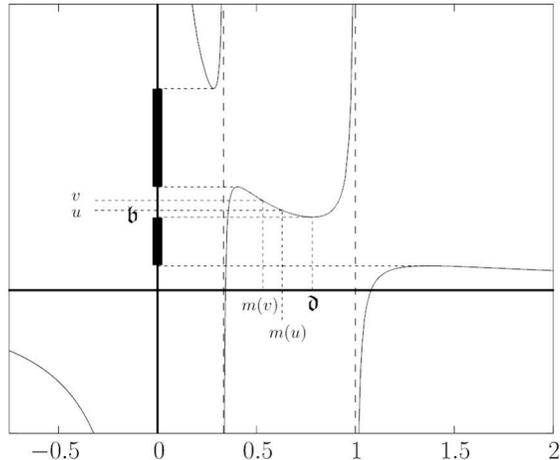}

\caption{Plot\vspace*{1pt} of $g\dvtx D\to\mathbb{R}$ for $\gamma=0.1$ and $\nu=
0.7\delta_1 + 0.3\delta_3$. In this case, $
D=(-\infty,0)\cup (0,\frac{1}3)\cup (\frac{1}3,1)\cup(1,\infty
) $.
%\mathbb{R}\setminus\{0,1/3,1\}$.
The two thick segments on the vertical axis represent $\supp(\mu)$.
The right edge $\frak b$ of the measure $\mu$ satisfies property \textup{(a)}
of Proposition \protect\ref{endpoints}.}
\label{fignospike}
\end{figure}

%re2.2 #&#
\begin{Remark}
 This proposition has the following practical
importance: in order to find $\supp(\mu)$, plot the function $g$ on
$D$; whenever $g$ is decreasing ($g'(x)<0$),
remove the corresponding points $g(x)$ from the vertical axis. What is
left is precisely $\supp(\mu)$.
\end{Remark}
As an example, a plot of the function $g$ is provided in Figure~\ref
{fignospike} along with $\supp(\mu)$ in the case where $\nu$ is the
weighted sum of
two Dirac measures and $\gamma<1$.

The soft edges of the bulk are described more precisely by the
next proposition.

%pr2.3 #&#
\begin{proposition}[(Silverstein and Choi \cite{silverstein-choi-1995})]\label{endpoints}
Any soft left edge $\frak a$ satisfies one
of the two following properties:
\begin{longlist}[(a)]
\item[(a)] There exists a unique $\frak c \in D$ such that
$\frak a = g(\frak c)$, $g'(\frak c) = 0$ and $g''(\frak c) < 0$.
\item[(b)] There exists a unique $\frak c \in\partial D$ such that
$(\frak c, \frak c + \varepsilon) \subset D$ for some
$\varepsilon> 0$ small enough, the function $g$ is decreasing on
$(\frak c, \frak c + \varepsilon)$, and
$\frak a = \lim_{x\downarrow\frak c} g(x)$. In this case, we write
$\frak a = g(\frak c)$.
\end{longlist}
Conversely, for any point $\frak c$ satisfying one of these properties,
$\frak a = g(\frak c)$ is a soft left edge.

Similar, any (soft) right edge $\frak b$ of the measure $\mu$
satisfies one of the
two following properties:
\begin{longlist}[(a)]
\item[(a)] There exists a unique $\frak d \in D$ such that
$\frak b = g(\frak d)$, $g'(\frak d) = 0$ and $g''(\frak d) >0$.
\item[(b)] There exists a unique $\frak d \in\partial D$ such that
$(\frak d-\varepsilon, \frak d) \subset D$ for some
$\varepsilon> 0$ small enough, the function $g$ is decreasing on
$(\frak d-\varepsilon, \frak d)$, and $\frak b = \lim_{x\uparrow
\frak d} g(x)$.
In this case, we write $\frak b = g(\frak d)$.
\end{longlist}
Conversely, for any point $\frak d$ satisfying one of these properties,
$\frak b = g(\frak d)$ is a right edge of the measure $\mu$.
\end{proposition}

Hence any soft edge of the bulk coincides with a unique
extremum $\frak c$ of the function $g$, and it reads $g(\frak c)$.
These extrema may or may not be attained
on $D$. In case they are, the second derivative of $g$ is never equal
to zero
there, and it has been proved in \cite{silverstein-choi-1995} that the
density vanishes like a square root at the associated edges. We shall
see later that the Tracy--Widom fluctuations appear in this case.
A~right edge $\frak b=g(\frak d)$ together with its preimage $\frak d$
are plotted in Figure~\ref{fignospike}.

%
%f1 #&#

The next proposition provides additional information on the bulk
that will be useful in the sequel. Its proof is in Appendix~\ref{anx-SB-gamma}.

%pr2.4 #&#
\begin{proposition}
\label{1111}
Let Assumption~\ref{assnu} hold true. Let $\frak a $ be the leftmost
edge of the bulk. The following facts hold true:
\begin{longlist}[(a)]
\item[(a)]
If $\gamma> 1$, then $\frak a > 0$. Moreover, the function $g(x)$
increases from zero to $\frak a$ then decreases from $\frak a$ to
$-\infty$
as $x$ increases from $-\infty$ to zero. In particular, %if $
%\gamma>1$, then
$\frak a$ is the unique maximum of $g$ on $(-\infty,0)$.
\item[(b)]
If $\gamma\leq1$, the function $g$
is negative and decreasing on $(-\infty, 0)$.
\item[(c)]
If $\gamma< 1$, then $\frak a > 0$. Moreover, if we set
$\eta=\inf\supp(\nu)>0$, then $\frak a = g(\frak c)$ is the
supremum of
$g$ on $(1/\eta, \infty)$. In addition, $g$ increases to $\frak a$ on
$(1/\eta, \frak c)$ whenever this interval is nonempty, then decreases
from $\frak a$ to zero on $(\frak c, \infty)$.
\end{longlist}
Let $\frak b = g(\frak d)$ be a right edge of the bulk. Then the following
facts hold true:
\begin{longlist}[(a)]
\item[(d)]
$[ \frak d, \infty) \not\subset D$.
\item[(e)]
Assume $\frak b$ is the rightmost edge of the bulk. For any
$\gamma\in(0,\infty)$, if we set $\xi=\sup\supp(\nu)<\infty$,
then $g$
decreases from infinity to $\frak b$ on $(0, \frak d)$ and increases on
$(\frak d, 1/\xi)$ if this last interval is not empty. In particular,
$\frak d$ is the unique extremum of $g$ on $(0,1/\xi)$.
\end{longlist}
\end{proposition}

Fact (a) shows that when $\gamma> 1$, the study of $g$ on $(-\infty, 0)$
allows us to locate the leftmost edge $\frak a$ and this edge only.
Facts (a) and (b) show that if $\gamma\leq1$, then it suffices to
study $g$
on $D \cap(0,\infty)$ to locate the edges of the bulk.
%Facts (a) and (b) show that except for
%this case, one needs to study $g$ on $D \cap(0,\infty)$ to locate the
%edges
%of the bulk.
In particular, if $\gamma< 1$, fact~(c) shows that the
location of $\frak a$ is provided by the study of $g$ on $(1/\eta,
\infty)$.
This is illustrated by Figure~\ref{fignospike}, where $\frak a$ is the
rightmost maximum of the function $g$. Fact~(d) shows that when
$\frak b = g(\frak d)$ is a right edge of the bulk, then $\frak d$ cannot
belong to the unbounded connected component of $D$ in $(0,\infty)$.
Finally, the behavior of $g$ described by~(e) is illustrated on
Figure~\ref{fignospike} by the plot of this function on the interval
$(0, 1/3)$.

%s2.2 #&#
\subsection{Regularity condition and its consequences}
So far, we have thoroughly described the edges of the limiting
eigenvalue distribution. Remember, however, that BBP phase transition
\cite{BBP-2005}
may occur regardless of the limiting spectral distribution (which is
always \v MP distribution in \cite{BBP-2005}). As we shall see later,
the notion of a \textit{regular endpoint} captures a joint condition on
the limiting spectral distribution $\mu$ and on the convergence $\nu
_N\to\nu$, which will guarantee Tracy--Widom fluctuations; cf.
Theorem~\ref{thfluctuations-TW}.
% We now introduce the notion of regular endpoints that plays a key
%role in this paper:

%
%de2.5 #&#
\begin{definition}[(Regular edge)]\label{defreg}
Recall that the $\lambda_i$'s are the eigenvalues of matrix ${\bolds
\Sigma}_N$; a soft edge $\frak a=g(\frak c)$ is \emph{regular} if
%
%e12 #&#
\begin{equation}
\label{eqdef-regular} \liminf_{N\rightarrow\infty}\min_{j=1}^n
\bigl\llvert \frak c-\lambda _j^{-1}\bigr\rrvert >0.
\end{equation}
In particular, $\frak c\in D$.
\end{definition}

%{\red[jam: j'ai rassembl{\chr"C3}{\chr"A9} les 2 remarques et la
%proposition pr{\chr"C3}{\chr"A9}c{\chr"C3}{\chr"A9}dentes dans une
%seule remarque et 3 items]}

%\begin{remark}
%\label{in-D}
%If $\frak a= g(\frak c)$ is a regular soft edge, then the weak
%convergence $\nu_N\rightarrow\nu$ stated in Assumption~\ref{assgauss}
%rules out
%the options labelled (b) in Proposition~\ref{endpoints}.
%\end{remark}
%
%\begin{remark}
%\label{in-D2}
%If $\frak a$ is an endpoint satisfying one of the options labelled (a)
%in Proposition~\ref{endpoints}, and if, furthermore, the distance
%$\dist( \lambda_i, \supp(\nu) )$ satisfies
%$$
%\max_{1\leq j\leq n} \dist(\lambda_j, \supp(\nu) )\mathop{
%\longrightarrow}_{N\to\infty} 0,
%$$
%then $\frak a$ is a regular endpoint of $\supp(\mu)$. However, this
%last condition is not necessary. Further comments will be made in
%Remark~\ref{outliers} below.
%\end{remark}
%
%\begin{proposition}
%\label{propa-reg}
%If $\gamma> 1$, then the leftmost edge is regular.
%\end{proposition}
%\begin{proof}
%Writing the leftmost edge as $g(\frak c)$, Proposition~
%\ref{SB-gamma}--(a)
%shows that $\frak c < 0$, which immediatly implies
%\eqref{eqdef-regular}.
%\end{proof}
%
%{\red[jam: tiens voil{\chr"C3}{\chr"A0} ce que {\chr"C3}{\chr"A7}a
%donnerait..]}

%re2.6 #&#
\begin{Remark}
\label{remark-with-3-items}
The following facts will illustrate the range of the definition:
\begin{longlist}[(a)]
\item[(a)] If $\frak a= g(\frak c)$ is a regular soft edge, then the weak
convergence $\nu_N\rightarrow\nu$ stated in Assumption~\ref{assnu}
rules out
the options labeled (b) in Proposition~\ref{endpoints}.
\item[(b)]
If $\frak a$ is an endpoint satisfying one of the options labeled (a)
in Proposition~\ref{endpoints}, and if, furthermore, the distance
$\dist( \lambda_i, \supp(\nu) )$ satisfies
\[
\max_{1\leq j\leq n} \dist\bigl(\lambda_j, \supp(\nu)
\bigr) \mathop{\longrightarrow}_{N\to\infty} 0, %
\]
then $\frak a$ is a regular endpoint of $\supp(\mu)$. However, this
last condition is not necessary. Further comments will be made in
Section~\ref{outliers} below.

\item[(c)] If $\gamma> 1$, then the leftmost edge is regular. [For a proof of
this fact, simply write
the leftmost edge as $g(\frak c)$, then Proposition~\ref{1111}(a)
shows that $\frak c < 0$, which immediately implies (\ref{eqdef-regular}).]
\end{longlist}
\end{Remark}

Let $\gamma_N=n/N$, and consider now the probability measure $\mu
(\gamma_N, \nu_N)$, which is the
unique solution of the fixed point equation (\ref{pt-fixe}) associated
with the data $\gamma_N,\nu_N$. It is a finite-$N$ deterministic
equivalent of the spectral measure of $\mathbf M_N$. Associated to $\mu
(\gamma_N, \nu_N)$ is the function
%
%e13 #&#
\begin{equation}
\label{gN} g_N(z) = \frac{1}z + \gamma_N
\int\frac{\lambda}{1 - z\lambda} \nu _N(\d\lambda) = \frac{1}z +
\frac{1}N \sum_{j=1}^n
\frac{\lambda_j}{1 - z\lambda_j};
\end{equation}
cf. (\ref{g}). Similarly to $\mu(\gamma,\nu)$,
the measure $\mu(\gamma_N, \nu_N)$ has a density on $(0,\infty)$
and its support can also be characterized with the help of Proposition
\ref{SC-support}
(by simply replacing $g$ by $g_N$).
%whose
%support can be characterized similarly to $\supp(\mu(\gamma,\nu))$
%after
%replacing $g$ with $g_N$.
We furthermore have the following proposition:

%pr2.7 #&#
\begin{proposition}
\label{gN->g}
Let Assumption~\ref{assnu} hold true.
Let $ g(\frak c)$ be a regular soft edge. Then for $N$ large enough:
\begin{longlist}[(a)]
\item[(a)] $g_N$ is analytic in a complex neighborhood of $\frak
c$ which
is independent of~$N$;
\item[(b)] $g_N$ converges to $g$ uniformly\vspace*{1pt} on the compact sets
of this
neighborhood, and so does its $k$th order derivative
$g_N^{(k)}$ to $g^{(k)}$, for any $k\geq1$;%of these derivatives.
\item[(c)] There exists a sequence of real numbers
$\frak c_N$, unique up to a finite number of terms,
%{\red[jam: {\chr"C3}{\chr"A7}a veut dire quoi exactement? comme $N$
%est sens{\chr"C3}{\chr"A9} {\chr"C3}{\chr"AA}tre grand, on peut virer
%"unique up..", non?]}
such that
$\frak c_N \to\frak c$, $g'_N(\frak c_N) = 0$ and
$g^{(k)}_N(\frak c_N) \to g^{(k)}(\frak c)$ as $N\to\infty$ for any~$k$.% of these derivatives.
\end{longlist}
\end{proposition}

This proposition shows in particular that when a soft edge $ g(\frak
c)$ is regular, there is
a sequence $g_N(\frak c_N)$ of endpoints of
$\supp(\mu(\gamma_N,\nu_N))$ that converge to $g(\frak c)$, and
$\frak c_N$
satisfies
%
%e14 #&#
\begin{equation}
\label{regul-gN} \liminf_{N\rightarrow\infty} \min_{j=1}^n
\bigl\llvert \frak c_N - \lambda _j^{-1}
\bigr\rrvert > 0.
\end{equation}

\begin{pf*}{Proof of Proposition~\ref{gN->g}}
Set
$\eta= \min( \llvert  \frak c\rrvert   / 2, \liminf_N \min_j \llvert   \lambda_j^{-1} -
\frak c \rrvert  )$,
and let $B = B(\frak c, \eta/2)$ be the open ball with center $\frak
c$ and
radius $\eta/2$. Since
\[
\frac{\lambda_j}{\llvert  1 - z\lambda_j\rrvert  } = \frac{1}{\llvert  \lambda_j^{-1} - z\rrvert  } \leq \frac{1}{\llvert  \lambda_j^{-1} - \frak c\rrvert   - \llvert   z - \frak c\rrvert  } \leq
\frac{3}{\eta}
\]
for $z\in B$ and for all $N$ large, the functions $g_N$ are analytic and
uniformly bounded on $B$ for all $N$ large. This establishes (a) in
particular.
Moreover, this yields that the family of analytic functions $g_N$ is uniformly
bounded on $B$. Thus by Montel's theorem, the family $g_N$ is normal.
It follows from the convergences $\gamma_N\to\gamma$ and $\nu_N\to
\nu$ provided by Assumption~\ref{assgauss} that $g_N$ converges
pointwise to
$g$ on $B$. Consequently, $g_N$ converges to $g$ uniformly on the
compact subsets of $B$, and the same is true for the convergence\vspace*{1pt}
of the $g^{(k)}_N$ to $g^{(k)}$ by \cite{book-rudin-real-complex}, Theorem~10.28.
Turning to (c), notice that $\frak c$ is a zero of $g'$ by the
regularity \mbox{assumption}; see Remark~\ref{remark-with-3-items}(a). Since
$g'_N$ converges
to $g'$ uniformly on the compact sets of $B$ and $g'$ is analytic
there, Hurwitz's theorem shows
that $g'_N$ has a zero $\frak c_N$ that converges to the zero $\frak c$
of $g'$ and that this zero is unique provided $N$ is large enough.
Moreover this zero is real since $g'_N(\overline{z})=\overline{g'_N(z)}$.
Write $\llvert  g^{(k)}_N(\frak c_N) - g^{(k)}(\frak c)\rrvert   \leq
\llvert   g^{(k)}_N(\frak c_N) - g^{(k)}(\frak c_N) \rrvert   +
\llvert   g^{(k)}(\frak c_N) - g^{(k)}(\frak c) \rrvert  $.
Since for any $k$, $g^{(k)}_N$ converge uniformly to $g^k$ on the compact
subsets of $B$, the first term at the right-hand side vanishes as
$N\rightarrow\infty$. The second
term vanishes as $N\rightarrow\infty$ by the continuity of $g^{(k)}$.
This establishes (c).
\end{pf*}

%s2.3 #&#
\subsection{Extremal eigenvalues and their convergence}

% {\red To make more precise
%
% \begin{theorem}[Extremal eigenvalues] \label{thdef-phi-N}
% \begin{longlist}[(2)] % \item[(a)] Assume that $\frak a$ is a soft edge which is
%regular in the sense of Definition~\ref{defreg}. Then, there exists a
%deterministic sequence $(\phi_N)_N$ of positive numbers such that,
%almost surely,
% \begin{equation}
% \label{leftmostev}
% \lim_{N\rightarrow\infty}x_{\phi_N}=\frak a,\qquad\limsup_{N
%\rightarrow\infty}x_{\phi_N-1}<\frak a.
% \end{equation}
% If $\lim_N\phi_N=1$, which happens if and only if $\frak a$ is the
%leftmost endpoint of $\supp(\mu)$ and
% \[
% \limsup_{N\rightarrow\infty}\frak c_N<0,
% \]
% then the statement involving $\phi_{N-1}$ in \eqref{leftmostev} is
%considered as empty.
% \item[(b)] Assume that $\frak b$ is a regular right edge of $
%\supp(\mu)$ in the sense of Definition~\ref{defreg} and let $(\frak
%c_N)_N$ be the sequence associated with $\frak b$ as in Proposition
%\ref{detequprop}.
% There exists a sequence $(\phi_N)_N$ where $\phi_N\in\{1,\ldots,N\}$
%such that, almost surely,
% \begin{equation}
% \label{rightmostev}
% \lim_{N\rightarrow\infty}x_{\phi_N}=\frak b,\qquad\liminf_{N
%\rightarrow\infty}x_{\phi_N+1}>\frak b.
% \end{equation}
% If $\lim_N\phi_N/N=1$, which happens if and only if $\frak b$ is the
%rightmost edge of $\supp(\mu)$ and
% \begin{equation}
% \label{cNisolatedmax}
% \limsup_{N\rightarrow\infty}\frak c_N\lambda_n<1,
% \end{equation}
% then the statement involving $\phi_{N+1}$ in \eqref{rightmostev} is
%considered as empty.
%
% \end{longlist}
% \end{theorem}
% }

Our purpose is now to locate the eigenvalues of $\mathbf M_N$ that
converge to a prescribed edge,\vspace*{1pt} or equivalently
those of $\widetilde{\mathbf M}_N$ (denoted by $\tilde x_1 \leq\cdots
\leq
\tilde x_n$).
%, by a deterministic indexation: $x_{\phi(N)} \to\frak a$.
%To that end, we use the following result due to Bai and Silverstein
%concerning the localisation of the eigenvalues of $\widetilde{\bv
%M}_N$.
The idea is the following: given an interval $(u,v)$ outside $\operatorname{Supp}(\mu(\gamma_N,\nu_N))$, its\vspace*{2pt} preimage $(m(v),m(u))$ by $g$ then
lies in between two groups of $\lambda_j^{-1}$'s, provided $N$ is
sufficiently large. Thus there is a unique integer\vspace*{1pt} $\phi(N)$ for which
$\lambda_{\phi(N)+1}^{-1} <m(v)<m(u)<\lambda_{\phi(N)}^{-1}$. This
$\phi(N)$ defines the deterministic index for which $x_{\phi(N)}$
converges a.s. toward the prescribed edge. Figure~\ref{fignospike}
illustrates this phenomenon. The following proposition formalizes this.

%re2.8 #&#
\begin{Remark}[(Convention)]\label{remconvention}
In the remaining, we shall systematically use the notational convention
$\lambda_0 = \tilde x_0 = 0$ and $\lambda_{n+1} = \tilde x_{n+1} =
\infty$.
\end{Remark}

%
%pr2.9 #&#
\begin{proposition}[(Bai and Silverstein \cite{bai-silverstein-98,silverstein-exact-separation-99})]
\label{propsep}
Let Assumptions~\ref{assgauss}~and~\ref{assnu} hold true. Assume that
$[u,v]$ with
$u > 0$ lies in an open interval outside $\supp(\mu(\gamma_N, \nu
_N))$ for $N$ large enough, and recall definition (\ref{pt-fixe}) of
the fixed-point solution $m$. Then the following facts hold true:
\begin{longlist}[(a)]
\item[(a)]
If $\gamma> 1$, then $\tilde x_{n-N+1} \to\frak a$ almost surely
as $N\to\infty$, where $\frak a > 0$ is the leftmost edge of the bulk.
\item[(b)]
In the following cases: \textup{(i)} $\gamma\leq1$ or \textup{(ii)} $\gamma> 1$ and
$[u,v] \not\subset[0, \frak a]$, it holds that $m(v) > 0$. Let $\phi
(N)$ be
the integer defined as
%
%e15 #&#
\begin{equation}
\label{sep-lambda} \lambda_{\phi(N)+1} > m(v)^{-1} \quad\mbox{and}\quad
\quad \lambda_{\phi(N)} < m(u)^{-1}.
\end{equation}
Then
%
%e16 #&#
\begin{equation}
\label{eqsep} \mathbb{P} (\tilde x_{\phi(N)+1} > v, \tilde
x_{\phi(N)} < u \mbox{ for all large } N ) = 1.
\end{equation}
\end{longlist}
\end{proposition}

%re2.10 #&#
\begin{Remark}\label{infinite-array}
Bianchi et al. \cite{bianchi-et-al-2011} established the
result for matrices ${\mathbf X}_N$ taken from a doubly infinite array of i.i.d.
random variables with finite fourth moment. If the entries are
Gaussian, one can relax the doubly infinite array assumption and establish
Proposition~\ref{propsep} by using the completely different tools
of~\cite{loubaton-vallet-11}.
\end{Remark}
%
%To be precise, this result was established in
%\cite{silverstein-exact-separation-99} in the case where
%the elements of $\bv X_N$ are taken from a doubly infinite array. When
%these
%elements are Gaussian (note that
%\cite{bai-silverstein-98,silverstein-exact-separation-99} are not
%restricted to this
%case), one can relax the doubly infinite array assumption and establish
%Proposition~\ref{propsep} by using the completely different tools
%of~\cite{loubaton-vallet-11}.
%\end{remark}

We are now in a position to properly state and prove part~(a)
of Theorem~\ref{thebigtheorem}.
%We can now state the result on the extremal eigenvalues of
%$\widetilde{\bv M}_N$.

%
%th2 #&#
\begin{theorem}[(Extremal eigenvalues)]\label{thdef-phi-N}
Let Assumptions~\ref{assgauss}~and~\ref{assnu} hold true.\footnote
{In view of Remark~\ref{infinite-array}, one can relax the Gaussianity
assumption in Theorem~\ref{thdef-phi-N}
and replace it by the fact that ${\mathbf X}_N$'s entries are extracted
from a doubly infinite array of i.i.d. random variables.}
\begin{longlist}[(a)]
\item[(a)]
If $\gamma> 1$ and $\frak a$ is the leftmost edge of the bulk, then
set $\varphi(N) = n-N+1$.
Otherwise, let $\frak a = g(\frak c)$ be a regular soft left edge, and let
$\varphi(N) = \min\{ j\dvtx    \lambda_j^{-1} < \frak c \}$. Then, almost
surely,
\[
\lim_{N\to\infty} \tilde x_{\varphi(N)} = \frak a \quad\mbox
{and}\quad \liminf_{N\rightarrow\infty} ( \frak a - \tilde x_{\varphi(N) - 1})
> 0.
\]
\item[(b)]
Let $\frak b = g(\frak d)$ be a regular right edge, and let
$\phi(N) = \max\{ j\dvtx    \lambda_j^{-1} > \frak d \}$. Then, almost
surely
\[
\lim_{N\to\infty} \tilde x_{\phi(N)} = \frak b \quad\mbox
{and}\quad \liminf_{N\to\infty} ( \tilde x_{\phi(N) + 1} - \frak b)
> 0.
\]
\end{longlist}
Eigenvalues $\tilde x_{\varphi(N)}$ and $\tilde x_{\phi(N)}$ are
called extremal eigenvalues.
\end{theorem}

\begin{pf}
We shall only prove the result for a right edge $\frak b$.
By Proposition~\ref{gN->g}, we can choose a compact neighborhood $B$
of $\frak d$ such that $g_N$ and $g'_N$ uniformly converge to $g$ and
$g'$. Let $p, q, r, s$ be real numbers such that
$p < q < r < s < \frak d$, $[p,s] \subset B$, and
$g'(x) < 0$ for $x \in[p,s]$. This last condition is made possible by
the fact
that $\frak b$ is a right edge of $\supp(\mu)$; cf. Figure~\ref{fignospike}.
Let $u= g(r)$ and $v = g(q)$. Since $g_N$ and $g'_N$
converge uniformly to $g$ and $g'$, respectively, on $[p,s]$, it holds that
$g'_N(x) < 0$ on $[p,s]$, and $[u,v] \subset[ g_N(s), g_N(p)]$ for all $N$
large.
Proposition~\ref{SC-support} applied to $\mu(\gamma_N,\nu_N)$ shows
then that
$[u,v]$ lies in an open set outside $\supp(\mu(\gamma_N,\nu_N))$
for all
$N$ sufficiently large.

Now the integer $\phi(N)$ defined in the statement is characterized by
the inequalities
\[
\lambda_{\phi(N)+1}^{-1} < \frak d < \lambda_{\phi(N)}^{-1}. %
\]
Since
no $\lambda_j^{-1}$'s belong to $B$ for $N$ large enough, we can
equivalently write
\[
\lambda_{\phi(N)+1}^{-1} < q = m(v) < r = m(u) < \lambda
_{\phi(N)}^{-1} %
\]
which is~(\ref{sep-lambda}). By Proposition~\ref{propsep}, we
get~(\ref{eqsep}).

Since $v > \frak b$, we have
$\liminf_N ( \tilde x_{\phi(N) + 1} - \frak b) > 0$ with probability one.
Moreover, we know that a.s., the number of $\tilde x_i$ in
$[\frak b-\varepsilon, \frak b]$ is nonzero for any $\varepsilon> 0$ and
for all large $N$.
Making $r \uparrow\frak d$, we get $u=g(r) \downarrow\frak b$.
Since $\tilde x_{\phi(N)} < u$ a.s.~for all large $N$, we get that
$\tilde x_{\phi(N)} \to\frak b$ a.s.~when $N \to\infty$.
\end{pf}

% Since the correspondance between the random eigenvalues $x_{N-
%\min(n,N)+1}\leq\cdots\leq x_{N}$ of $\bv M_N$ and the random
%eigenvalues $\tilde x_{n-\min(n,N)+1}\leq\cdots\leq\tilde x_n$ of $
%\widetilde{\bv M}_N$ reads
% \begin{equation}
% x_j=\tilde x_{j+n-N},\qquad N-\min(n,N)+1\leq j\leq N,
% \end{equation}
% Theorem~\ref{thdef-phi-N} easily follows Corollary~\ref{def-phi}.

%s2.4 #&#
\subsection{Summary of the properties of regular edges}

For the reader's convenience and constant use in the sequel, we gather
in the two following propositions some of the most important properties
of regular edges introduced above. Recall the convention in Remark~\ref
{remconvention}.

%
%pr2.11 #&#
\begin{proposition}[(Left regular soft edges)]\label
{propproperties-regular-left} Let Assumption~\ref{assnu} hold true.
Let $\frak a$ be a left edge.
\begin{longlist}[(a)]
\item[(a)] Consider first the case where $\frak a$ is the
leftmost edge:
\begin{longlist}[--]
\item[--] If $\gamma>1$, then $\frak a=g(\frak c)>0$ with $\frak c<0$, and
$\frak a$ is a regular soft edge.

\item[--] If $\gamma<1$, then $\frak a=g(\frak c)>0$ with $\frak c>0$;
$\frak a$ is a soft edge, but its regularity is a priori not granted.
% \item If $\gamma=1$ then $\frak a=0$ is a hard edge hence not regular.
\end{longlist}
\item[(b)] Assume now that $\frak a$ is a regular left soft
edge. Then
\begin{eqnarray*}
\frak a &=& g(\frak c)
\end{eqnarray*}
with
\[
\cases{ \displaystyle g'(
\frak c) = 0,
\cr
\displaystyle g''(\frak c) <0} \quad
\mbox{and}\quad \cases{ \displaystyle\frak c<0, &\quad if $\frak a$ is the
leftmost edge and $\gamma> 1$,
\cr
\displaystyle\frak c>0, &\quad otherwise.}
\]
For $N$ large enough, there exists a unique sequence $\frak c_N$ such
that $g_N'(\frak c_N)=0$ and
\[
\frak c_N\mathop{\longrightarrow}_{N\to\infty} \frak c,\qquad
g_N^{(k)}(\frak c_N) \mathop{
\longrightarrow}_{N\to\infty} g^{(k)}(\frak c)\qquad\mbox {for any } k
\geq0, %
\]
where by $g^{(0)}_N,g^{(0)}$ we mean $g_N,g$. Finally, there exists a
deterministic sequence $(\varphi(N))$ such that almost surely,
\[
\lim_{N\to\infty} \tilde x_{\varphi(N)} = \frak a \quad\mbox
{and}\quad \liminf_{N\rightarrow\infty} ( \frak a - \tilde x_{\varphi(N) - 1})
> 0. %
\]
\end{longlist}
\end{proposition}

%
%pr2.12 #&#
\begin{proposition}[(Right regular soft edges)]\label{propproperties-regular-right} Let Assumption~\ref{assnu} hold true, and
assume that $\frak b$ is a regular right soft edge. Then
\[
\frak b = g(\frak d)\qquad\mbox{with } \cases{ \displaystyle g'(
\frak d) =0,
\cr
\displaystyle g''(\frak d) > 0}
\mbox{ and }\frak d>0. %
\]
For $N$ large enough, there exists a unique sequence $\frak d_N$ such
that $g_N'(\frak d_N)=0$ and %(by $g^{(0)}$ we mean $g$)
\[
\frak d_N \mathop{\longrightarrow}_{N\to\infty} \frak d,\qquad
g_N^{(k)}(\frak d_N) \mathop{
\longrightarrow}_{N\to\infty} g^{(k)}(\frak d)\qquad\mbox {for any } k
\geq0. %
\]
Finally, there exists a deterministic sequence $(\phi(N))$ such that
almost surely,
\[
\lim_{N\to\infty} \tilde x_{\phi(N)} = \frak b \quad\mbox
{and}\quad \liminf_{N\rightarrow\infty} ( \tilde x_{\phi(N) + 1}- \frak b)
> 0. %
\]
\end{proposition}

%s3 #&#
\section{Fluctuations around the edges}\label{secmain-results}
In this section, we state the main results of the paper, namely the
fluctuations of the extremal eigenvalues and their asymptotic independence.
Parts (b), (c) and (d) of Theorem~\ref{thebigtheorem} are respectively
formalized in Theorem~\ref{thfluctuations-TW} (Section~\ref
{secdescription-TW}),
Theorem~\ref{thasymptindep} (Section~\ref{secAIresults}) and Theorem
\ref{thBessel} (Section~\ref{secdescription-BE}). We also provide a
discussion on nonregular edges
and spikes phenomena with graphical illustrations.

As an application, we obtain in Section~\ref{conditionnumber} new
results for the asymptotic behavior of the condition number of complex
correlated Wishart matrices.

%In this section, we state our results concerning the fluctuations of
%the
%extremal eigenvalues and their asymptotic independence. The precise
%versions of
%the informal Theorem (b), (c) and (d) from the introduction are Theorem
%\ref{thfluctuations-TW}, Theorem~\ref{thasymptindep} and Theorem
%\ref{th%Bessel} respectively. We also provide a discussion on the
%nonregular edges
%and spikes phenomena with graphical illustrations.
%
%s3.1 #&#
\subsection{Tracy--Widom fluctuations at the regular soft edges}\label{secdescription-TW}
We first introduce the Tracy--Widom distribution.
%Our first result states that at every regular soft edge the
%fluctuations for associated extremal eigenvalue follow the Tracy-Widom
%distribution, that we introduce now.
The Airy function $\Aii$ is the unique solution of the differential
equation $ \Ai''(x)=x\Ai(x)$ which satisfies the asymptotic behavior
\[
\Ai(x)=\frac{1}{2\sqrt\pi x^{1/4}}e^{-\sklfrac{2}{3}x^{3/2}} \bigl(1+o(1) \bigr),\qquad x
\rightarrow+\infty.
\]
With a slight abuse of notation, denote by $\K_\Aii$ the integral
operator associated with the Airy kernel
%
%e17 #&#
\begin{equation}
\label{Kai} \K_{\Aii}(x,y)=\frac{\Ai(x)\Ai'(y)-\Ai(y)\Ai'(x)}{x-y}.
\end{equation}
A real-valued random variable $X$ is said to have Tracy--Widom distribution~if
\[
\mathbb{P}( X\leq s ) = \det (I-\K_{\Aii} )_{L^2(s,\infty
)},\qquad s\in
\R, %
\]
where the right-hand side stands for the Fredholm determinant of the
restriction to $L^2(s,\infty)$ of the operator $\K_{\Aii}$ (see also
Section~\ref{determinantalprocess}). Tracy and Widom \cite
{tracy-widom-airy-94} established the famous representation
\[
\det (I-\K_{\Aii} )_{L^2(s,\infty)}=\exp \biggl({-\int
_{s}^\infty}(x-s)q(x)^2\,\d x \biggr),
\]
where $q$ is the Hastings--McLeod solution of the Painlev\'e II
equation, namely the unique solution of $q''(x)=2q(x)^3+xq(x)$ with
boundary condition $q(x)\sim\Ai(x)$ as $ x\rightarrow\infty$.

%Before stating Tracy-Widom fluctuations for the soft edges, we first
%gather all the required informations. Let Assumption~\ref{assgauss}
%hold true and recall definition~\ref{defreg} of regular edges. Given a
%left regular soft edge $\frak a$, there exists an associated extremal
%eigenvalue $\tilde x_{\varphi(N)}$ by Theorem~\ref{thdef-phi-N} such
%that $\lim_{N\to\infty} \tilde x_{\varphi(N)} = \frak a$. Let $\frak
%c$ be such that $\frak a = g(\frak c)$, then there exists by
%Proposition~\ref{gN->g} a sequence $\frak c_N$ which converge to $
%\frak c$ such that $$g_N^{(k)}(\frak c_N)\mathop{\longrightarrow}_{N\to
%\infty} g^{(k)}(\frak c),$$ where $g_N$ has been introduced in
%\eqref{gN}. In short:
%\begin{equation}\label{leftsoftedge}
%{\frak a} \textrm{left regular soft edge} \quad\Rightarrow\quad
%\left\{
%\begin{array}{lc}
%\exists  \frak c, & \frak a = g(\frak c)\\
%\exists  \tilde x_{\varphi(N)}\in\textrm{eig}(\widetilde{\bv M}_N),&
%\tilde x_{\varphi(N)}\to\frak a\\
%\exists  \frak c_N \to\frak c, & g^{(k)}_N(\frak c_N) \to g^{(k)}(
%\frak c)
%\end{array}
%\right.,
%\end{equation}
%where $\textrm{eig}(\widetilde{\bv M}_N)$ stands for the spectrum of $
%\widetilde{\bv M}_N$. Similarly,
%\begin{equation}\label{rightsoftedge}
%{\frak b} \textrm{right regular soft edge} \quad\Rightarrow\quad
%\left\{
%\begin{array}{lc}
%\exists  \frak d, & \frak b = g(\frak d)\\
%\exists  \tilde x_{\phi(N)}\in\textrm{eig}(\widetilde{\bv M}_N),&
%\tilde x_{\phi(N)}\to\frak b\\
%\exists  \frak d_N \to\frak d, & g^{(k)}_N(\frak d_N) \to g^{(k)}(
%\frak d)
%\end{array}
%\right..
%\end{equation}
We are now in position to state our result concerning the Tracy--Widom
fluctuations.
Recall that $g_N$ has been introduced in (\ref{gN}).

%th3 #&#
\begin{theorem} \label{thfluctuations-TW}
Let Assumptions~\ref{assgauss}~and~\ref{assnu} hold true.
\begin{longlist}[(a)]
\item[(a)] Let $\frak a$ be a left regular soft edge, and
$\tilde x_{\varphi(N)}$ and $(\frak c_N)_N$ be as in Proposition~\ref
{propproperties-regular-left}. Set
%
%Consider a left soft edge $\frak a$ and assume it is regular
%in the sense of Definition~\ref{defreg}. Let $\tilde x_{\varphi(N)}$
%be the
%extremal eigenvalue of $\widetilde{\bv M}_N$ associated to $\frak a$
%as in the
%statement of Theorem~\ref{thdef-phi-N}--{\rm(a)}, and $(\frak c_N)_N$
%be the
%sequence associated with $\frak a$ coming from Proposition~
%\ref{gN->g}. If we
%set
%
\[
\frak a_N=g_N(\frak c_N),\qquad
\sigma_N= \biggl(\frac{2}{-g_N''(\frak
c_N)} \biggr)^{1/3}.
\]
Then, for every $s\in\R$,
%
%e18 #&#
\begin{equation}
\label{TWleft} \lim_{N\rightarrow\infty} \p \bigl( N^{2/3}
\sigma_N (\frak a_N-\tilde x_{\varphi(N)} )\leq s
\bigr) =\det (I-\K_{\Aii} )_{L^2(s,\infty)}.
\end{equation}

\item[(b)] Let $\frak b$ be a right regular soft edge, and
$\tilde x_{\phi(N)}$ and $(\frak d_N)_N$ be as in Proposition~\ref
{propproperties-regular-right}. Set
%
%be the
%extremal eigenvalue of $\widetilde{\bv M}_N$ associated to $\frak b$
%as in
%the statement of Theorem~\ref{thdef-phi-N}--{\rm(b)}, and
%$(\frak c_N)_N$ be the sequence associated with $\frak b$ coming from
%Proposition~\ref{gN->g}. If we set
%
\[
\frak b_N=g_N(\frak d_N),\qquad
\delta_N= \biggl(\frac{2}{g_N''(\frak
d_N)} \biggr)^{1/3}.
\]
Then, for every $s\in\R$,
%
%e19 #&#
\begin{equation}
\label{TWright} \lim_{N\rightarrow\infty} \p \bigl( N^{2/3}
\delta_N (\tilde x_{\phi(N)}-\frak b_N )\leq s
\bigr) =\det (I-\K_{\Aii} )_{L^2(s,\infty)}.
\end{equation}
\end{longlist}
\end{theorem}

The proof is deferred to Section~\ref{SectionTW}, and an outline is
provided in Section~\ref{outlineproof}.

\subsubsection*{Connexion with El Karoui's result}

Let us first comment the last theorem in the light of El Karoui's
result \cite{elkaroui-07-TW}; see also Onatski's work \cite
{onatski-2008-TW}. If we assume
%
%e20 #&#
\begin{equation}
\label{EKregularity} \liminf_{N\rightarrow\infty}\frak d_N
\lambda_n<1,
\end{equation}
then, as a consequence of the analysis provided in Section~\ref
{secbackground}, the sequence $(\frak d_N)_N$ is associated with the
rightmost\vspace*{1pt} edge $\frak b$, and the\vspace*{1pt} associated extremal eigenvalue has to
be the largest eigenvalue of $\widetilde{\mathbf M}_N$ (or
equivalently of $\mathbf M_N$). Moreover, (\ref{EKregularity}) implies
that $\frak b$ is regular, so that Theorem~\ref{thfluctuations-TW}
applies. This is the result of El Karoui, announced in the
\hyperref[sec1]{Introduction}, which he actually proves in a more general setting.

Indeed, in \cite{elkaroui-07-TW} the weak convergence of $\nu_N$
toward some limiting probability distribution and the convergence of
$n/N$ to some limit were not assumed; it is only assumed that $n/N$
stays in a bounded set of $(0,1]$ (actually of $(0,+\infty)$ after~\cite{onatski-2008-TW}) together with (\ref{EKregularity}). Let us
mention that only under these assumptions, by compactness, one can
always extract converging subsequences for $\nu_N$ and $n/N$ so that
our result applies along a subsequence.

Notice also that condition (\ref{EKregularity}) is stronger than our
regularity condition, since $\frak b$ can be regular with $\liminf_{N}\frak d_N\lambda_n>1$. In\vspace*{2pt} this case, the extremal eigenvalue
associated with the rightmost edge is no longer the largest eigenvalue of
$\widetilde{\mathbf M}_N$; this entails the presence of outliers, as
we shall explain in the next paragraph. Our result then states that the
largest eigenvalue which actually converges to the rightmost edge
$\frak b$ fluctuates for large $N$, according to the Tracy--Widom law.

%{\red[jam: il me semble n{\chr"C3}{\chr"A9}anmoins que pour la premi{
%\chr"C3}{\chr"A8}re valeur propre "r{\chr"C3}{\chr"A9}guli{\chr"C3}{
%\chr"A8}re", disons $\lambda_{n'}$, on aurait alors
%$$
%\liminf_{} \frak d_n \lambda_{n'}<1
%$$
%non?]}

\subsubsection*{Nonregular edges and spikes phenomena}
\label{outliers}
\begin{figure}

\includegraphics{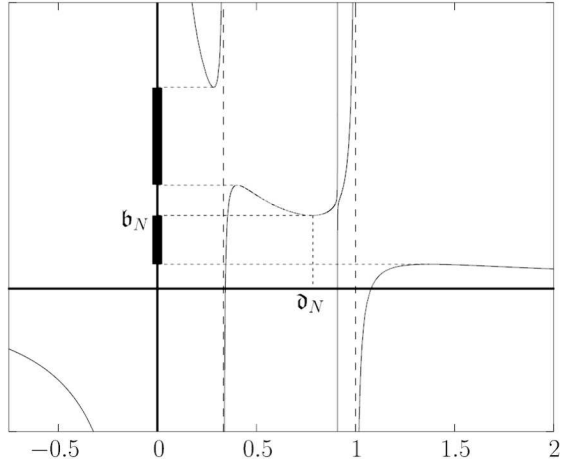}

\caption{Plot of $g_N(x)$ for $n = 300$, $\gamma_N = 0.1$ and
$\nu_N = \frac{1}{300} \delta_{1.1} + \frac{209}{300}\delta_1
+ \frac{90}{300}\delta_3$. The spike $\zeta=1.1$ does not produce an outlier.
The asymptote at $1/\zeta$ is not shown for better visibility.}\label{fignoout}
\end{figure}

In Remark~\ref{remark-with-3-items}(b), we explained that when a soft
edge reads
$\frak b = g(\frak d)$ with $\frak d \notin\partial D$, and when the
Hausdorff distance between
$\supp(\nu_N)$ and $\supp(\nu)$ converges to zero, then the
endpoint $\frak b$
is regular. Still assuming that $\frak d \notin\partial D$, let us now assume
instead that
\[
\nu_N = \frac{k}n \delta_\zeta+ \tilde
\nu_N,
\]
where $k$ is a fixed positive integer, $\zeta>0$ is fixed and lies outside
$\supp(\nu)$ and the Hausdorff distance between $\supp(\tilde\nu
_N)$ and
$\supp(\nu)$ converges to zero. The eigenvalue $\zeta$ of ${\bolds
\Sigma}_N$ with
multiplicity $k$ is often called a \textit{spike}. Assume without loss of
generality
that $\frak b$ is a right edge and that $1/\zeta$ belongs to the same connected component of $D$ as
$\frak d$.
Three situations that we describe without formal proofs are of interest:
\begin{longlist}[(2)]
\item[(1)] The spike $\zeta$ satisfies $g'(1/\zeta) < 0$. This can
only happen
if $1/\zeta< \frak d$, as shown in Figure~\ref{figoutlier}. In this case,
$\zeta$ produces $k$ \textit{outliers}, that is, eigenvalues
of $\widetilde{\mathbf M}_N$ which converge to a value outside the bulk;
see~\cite{baik-silverstein-2006,benaych-rao-2011}. In terms of the
support of
$\mu(\gamma_N, \nu_N)$, the location of these outliers corresponds
to a
small interval in $\supp(\mu(\gamma_N, \nu_N))$ (see Figure~\ref
{figoutlier})
which is absent from $\supp(\mu(\gamma,\nu))$. The width of this
new interval
is of order $N^{-1/2}$.

%
%f2 #&#
\begin{figure}

\includegraphics{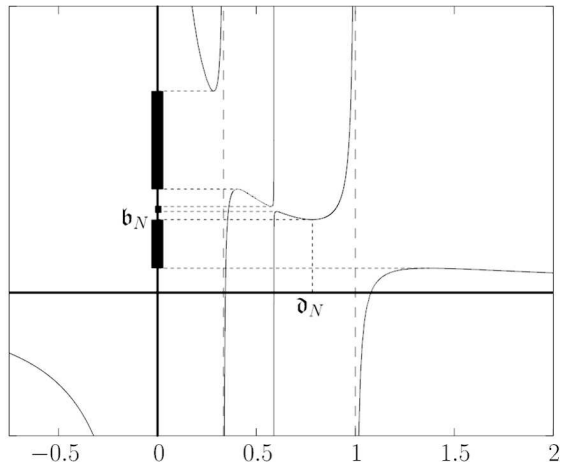}

\caption{Plot of $g_N(x)$ for $n = 300$, $\gamma_N = 0.1$ and
$\nu_N = \frac{1}{300} \delta_{1.7} + \frac{209}{300}\delta_1
+ \frac{90}{300}\delta_3$. The spike $\zeta=1.7$ produces an
outlier. The
asymptote at $1/\zeta$ is not shown for better visibility.}\label{figoutlier}
\end{figure}

Since $1/\zeta< \frak d$, the regularity condition still holds for
$\frak b$,
and Tracy--Widom fluctuations around $\frak b_N = g_N(\frak d_N)$ will be
observed.\vadjust{\goodbreak}

Let us say a few words on the fluctuations of the outliers. Notice that
$\zeta$
incurs the presence of a local minimum and a new local maximum in $g_N$
which are absent from $g$; see Figures~\ref{fignospike} and~\ref
{figoutlier}. Considering, for example, the minimum reached at, say
$\frak d'_N$, one can
show that $\llvert  1/\zeta- \frak d'_N\rrvert  $ is of order $N^{-1/2}$.
In particular, the regularity assumption~(\ref{regul-gN}) is not
satisfied for $\frak d'_N$.
In fact, it is known that when they are scaled by $N^{1/2}$, the $k$
outliers asymptotically fluctuate up to a multiplicative constant as the
eigenvalues of a $k\times k$ matrix taken from the GUE ensemble;
see~\cite{BBP-2005,bai-yao-ihp-08,benaych-guionnet-maida-11} among others.

\item[(2)] The spike $\zeta$ satisfies $g'(1/\zeta) > 0$. The case where
$1/\zeta> \frak d$ is shown in Figure~\ref{fignoout}. Here, the spike
$\zeta$ does not create an outlier, and the regularity condition on
$\frak b$ is still satisfied.
Tracy--Widom fluctuations around $\frak b_N = g_N(\frak d_N)$ will be
also observed here.

\item[(3)] The spike depends generally on $N$ and satisfies $1/\zeta
\to\frak d$
as $N\to\infty$.
Here, we are at the crossing point of the phase transition discovered in
\cite{BBP-2005} between the ``Tracy Widom regime'' and the ``GUE regime.''
More specifically, under an additional condition [see (\ref
{speedcondition})] we shall briefly outline in Appendix~\ref{appBBP}
that at the scale
$N^{2/3}$ the asymptotic fluctuations are described by the so-called deformed
Tracy--Widom law whose distribution function $F_k$ is defined
in~\cite{BBP-2005}, equation~(17). One can also be interested in the
regime where $k=k(N)\to\infty$ as $N\to\infty$. In the setting of
additive perturbations of Wigner matrices, this situation has been
considered by P\'ech\'e when $k/N\to0$, and she proved Tracy--Widom
fluctuations arise; see \cite{peche-ptrf-2006}, Theorem 1.5. We do not
pursue this direction here.
\end{longlist}

All these arguments can be straightforwardly generalized to the case
where a
finite number of different spikes are present.

As explained in the third point above and in Appendix~\ref{appBBP},
we can tackle the situation where an edge satisfies a weak kind of
nonregularity. Nevertheless, our approach breaks down in the case of a
limiting measure $\nu$ for which Proposition~\ref{endpoints}(b) occurs.\vspace*{-1pt}

%
%f3 #&#

%s3.2 #&#
\subsection{Asymptotic independence}
\label{secAIresults}

Our next result states that the fluctuations of the extremal
eigenvalues associated with any finite number of regular soft edges are
asymptotically independent.

%
%th4 #&#
\begin{theorem}
\label{thasymptindep}
Let Assumptions~\ref{assgauss}~and~\ref{assnu} hold true, and let $I$
and $J$ be finite sets of indices. Denote by $(\frak a_i)_{i\in I}$
left regular soft edges and by
$(\frak b_j)_{j\in J}$ right regular soft edges.

Let $\tilde x_{\varphi_i(N)}$ and $\frak c_{i,N}$ be associated to
$\frak a_i$ as in Proposition~\ref{propproperties-regular-left}, and set
\[
\frak a_{i,N}=g_N(\frak c_{i,N}),\qquad
\sigma_{i,N}= \biggl(\frac
{2}{-g_N''(\frak c_{i,N})} \biggr)^{1/3}.
\]
Similarly, let $\tilde x_{\phi_j(N)}$ and $\frak d_{j,N}$ be
associated to $\frak b_j$ as in Proposition~\ref
{propproperties-regular-right}, and
\[
\frak b_{j,N}=g_N(\frak d_{j,N}),\qquad
\delta_{j,N}= \biggl(\frac
{2}{g_N''(\frak d_{j,N})} \biggr)^{1/3}.
\]
Then, for every real numbers $(s_i)_{i\in I}$, $(t_j)_{j\in J}$, we have
\begin{eqnarray*}
&& \lim_{N\rightarrow\infty} \p \bigl( N^{2/3}\sigma_{i,N}
(\frak a_{i,N}-x_{\varphi_i(N)} )\leq s_i,
\\[-4pt]
&&\hspace*{37pt} i\in
I, N^{2/3}\delta_{j,N} (x_{\phi_j(N)} -\frak
b_{j,N} )\leq t_j,   j \in J \bigr)
\\[-1pt]
&&\qquad= \prod_{i\in I}\det (I-\K_{\Aii}
)_{L^2(s_i,\infty
)}\prod_{j\in J}\det (I-
\K_{\Aii} )_{L^2(t_j,\infty)}.
\end{eqnarray*}
\end{theorem}

We prove Theorem~\ref{thasymptindep} in Section~\ref
{secasymptotic-independence}. Our strategy is to build on the
operator-theoretic proof of Bornemann in the case of the smallest and
largest eigenvalues of the GUE \cite{bornemann-2010}; it essentially
amounts to proving that the off-diagonal entries of a two-by-two
operator valued matrix decay to zero in the trace class norm. In our
setting, the problem involves a larger operator valued matrix, and we
show that obtaining the decay to zero for the off-diagonal entries in
the Hilbert--Schmidt norm is actually sufficient. We establish the
latter by using the estimates established in Section~\ref{SectionTW}.

\subsubsection*{A comment on universality} The results presented in this
paper rely on the fact that the entries of $\mathbf X_N$ are complex
Gaussian random variables,
a key assumption in order to take advantage of the determinantal
structure of the eigenvalues of the model under study. A recent work by
Knowles and Yin \cite{knowles-yin-2014-preprint} enables one to
transfer the results presented here (except the hard edge fluctuations;
see Theorem~\ref{thBessel} below) to the case of complex, but not
necessarily Gaussian, random variables. Indeed, by combining the local
convergence to the limiting distribution established in~\cite
{knowles-yin-2014-preprint} together with Theorems~\ref
{thfluctuations-TW}~and~\ref{thasymptindep}, one obtains Tracy--Widom fluctuations and
asymptotic independence in this more general setting, provided that the
entries of matrix $\mathbf{X}_N$ fulfill some moment condition. This
also provides a similar generalization of our Proposition~\ref
{propcondition-number}, describing the asymptotic behavior for the
condition number of $\mathbf M_N$ when $\gamma>1$. Let us stress that
the case of
real Gaussian random variables (except the largest one covered in \cite
{lee-schnelli-preprint}), of great importance in statistical
applications, remains open.

%Let us also mention the case of covariance matrices with
%real entries, of important interest in statistical applications. In
%the recent work \cite{lee-schnelli-preprint}, Lee and Schnelli prove
%(GOE) Tracy-Widom fluctuations for the maximal eigenvalue under rather
%weak assumptions. However, a complete description for the fluctuations
%at each edge, and in particular at the hard edge, seems to remain open.

%s3.3 #&#
\subsection{Fluctuations at the hard edge}\label{secdescription-BE}

Proposition~\ref{1111} shows that when the leftmost edge is a hard edge,
$\gamma= 1$. (Actually, one can show that this is an equivalence.)
In order to study the smallest random eigenvalue fluctuations at the hard
edge, we restrict ourselves to the case where $n=N+\alpha$, where
$\alpha\in\mathbb Z$ is independent of $N$. Thus the
smallest random eigenvalue of $\mathbf M_N$ is
\[
x_{\min}= \cases{ x_1=\tilde x_{\alpha+1}, &\quad if
$\alpha\geq0$,
\cr
x_{1-\alpha}=\tilde x_{1}, &\quad if $
\alpha<0$.}
\]
We shall prove that the fluctuations of $x_{\min}$ around the origin
are described by mean of the Bessel kernel with parameter $\alpha$,
%follows the hard-edge Tracy-Widom distribution of parameter $\alpha$,
that we introduce now.

The Bessel function of the first kind $J_\alpha$ with parameter
$\alpha\in\mathbb Z$ is defined~by
%
%e21 #&#
\begin{equation}
\label{seriesrepBessel} J_\alpha(x) = \biggl( \frac{x}2
\biggr)^\alpha\sum_{n=0}^\infty
\frac
{(-1)^n}{n!  \Gamma(n+\alpha+1)} \biggl( \frac{x}2 \biggr)^{2n},\qquad x>0.
\end{equation}
Note that when $\alpha<0$, the first $\llvert  \alpha \rrvert  $ terms in the series
vanish since the Gamma function $\Gamma$ has simple poles on the
nonpositive integers. Denote by $\K_{\mathrm{Be},\alpha}$ the Bessel kernel
%
%e22 #&#
\begin{equation}
\label{Besselkernel} \K_{\mathrm{Be},\alpha}(x,y)=\frac{\sqrt{y}J_\alpha(\sqrt
{x})J_\alpha'(\sqrt{y})-\sqrt{x}J_\alpha'(\sqrt{x})J_\alpha(\sqrt
{y})}{2(x-y)},
\end{equation}
and by extension, $\K_{\mathrm{Be},\alpha}$, the associated integral
operator. Given a nonnegative real-valued random variable $X$, the following
probability distribution will be of particular interest:\vspace*{-2pt}
\[
\mathbb{P}(X\geq s) = \det (I-\K_{\mathrm{Be},\alpha} )_{L^2(0,s)},\qquad s>0,
\]
where the right-hand side stands for the Fredholm determinant of the
restriction to $L^2(0,s)$ of the integral operator $\K_{\mathrm{Be},
\alpha}$. When $\alpha=0$, this is actually the distribution of an
exponential law of parameter $1$, namely $\det (I-\K_{\mathrm{Be},0} )_{L^2(0,s)}=e^{-s}$. Also\vspace*{1pt} of interest is the
%As for the name, it comes from the following
alternative representation due to Tracy and Widom \cite{tracy-widom-bessel-94},
\[
\det (I-\K_{\mathrm{Be},\alpha} )_{L^2(0,s)}=\exp \biggl(-\frac
{1}{4}\int
_0^s(\log s -\log x) q(x)^2\,\d x
\biggr),
\]
where $q$ is the solution of a differential equation which is reducible
to a particular case of the Painlev\'e V equation (involving $\alpha$
in its parameters) and boundary condition $q(x)\sim J_\alpha(\sqrt x)$
as $x\to0$.

Let us now state our result for the fluctuations around the hard edge.

%th5 #&#
\begin{theorem}
\label{thBessel}
Let Assumptions~\ref{assgauss}~and~\ref{assnu} hold true; assume
moreover that $n=N+\alpha$, where $\alpha\in\mathbb Z$ is
independent of $N$. Set\vspace*{-2pt}
%
%e23 #&#
\begin{equation}
\label{varBessel} \sigma_N=\frac{4}{N}\sum
_{j={1}}^n\frac{1}{\lambda_j}.
\end{equation}
Then, for every $s>0$, we have\vspace*{-2pt}
%
%e24 #&#
\begin{equation}
\lim_{N\rightarrow\infty}\p \bigl(N^2\sigma_N
x_{\min}\geq s \bigr)=\det (I-\K_{\mathrm{Be},\alpha} )_{L^2(0,s)}.
\end{equation}
In particular, if $N=n$, then we have for every $s>0$,
%
%e25 #&#
\begin{equation}
\lim_{N\rightarrow\infty}\p \bigl(N^2\sigma_N
x_{\min}\geq s \bigr)=e^{-s}.
\end{equation}
\end{theorem}

%re3.1 #&#
\begin{Remark}
The assumption that $\nu_N$ converges weakly toward some limit $\nu$
is actually not used in the proof of Theorem~\ref{thBessel}. Namely,
this result holds true under Assumptions~\ref{assgauss}~and~\ref{assnu}(2) only.
\end{Remark}

We provide a proof for Theorem~\ref{thBessel} in Section~\ref
{Besselsection}. It is also based on an asymptotic analysis for the
rescaled kernel; the key observation here is that when an edge is the
hard edge, the associated critical point $\frak c$ should be located at
infinity (when embedding the complex plane into the Riemann sphere).

\subsection{Application: Condition numbers}
\label{conditionnumber}

The condition number of the matrix ${\mathbf M}_N$ with eigenvalues $0\leq
x_1\leq\cdots\leq x_N$
is defined by
\[
\kappa_N=\frac{x_N}{x_1}, %
\]
provided it is finite, that is, $n/N \geq1$. If $n/N < 1$, one may
instead\vspace*{1pt} consider the condition number associated to $\widetilde{\mathbf
M}_N$, defined as $\tilde\kappa_N= \tilde x_n / \tilde x_1$. The
study of condition numbers is important in numerical linear algebra
\cite{von-neumann-goldstine-47,von-neumann-goldstine-51}, and random
matrix theory has already provided interesting theoretical \cite
{edelman-1988-SIAM,basor-et-al-2012} and applied \cite
{mckay-condition-number-2010,bianchi-et-al-2011} results. As a
consequence of our former results, we provide an asymptotic study for
$\kappa_N$. (One can easily derive similar results for $\tilde\kappa_N$.)

\subsubsection*{Notation} We use the notation $\stackrel{\mathcal D}{\to
}$ for the convergence in distribution of random variables.

%
%pr3.2 #&#
\begin{proposition} \label{propcondition-number}
Let Assumptions~\ref{assgauss}~and~\ref{assnu} hold true and $\gamma
>1$. Let $\frak a$ be the leftmost edge, assume it is regular and let
$(\frak c_N)_N$ and $\frak c$ be as in Proposition~\ref
{propproperties-regular-left}. Let $\frak b$ be the rightmost edge,
assume it is regular and let $(\frak d_N)_N$ and $\frak d$ be as in
Proposition~\ref{propproperties-regular-right}.
Set
\begin{eqnarray*}
\frak a_N &=& g_N(\frak c_N),\qquad
\sigma_N= \biggl(\frac
{2}{-g_N''(\frak c_N)} \biggr)^{1/3},
\\
\frak b_N &=& g_N(\frak d_N),\qquad
\delta_N= \biggl(\frac
{2}{g_N''(\frak d_N)} \biggr)^{1/3}.
\end{eqnarray*}
Assume moreover that $ x_1\to\frak a$ and $ x_N\to\frak b$ a.s. Then
\[
\kappa_N \mathop{\longrightarrow}_{N\to\infty}^{a.s.}
\frac{\frak
b}{\frak a}\quad\mbox{and}\quad N^{2/3} \biggl(
\kappa_N - \frac{\frak b_N}{\frak a_N} \biggr) \mathop{\longrightarrow}_{N\to\infty}^{\mathcal D}
\frac
{X} {\delta\frak a} + \frac{\frak b Y}{\sigma\frak a^2},
\]
where $X$ and $Y$ are two independent Tracy--Widom distributed random
variables, and where
\[
\sigma= \biggl(\frac{2}{-g''(\frak c)} \biggr)^{1/3}=\lim
_{N\to
\infty} \sigma_N \quad\mbox{and} \quad\delta=
\biggl(\frac
{2}{g''(\frak d)} \biggr)^{1/3}=\lim_{N\to\infty}
\delta_N. %
\]
%
%where $\frak c$ and $\frak d$ are defined in \eqref{leftsoftedge} and
%\eqref{rightsoftedge}.
\end{proposition}

\begin{Remark}
The condition that $ x_1\to\frak a$ and $ x_N\to
\frak b$ a.s. imposes that neither $x_N$ nor $x_1$ are outliers;
otherwise their fluctuations (together with those of $\kappa_N$) would
be of order $N^{1/2}$, and a different (somewhat easier) asymptotic
analysis should be conducted. We do not pursue in this direction here.
\end{Remark}

\begin{pf*}{Proof of Proposition \ref{propcondition-number}}
Only the convergence in distribution requires an argument. Write
\begin{eqnarray*}
N^{2/3} \biggl( \kappa_N - \frac{\frak b_N}{\frak a_N} \biggr) &=&
N^{2/3} \biggl( \frac{ x_N}{x_1} - \frac{\frak
b_N}{\frak a_N} \biggr) =
N^{2/3} \biggl( \frac{\frak a_N x_N - \frak b_N x_1}{x_1\frak a_N} \biggr)
\\
&=& \frac{N^{2/3}}{x_1 \frak a_N} \bigl\{ \frak a_N ( x_N - \frak
b_N ) - \frak b_N ( x_1 - \frak
a_N ) \bigr\}
\\
&=& \frac{1}{ x_1 \delta_N} N^{2/3} \delta_N ( x_N -
\frak b_N ) +\frac{\frak b_N}{x_1 \frak a_N \sigma_N} N^{2/3} \sigma
_N (\frak a_N - x_1 ).
\end{eqnarray*}
Using the asymptotically independent Tracy--Widom fluctuations of
$N^{2/3}\times\break  \delta_N  ( x_N - \frak b_N )$ and $ N^{2/3}
\sigma_N  (\frak a_N - x_1 ) $ (cf. Theorems~\ref
{thfluctuations-TW}~and~\ref{thasymptindep}) together with the a.s.
convergence $x_1\to\frak a$ and the convergences $\frak a_N \to\frak
a$, $\frak b_N\to\frak b$, $\delta_N \to\delta$ and $ \sigma_N \to
\sigma$ (cf. Proposition~\ref{gN->g}), one can conclude using
Slutsky's lemma \cite{book-vandervaart}, Lemma 2.8.
\end{pf*}
We now handle the case where $\gamma=1$.

%
%pr3.4 #&#
\begin{proposition}\label{propcondition-number-gamma-1}
Let Assumptions~\ref{assgauss}~and~\ref{assnu} hold true, and let
$n=N+\alpha$ where $\alpha\in\mathbb{N}$ is independent of $N$. Let
\[
\sigma_N=\frac{4}N \sum_{j=1}^n
\frac{1}{\lambda_j}\quad\mbox {and}\quad\sigma=4\int\frac{1}{x} \,\d\nu(x) =
\lim_{N\to\infty} \sigma_N.
\]
Assume that a.s. $x_N\to\frak b$ for some $\frak b>0$. Then
\[
\frac{1}{N^2} \kappa_N \mathop{\longrightarrow}_{N\to\infty
}^{\mathcal D}
\frac{\frak b \sigma} {X},
\]
where $X$ is a random variable with distribution
\[
\p (X \geq s )=\det (I-\K_{\mathrm{Be},\alpha} )_{L^2(0,s)},\qquad s> 0.
\]
\end{proposition}

\begin{pf}
Write
\[
\frac{\kappa_N}{N^2} = \frac{\sigma_N(x_N -\frak b)}{N^2 \sigma_N
x_1} +\frac{\sigma_N \frak b}{N^2 \sigma_N x_1}. %
\]
Since by assumption $x_N - \frak b\to0$ a.s. and by Theorem~\ref
{thBessel} $(N^2 \sigma_N x_1)^{-1} \to X^{-1}$ in distribution, where
$X$ has the distribution specified in the statement, we have
\[
\frac{\sigma_N(x_N -\frak b)}{N^2 \sigma_N x_1} \mathop {\longrightarrow}_{N\to\infty}^{\mathcal D} 0.
\]
By Slutsky's lemma, $N^{-2} \kappa_N$ then converges toward $\frak b
\sigma X^{-1}$ in distribution.
\end{pf}
%
%Proof of Proposition~\ref{propcondition-number-gamma-1} is left to the
%reader.

%re3.5 #&#
\begin{Remark}
Interestingly, in the square case where $\gamma=1$, the fluctuations
of the
largest eigenvalue $x_N$ (either of order $N^{1/2}$ if $x_N$ is an
outlier or of
order $N^{2/3}$ in the Tracy--Widom regime) have no influence on the
fluctuations of $\kappa_N$ as these are imposed by the limiting distribution
of $x_1$ at the hard edge.
\end{Remark}

%s4 #&#
\section{Proof of Theorem \texorpdfstring{\protect\ref{thfluctuations-TW}}{3}: Tracy--Widom fluctuations}\label{SectionTW}
This section is devoted to the proof of Theorem~\ref{thfluctuations-TW}.

%s4.1 #&#
\subsection{Outline of the proof}
\label{outlineproof}

\subsubsection*{Step 1 (preparation)} As in \cite{BBP-2005} and \cite
{elkaroui-07-TW}, the starting point to establish Tracy--Widom
fluctuations is that the random eigenvalues of $\mathbf M_N$ or
$\widetilde{\mathbf M}_N$ form a determinantal\vadjust{\goodbreak} point process, so that
the gap probabilities can be expressed as Fredholm determinants of an
integral operator $\K_N$ with kernel $\K_N(x,y)$. We provide all the
necessary material from operator theory in Section~\ref
{determinantalprocess}. In Section~\ref{kernelofWishart} we first
recall the double contour integral formula for $\K_N(x,y)$ obtained in
\cite{BBP-2005,onatski-2008-TW}. Next, we show using Theorem~\ref
{thdef-phi-N} that one can represent the cumulative distribution
functions for the extremal eigenvalues as Fredholm determinants
involving $\K_N$ asymptotically. As a consequence, proving the
Tracy--Widom fluctuations boils down to establishing the appropriate
convergence of rescaled versions $\widetilde\K_N(x,y)$ of the kernel
$\K_N(x,y)$ toward the Airy kernel. To this end, we split $\widetilde
\K_N(x,y)$ into two parts, $\K_N^{(0)}(x,y)$ and $\K_N^{(1)}(x,y)$,
each involving different integration contours.

\subsubsection*{Step 2 (contours deformations)} Anticipating the
forthcoming asymptotic analysis, we focus in Section~\ref
{contourssection} on right edges\vspace*{1pt} and prove the existence of appropriate
integration contours coming with $\K_N^{(0)}(x,y)$ and $\K
_N^{(1)}(x,y)$; the case of a left edge is deferred to Section~\ref
{variationsleftedges}. Obtaining appropriate explicit contours is
usually the hard part in the asymptotic analysis; see, in particular,
\cite{elkaroui-07-TW}. Here, we instead provide a nonconstructive
proof for the existence of appropriate contours by mean of the maximum
principle for subharmonic functions, which has the advantage to work
for every regular edge up to minor modifications.

\subsubsection*{Step 3 (asymptotic analysis)} Still focusing on the right
edge setting, we prove in Section~\ref{K0playsnorole} that $\K
_N^{(0)}(x,y)$ does not contribute in the large $N$ limit. Moreover, we
prove the convergence of kernel $\K_N^{(1)}(x,y)$ to the Airy kernel
in an appropriate sense and then complete the proof of Theorem~\ref
{thfluctuations-TW}(b). For this last step, we use a different approach
than in \cite{BBP-2005,elkaroui-07-TW}: instead of relying on a
factorization trick and the H\"older inequality to obtain the trace
class convergence, we use an argument involving the regularized
Fredholm determinant $\det_2$ to show the convergence of the Fredholm
determinants. Finally, in Section~\ref{variationsleftedges}, we adapt
the arguments to the left edge setting and complete the proof of
Theorem~\ref{thfluctuations-TW}.

%s4.2 #&#
\subsection{Operators, Fredholm determinants and determinantal processes}\label{determinantalprocess}

% {\red[jam: me semble {\chr"C3}{\chr"A0} peu pr{\chr"C3}{\chr"A8}s
%clean, peut-{\chr"C3}{\chr"AA}tre rajouter un rappel sur l'expression
%explicite de la trace $\int K(x,x)  \,dx$ d'un op{\chr"C3}{
%\chr"A9}rateur {\chr"C3}{\chr"A0} noyau $K(x,x)$]}

\subsubsection*{Trace class operators and Fredholm determinants}
We provide hereafter a few elements of operator theory; for classical
references, see \cite{book-davies-operator,book-gohberg-2000,book-simon-2005}.
Consider a compact linear operator $\A$ acting on a separable Hilbert
space ${\mathcal H}$ [we write $\A\in L({\mathcal H})$], and denote by
$(s_n)_{n=1}^\infty$ the \textit{singular values} of $\A$ repeated
according to their multiplicities, that is, the eigenvalues of $(\A\A
^*)^{1/2}$. The set
\[
{\mathcal J}_1= \Biggl\{ \A\in L({\mathcal H}), \sum
_{n=1}^\infty s_n<\infty \Biggr\}
\]
is the (sub-)algebra of \textit{trace class operators} and endowed with
the norm
$
{\llVert  \A\rrVert  }_1=\sum_{n=1}^\infty s_n
$; $({\mathcal J}_1,\llVert  \cdot\rrVert  _1)$ is complete. If $\A\in{\mathcal
J}_1$ with eigenvalues $(a_n)_{n=1}^\infty$ (repeated according to
their multiplicities), then the \textit{trace} and the \textit{Fredholm
determinant} of~$\A$,
\[
\Tr(\A)=\sum_{n=1}^\infty a_n
\quad\mbox{and}\quad\operatorname {det}(I-\A) = \prod
_{n=1}^\infty(1-a_n), %
\]
are well defined and finite (Lidskii's trace theorem). The maps $\A
\mapsto\Tr(\A)$ and $\A\mapsto\operatorname{det}(I-\A)$ are
continuous on $({\mathcal J}_1,{\llVert   \cdot\rrVert  }_1)$. If both $\A\B$ and
$\B\A$ are trace class, then we have the useful identity
%
%e26 #&#
\begin{equation}
\label{switchidentity} \det(I-\A\B)=\det(I-\B\A).
\end{equation}
Similarly, let
\[
{\mathcal J}_2= \Biggl\{ \A\in L({\mathcal H}), \sum
_{n=1}^\infty s^2_n<\infty
\Biggr\} %
\]
be the (sub-)algebra of \textit{Hilbert--Schmidt operators} endowed with
the norm
$
{\llVert  A\rrVert  }_2= \{ \sum_{n=1}^\infty s^2_n \} ^{1/2}
$. The set $({\mathcal J}_2,\llVert   \cdot\rrVert  _2)$ is complete. If $\A\in
{\mathcal J}_2$ with eigenvalues $(a_n)_{n=1}^\infty$ (repeated
according to their multiplicities), then the \textit{regularized} 2-\textit{determinant} of~$\A$,
%
%e27 #&#
\begin{equation}
\label{det2} \operatorname{det}_2(I-\A) = \prod
_{n=1}^\infty(1-a_n)e^{a_n},
\end{equation}
is well defined and finite. Moreover, the map $\A\mapsto\operatorname
{det}_2(I-\A)$ is continuous on $({\mathcal J}_2,\llVert  \cdot\rrVert  _2)$.

The inclusion ${\mathcal J}_1\subset{\mathcal J}_2$ is
straightforward. The H\"older inequality $\llVert   \A\B\rrVert  _1\le\llVert   \A\rrVert  _2 \llVert
\B\rrVert  _2$ yields that if $\A,\B$ are Hilbert--Schmidt, then both $\A
\B$ and $\B\A$ are trace class. The following simple property will
play a key role in the sequel:

%pr4.1 #&#
\begin{proposition}\label{propdet2}
If $\A\in{\mathcal J}_1$, then
\[
\mathrm{det}_2(I-\A) = \mathrm{det}(I-\A) e^{\Tr(\A)}.
\]
As a consequence, if the operators $\A_n,\A\in{\mathcal J}_1 $ are
such that $\Tr(\A_n) \to\Tr(\A)$ and ${\llVert   \A_n-\A\rrVert  }_2\to0$ as
$n\to\infty$, then
\[
\mathrm{det}(I-\A_n) \mathop{\longrightarrow}_{n\to\infty}
\mathrm{det}(I-\A). %
\]
\end{proposition}

\begin{longlist}
\item[\textit{Integral operators}.]
When working on $\mathcal H=L^2(\R)$, we identify a given kernel
$(x,y)\mapsto\K(x,y)$ with its associated integral operator $\K f =
\int\K(\cdot,y)f(y)  \,\d y$ acting on $L^2(\R)$, provided the
latter makes sense.
Let $J\subset\R$ be a Borel set and $\mathbf1_J$ be the orthogonal
projection of $L^2(\R)$ onto $L^2(J)$. The restriction ${\K\mid}_J$ of
$\K$ to $L^2(J)$ is defined by
\[
{\K\mid}_J f(x) = {\mathbf1}_J(x) \int
_J \K(x,y) f(y) \,\d y, \qquad f \in L^2(J),
\]
and is associated to the kernel $(x,y)\mapsto{\mathbf1}_J(x) \K(x,y)
{\mathbf1}_J(y)$, namely ${\K\mid}_J={\mathbf1}_J \K{\mathbf1}_J$. In
order to keep track of these projections when dealing with Fredholm
determinants, we shall often write $\det(I-\K)_{L^2(J)}$ for $\det
(I-\mathbf1_J \K\mathbf1_J)$.

Given a measurable kernel $\K\dvtx \R\times\R\to\R$, the associated
integral operator $\K$ on $L^2(\R)$ is Hilbert--Schmidt if and only if
\[
\int_{\R}\int_{\R}
\K(x,y)^2 \,\d x\,\d y <\infty,
\]
and in this case we have
%
%e28 #&#
\begin{equation}
\label{defHS-kernel} {\llVert \K\rrVert }_2= \biggl(\int
_\R\int_{\R} \K(x,y)^2 \,\d
x\,\d y \biggr)^{1/2}.
\end{equation}
We finally recall (cf. \cite{book-gohberg-2000}, Theorem 8.1) that if
$\K\dvtx [a,b]\times[a,b]\to\R$ is a continuous kernel whose associated
operator $\mathbf 1_{(a,b)}\K\mathbf 1_{(a,b)}$ is trace
class\footnote{See,
for instance, \cite{book-gohberg-2000}, Theorem 8.2, for sufficient
conditions on $\K$ to be trace class.} on $L^2(\R)$, then
%
%e29 #&#
\begin{equation}
\label{defTrace-kernel} \Tr (\mathbf 1_{(a,b)}\K\mathbf 1_{(a,b)} )
= \int_a^b \K(x,x) \,\d x.
\end{equation}
\end{longlist}

\begin{longlist}
\item[\textit{Convention}.] From this point forward, the trace $\Tr$ and
the Hilbert--Schmidt norm ${\llVert   \cdot \rrVert  }_2$ will always refer to the
Hilbert space $L^2(\R)$.
\end{longlist}

\begin{longlist}
\item[\textit{Determinantal point process}.]
Real random variables $x_1,\ldots,x_m$ are said to form a
determinantal point process with kernel $\K\dvtx \R\times\R\rightarrow
\R$ (and Lebesgue measure for reference measure) if its gap
probabilities are expressed as Fredholm determinants; namely, for any
Borel set $J\subset\R$, we have
\[
\p \bigl(\sharp \{1\leq k \leq m\dvtx  x_{k}\in J \}=0 \bigr)=\det (I-\K
)_{L^2(J)},
\]
provided that the right-hand side makes sense; the latter stands for
the Fredholm determinant of the restriction to $L^2(J)$ of the integral
operator with kernel $\K(x,y)$.
\end{longlist}

%s4.3 #&#
\subsection{The kernel of a correlated Wishart matrix and its properties}
\label{kernelofWishart}

The next proposition will be of fundamental use in this paper.

%
%pr4.2 #&#
\begin{proposition} \label{propdet-representation}
Let\vspace*{1.5pt} Assumption~\ref{assgauss} hold true. Then, for every $N$, the
$\min(n,N)$ random eigenvalues of $\widetilde{\mathbf M}_N$ (and
equivalently of $\mathbf M_N$) form a determinantal point process
associated with the kernel
%
%e30 #&#
\begin{eqnarray}
\label{kernel}
\K_{N}(x,y)
&=& \frac{ N}{(2i\pi)^2}\oint_{\Gamma}\d z\oint_{\Theta} \d w\, e^{-
Nx(z-q) +Ny(w-q)}
\nonumber\\[-8pt]\\[-8pt]\nonumber
&&\hspace*{71pt}{}\times \frac{1}{w-z} \biggl(\frac{z}{w}
\biggr)^{ N}\prod_{j=1}^n
\biggl(\frac{w-\lambda_j^{-1}}{z-\lambda_j^{-1}} \biggr),
\end{eqnarray}
where the real $q\in(0,\lambda_n^{-1})$ is a free parameter, and we
recall that the $\lambda_i$'s are the eigenvalues of $\bolds\Sigma_N$.
$\Gamma$ and $\Theta$ are disjoint closed contours, both oriented
counterclockwise, such that $\Gamma$ encloses the $\lambda_j^{-1}$'s
and lies\vspace*{2pt} in $\{z\in\C\dvtx   \re z>q\}$, whereas $\Theta$ encloses the
origin and lies in $\{z\in\C\dvtx   \re z<q\}$.
\end{proposition}

\begin{figure}

\includegraphics{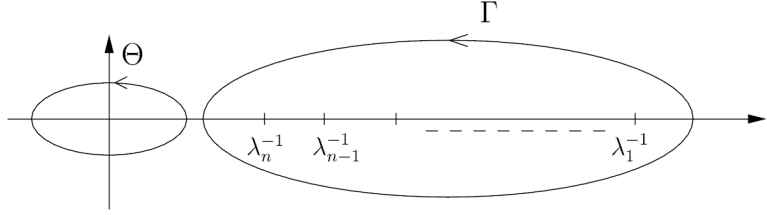}

\caption{The contours of integration.}\label{figpaths}
\end{figure}

By convention, all the contours we shall consider will be assumed to be
simple and oriented counterclockwise. The integration contours are
shown in Figure~\ref{figpaths}.%{\red[jam: possible (mais pas forc{
%\chr"C3}{\chr"A9}ment indispensable) de faire figurer $q$ sur la
%figure 4?]}\\

This proposition can be found in \cite{BBP-2005} ($n/N\leq1$) where
it is attributed to Johansson, and in \cite{onatski-2008-TW} ($n/N>
1$). Notice that
since the pioneering work of Br\'ezin and Hikami \cite
{brezin-hikami-98}, many such double integral representations appeared
for determinantal point processes.

%
%f4 #&#

%re4.3 #&#
\begin{Remark}%[on the localization of the free parameter $q$]
\label{freedomq}
The assumption over $q$, that is, $q\in(0,\lambda_n^{-1})$, ensures
that $\K_N$ with kernel (\ref{kernel}) is trace class on $L^2(\R)$.
In the sequel, we shall only need $\K_N$ to be \textit{locally trace
class}, that is, trace class on $L^2(J)$ for every compact subset
$J\subset\R$.
As an important consequence, we can choose $q\in\mathbb{R}$ with no
further restriction. In fact, let $q\in(0,\lambda_n^{-1})$, $q'\in\R
$ and $J\subset\R$ be a compact set. Then the multiplication operator
$\E\dvtx f(x)\mapsto e^{(q'-q)Nx}f(x)$ and its inverse $\E^{-1}$ are trace
class on $L^2(J)$. Write $\K_N = \K_N \E^{-1} \E$, and use (\ref
{switchidentity})
to get
\[
\det(I-\K_N)_{L^2(J)}=\det\bigl(I-\E\K_N
\E^{-1}\bigr)_{L^2(J)}. %
\]
The kernel of $\E\K_N \E^{-1}$ is simply obtained by (\ref{kernel})
where $q$ has been replaced by~$q'$, and our claim follows.
\end{Remark}

\subsubsection*{Asymptotic determinantal representation for the law of extremal eigenvalues}
Recall that to prove Tracy--Widom fluctuations for the maximal
eigenvalue $\tilde x_{n}$ of $\widetilde{\mathbf M}_N$,
a classical way\vspace*{2pt} to proceed is to identify the events
$
\{ N^{2/3}\sigma_N(\tilde x_{n}- \frak b_N)\leq s\} = \{ \mbox{no }
\tilde x_i\mbox{'s in } (\frak b_N+s/(N^{2/3}\sigma_N),\infty)\} $,
to use the determinantal representation
\[
\p \bigl(N^{2/3}\sigma_N(\tilde x_{n}- \frak
b_N)\leq s \bigr)= \det(I-\K_N)_{L^2(\frak b_N +s/(N^{2/3}\sigma_N),\infty)}
\]
and to prove the convergence of operator $\K_N$ to the Airy operator
$\K_{\Aii}$ after the rescaling
$x\mapsto\frak b_N+x/(N^{2/3}\sigma_N)$ for the trace class topology.
This would yield the desired result since the Fredholm determinant is
continuous for that topology.\vspace*{1pt}

Since the probabilities of interest $\p(N^{2/3}\sigma_N(\tilde
x_{\phi(N)}-\frak b_N)\leq s)$ and\break $\p(N^{2/3}\sigma_N(\frak
a_N-\tilde x_{\varphi(N)})\leq s)$
can no longer
be expressed as gap probabilities in general, we provide below an
asymptotic Fredholm determinant representation as $N\to\infty$ for these.

%
%pr4.4 #&#
\begin{proposition}
\label{asymptoticFredholm}
Consider the setting of Theorem~\ref{thfluctuations-TW}, and recall
that by convention $\tilde x_0=0$ and $\tilde x_{n+1}=+\infty$. Then
the following
facts hold true:
\begin{longlist}[(a)]
\item[(a)]
For every $\varepsilon>0$ small enough and for every sequence $(\eta
_N)_N$ of
positive numbers satisfying $\lim_N \eta_N=+\infty$,
%
%e31 #&#
\begin{equation}
\label{leftFredholmrep} \p \bigl(\eta_N (\frak a_N-\tilde
x_{\varphi(N)} )\leq s \bigr)= \det (I-\K_N )_{L^2(\frak a_N-\varepsilon, \frak a_N -
s/{\eta_N})}+o(1)
\end{equation}
as $N\to\infty$.
\item[(b)]
For every $\varepsilon> 0$ small enough and for every sequence $(\eta_N)_N$
of positive numbers satisfying $\lim_N \eta_N=+\infty$,
%
%e32 #&#
\begin{equation}
\label{rightFredholmrep} \p \bigl(\eta_N (\tilde x_{\phi(N)}-\frak
b_N )\leq s \bigr) =\det (I-\K_N )_{L^2(\frak b_N+s/{\eta_N}, \frak
b_N+\varepsilon)}+o(1)
\end{equation}
as $N\to\infty$.
\end{longlist}
\end{proposition}

\begin{pf}
We only prove (b), proof of (a) being similar.
%vide a proof for Proposition~\ref{asymptoticFredholm}--(b) since
%Proposition~\ref{asymptoticFredholm}--(a) can be established by
%following the
%same lines.
Observe that Theorem~\ref{thdef-phi-N}(b) and the convergence $\frak
b_N\to\frak b$ yield together the
existence of $\varepsilon>0$ small enough such that
\begin{eqnarray}
\label{A1}
&& \p \bigl(\eta_N (\tilde x_{\phi(N)}-\frak
b_N )\leq s \bigr)
\nonumber\\[-8pt]\\[-8pt]\nonumber
&&\qquad = \p \bigl(\eta_N (\tilde x_{\phi(N)}-\frak b_{N} )
\leq s, \tilde x_{\phi(N)+1} \geq\frak b_N+\varepsilon
\bigr)+o(1)
\nonumber
\end{eqnarray}
as $N\rightarrow\infty$. Now, $\varepsilon$ being fixed, use the
determinantal representation to write
%
%e33 #&#
\begin{eqnarray}
\label{A2}
&& \det ( I- \K_N )_{L^2(\frak b_N+s/\eta_N,   \frak
b_N+\varepsilon)}
\nonumber\\[-8pt]\\[-8pt]\nonumber
&&\qquad = \p \bigl(\sharp \{\ell\leq k \leq n\dvtx  \frak b_N+s/
\eta_N \leq \tilde x_{k}\leq\frak b_N+
\varepsilon \}=0 \bigr),
\end{eqnarray}
where
$\ell= n - \min(N,n) +1$. Recall the notational convention in Remark
\ref{remconvention}; we obtain by splitting along disjoint events
%
%e34 #&#
\begin{eqnarray}
\label{A3} && \p \bigl(\sharp \{\ell\leq k \leq n\dvtx  \frak b_N+s/
\eta_N \leq\tilde x_{k}\leq\frak b_N+
\varepsilon \}=0 \bigr)
\nonumber
\\
&&\qquad= \p \bigl(\eta_N (\tilde x_{\phi(N)}-\frak
b_N ) \leq s, \tilde x_{\phi(N)+1} \geq\frak b_N+
\varepsilon \bigr)
\nonumber\\[-8pt]\\[-8pt]\nonumber
&&\quad\qquad{} + \p (\tilde x_\ell\geq\frak b_N+
\varepsilon )
\\
&&\quad\qquad{} + \sum_{k=\ell,  k  \neq  \phi(N)}^{n}
\p ( \tilde x_k\leq\frak b_N+s/\eta_N,
\tilde x_{k+1} \geq\frak b_N+\varepsilon ).
\nonumber
\end{eqnarray}
Since we have the upper bounds
\begin{eqnarray*}
\sum_{k=\ell}^{\phi(N)-1}\p ( \tilde x_k
\leq\frak b_N+s/\eta _N, \tilde x_{k+1} \geq
\frak b_N+\varepsilon )&\leq& \p (\tilde x_{\phi(N)}\geq\frak
b_N+\varepsilon ),
\\
\sum_{k=\phi(N)+1}^{n}\p ( \tilde x_k
\leq\frak b_N+s/\eta_N, \tilde x_{k+1} \geq
\frak b_N+\varepsilon ) & \leq& \p (\tilde x_{\phi(N)+1} \leq\frak
b_N+s/\eta_N ),
\end{eqnarray*}
we obtain from (\ref{A2}), (\ref{A3}), Theorem~\ref{thdef-phi-N}(b)
and the convergence $\frak b_N\to\frak b$ that
\begin{eqnarray}
\label{A4}
&& \det (I- \K_N )_{L^2(\frak b_N+s/\eta_N, \frak
b_N+\varepsilon)}
\nonumber\\[-8pt]\\[-8pt]\nonumber
&&\qquad =\p \bigl(\eta_N (\tilde x_{\phi(N)}-\frak b_N )
\leq s, \tilde x_{\phi(N)+1}\geq\frak b_N+\varepsilon
\bigr)+o(1).
\nonumber
\end{eqnarray}
Finally, (\ref{rightFredholmrep}) follows by combining (\ref{A1}) and
(\ref{A4}).
\end{pf}

\subsubsection*{Rescaling and splitting the kernel $\K_N$}
We introduce
hereafter the rescaled kernel $\widetilde{\K}_N$ and provide an
alternative integral representation with new contours. The aim is to
prepare the forthcoming asymptotic analysis for right regular edges.

Let $\frak b$ be a soft regular right edge. By Proposition~\ref
{propproperties-regular-right}, there exist $\frak d>0$ such that
%
%e35 #&#
\begin{equation}
\label{critg} \frak b=g(\frak d),\qquad g'(\frak d)=0,\qquad
g''(\frak d)>0,
\end{equation}
and an associated sequence $(\frak d_N)$ such that $g_N^{(k)}(\frak
d_N) \to g^{(k)} (\frak d)$. Set
%
%e36 #&#
\begin{equation}
\label{notations} \frak b_N=g_N(\frak d_N),
\qquad\delta_N= \biggl(\frac{2}{g''_N(\frak
d_N)} \biggr)^{1/3},
\end{equation}
so that we have
%
%e37 #&#
\begin{eqnarray}\label{infocN}
g'_N(\frak d_N) &=& 0,\qquad
\lim_{N\rightarrow\infty}\frak d_N=\frak d,
\nonumber\\[-8pt]\\[-8pt]\nonumber
\lim _{N\rightarrow\infty} \frak b_N&=&\frak b,\qquad\lim
_{N\rightarrow\infty}\delta_N= \biggl(\frac{2}{g''(\frak d)}
\biggr)^{1/3}.
\end{eqnarray}
In particular $\frak c_N$, $g_N''(\frak c_N)$ and $\sigma_N$ are positive
numbers for every $N$ large enough, and $(\sigma_N)_N$ is a bounded sequence.

It follows from the definition of the extremal eigenvalue $\tilde
x_{\phi(N)}$
(see Theorem~\ref{thdef-phi-N} and Proposition~\ref
{asymptoticFredholm}) that
for every $\varepsilon>0$ small enough,
\begin{eqnarray}
\label{changevariables}
&& \p \bigl(N^{2/3}\delta_N (\tilde
x_{\phi(N)}-\frak b_{N} )\leq s \bigr)
\nonumber\\[-8pt]\\[-8pt]\nonumber
&&\qquad = \det (I-\K_N )_{L^2(\frak b_{N}+s/{(N^{2/3}\delta_N)},
\frak b_{N}+\varepsilon)} +o(1)
\nonumber
\end{eqnarray}
as $N\rightarrow\infty$. By a change of variable, we can write
\begin{eqnarray}
\label{resc}
&& \det (I-\K_N )_{L^2(\frak b_{N}+s/{(N^{2/3}\delta_N)},
\frak b_{N}+\varepsilon)}
\nonumber\\[-8pt]\\[-8pt]\nonumber
&&\qquad = \det (I-\mathbf 1_{(s,
\varepsilon
N^{2/3}\delta_N)}\widetilde{\K}_N
\mathbf 1_{(s, \varepsilon
N^{2/3}\delta_N)} )_{L^2(s,\infty)},
\nonumber
\end{eqnarray}
where the scaled integral operator $\widetilde{\K}_{N}$ has kernel
%
%e38 #&#
\begin{equation}
\label{rescaledkernel} \widetilde{\K}_N(x,y)=\frac{1}{N^{2/3}\delta_N}
\K_N \biggl(\frak b_{N}+\frac{x}{N^{2/3}\delta_N},\frak
b_{N}+\frac{y}{N^{2/3}\delta
_N} \biggr)
\end{equation}
with $\K_N(x,y)$ introduced in (\ref{kernel}). Consider the map
%
%e39 #&#
\begin{equation}
\label{fN} f_N(z)=-\frak b_{N}(z-\frak
d_N)+\log(z)-\frac{1}{N}\sum_{i=1}^n
\log (1-\lambda_i z).
\end{equation}

%re4.5 #&#
\begin{Remark}
In order to fully define $f_N$, one needs to specify the determination
of the logarithm. This will be done when needed. Notice, however, that
functions $\re f_N$, $\exp(f_N)$ and the derivatives $f_N^{(k)}$ are
always well defined.
\end{Remark}

By taking $q=\frak d_N$ in (\ref{kernel}), which is possible according
to Remark~\ref{freedomq}, we have
\begin{eqnarray}
\label{KNversion}
\K_{N}(x,y)
&=& \frac{ N}{(2i\pi)^2}\oint_{\Gamma}\d z\oint_{\Theta
} \d w\, e^{- Nx(z-\frak d_N) +Ny(w-\frak d_N)}
\nonumber\\[-8pt]\\[-8pt]\nonumber
&&\hspace*{71pt}{}\times
\frac{1}{w-z} \biggl(\frac{z}{w}
\biggr)^{ N}\prod_{j=1}^n
\biggl(\frac{1-\lambda
_jw}{1-\lambda_jz} \biggr),
\nonumber
\end{eqnarray}
where we recall that the contour $\Gamma$ encloses the $\lambda
_j^{-1}$'s whereas the contour $\Theta$ encloses the origin and is
disjoint from $\Gamma$.
It then follows from definition (\ref{rescaledkernel}) of $\widetilde
\K_N$ that
%
%e40 #&#
\begin{eqnarray}
\label{kernelright}
\widetilde\K_{N}(x,y)
&=&\frac{N^{1/3}}{(2i\pi)^2\delta_N}\nonumber
\\
&&{}\times \oint _{\Gamma} \d z \oint_{\Theta} \d w
\frac{1}{w-z}
\\
&&\hspace*{49pt}{} \times  e^{-
N^{1/3}x\sklfrac{(z-\frak d_N)}{\delta_N}
+N^{1/3}y\sklfrac{(w-\frak d_N)}{\delta_N}+N f_N(z)-N f_N(w)}.
\nonumber
\end{eqnarray}
The key observation here is the identity
%
%e41 #&#
\begin{equation}
\label{keyidentity} f'_N(z)=g_N(z)-g_N(
\frak d_N),
\end{equation}
which follows from (\ref{gN}) and (\ref{notations}). As a byproduct,
(\ref{infocN}) yields that $\frak d_N$ is a root of multiplicity two
for $ f'_N$, and more precisely,
%
%e42 #&#
\begin{equation}
\label{csqkeyidentity} f'_N(\frak d_N)=f''(
\frak d_N)=0,\qquad f^{(3)}_N(\frak
d_N)=g''_N(\frak
d_N)>0.
\end{equation}

The aim is to perform a saddle point analysis for $f_N$ around its
critical point~$\frak d_N$. To this end, we deform the contours $\Gamma$ and
$\Theta$ in a way that they pass near $\frak d_N$.

If\vspace*{-1pt} $\frak d_N$ is smaller than all the $\lambda_j^{-1}$'s, as it is
the case in \cite{elkaroui-07-TW} when dealing with the maximal eigenvalue,
then go directly to Section~\ref{contourssection}, set $\Gamma
^{(1)}=\Gamma$, $\K_N^{(1)}= \widetilde\K_N$ and disregard every
statement related to $\Gamma^{(0)}$.

%
%f5 #&#
\begin{figure}

\includegraphics{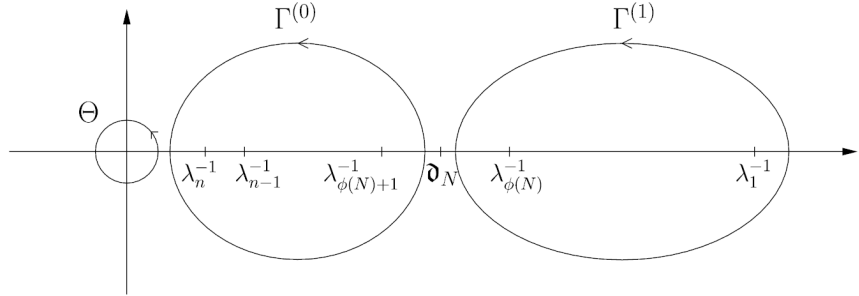}

\caption{The new contours $\Gamma^{(0)}$ and $\Gamma^{(1)}$.}
\label{figdeform-gamma}
\end{figure}

%
%f6 #&#
\begin{figure}[b]

\includegraphics{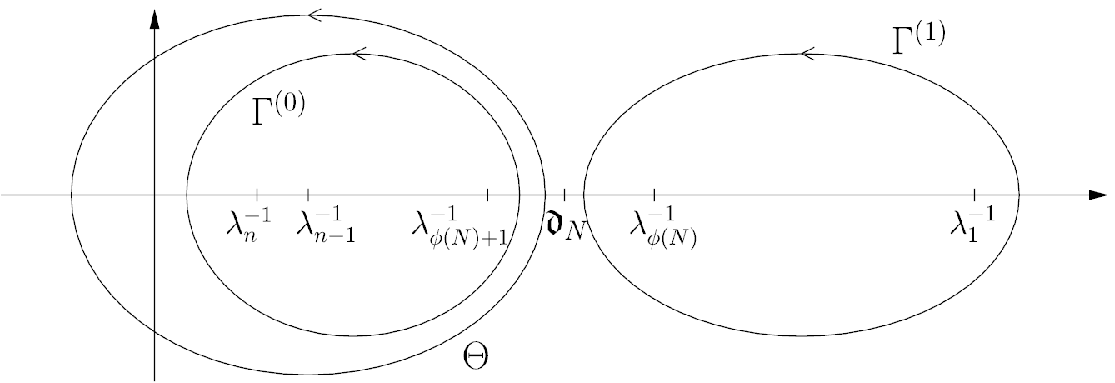}

\caption{The new contours for the kernel $\widetilde\K_N$.}
\label{fignewcontours}
\end{figure}

If not, then we proceed in two steps. First, we split $\Gamma$ into two
disjoint contours, $\Gamma^{(0)}$ and $\Gamma^{(1)}$, as shown in
Figure~\ref{figdeform-gamma}: the contour $\Gamma^{(0)}$ encloses the
$\lambda_j^{-1}$'s which are smaller than $\frak d_N$, while $\Gamma^{(1)}$
encloses the $\lambda_j^{-1}$'s which are larger that $\frak d_N$.
Notice that Proposition~\ref{1111}(d) applied
to the measure $\nu_N$ shows that the set $\{ j, 1\le j\le n\dvtx
\lambda_j^{-1} > \frak d_N \}$ is not empty. Therefore, the contour
$\Gamma^{(1)}$ is always well defined.\vadjust{\goodbreak} %\vspace*{1pt}

We now introduce for $\alpha\in\{0,1\}$ the kernels
%
%e43 #&#
\begin{eqnarray}
\label{kernelsplitdef}
&& \K_{N}^{(\alpha)}(x,y)\nonumber
\\
&&\qquad = \frac{N^{1/3}}{(2i\pi)^2\delta_N}
\nonumber\\[-8pt]\\[-8pt]\nonumber
&&\quad\qquad{}\times \oint _{\Gamma^{(\alpha)}} \d z \oint_{\Theta} \d w
\frac{1}{w-z}
\\
&&\hspace*{33pt}\hspace*{59pt}{} \times e^{-N^{1/3}x\sklfrac{(z-\frak d_N)}{\delta_N}
+N^{1/3}y\sklfrac{(w-\frak d_N)}{\delta_N}+N f_N(z)-N f_N(w)}.\nonumber
\end{eqnarray}
Then it follows from the residue theorem that
%
%e44 #&#
\begin{equation}
\label{split} \widetilde\K_N(x,y)= \K_N^{(0)}(x,y)+
\K_N^{(1)}(x,y),
\end{equation}
and a similar identity for the associated operators.

In the second step, we modify the contour $\Theta$ in order for it to surround
$\Gamma^{(0)}$ while remaining at the left of $\frak d_N$; cf.
Figure~\ref{fignewcontours}. This can be done with no harm for the kernel
$\K_N^{(1)}$. As for $\K_N^{(0)}$, this modification for the contours
yields a residue term, coming with the singularity $(w-z)^{-1}$ of the
integrand. The latter residue term equals
\[
\frac{N^{1/3}}{2i\pi\delta_N}\oint_{\Gamma^{(0)}}e^{N^{1/3}\sklfrac
{(y-x)}{\delta_N}(z-\frak d_N)}\,\d z
\]
and thus identically vanishes since the integrand is analytic.

%s4.4 #&#
\subsection{Contours deformations and subharmonic functions: The right
edge case}
\label{contourssection}

We now provide the existence of deformations for the contours $\Gamma
^{(0)}$, $\Gamma^{(1)}$ and $\Theta$~which are appropriate\vspace*{1pt} for the
asymptotic analysis. These new contours will be referred to as
$\Upsilon^{(0)}$, $\Upsilon^{(1)}$ and $\widetilde\Theta$.

%
%pr4.6 #&#
\begin{proposition}
\label{contoursprop}
For every $\rho>0$ small enough, there exists a contour $\Upsilon
^{(0)}$ independent
of $N$ and two contours $\Upsilon^{(1)}=\Upsilon^{(1)}(N)$ and
$\Ttilde=\Ttilde(N)$ which satisfy
for every $N$ large enough the following properties:
\begin{longlist}[(1)]
\item[(1)]
\begin{enumerate}[(a)]
\item[(a)]
$\Upsilon^{(0)}$ encircles the $\lambda_j^{-1}$'s smaller than $\frak d_N$;
\item[\hspace*{31pt}(b)]
$\Upsilon^{(1)}$ encircles all the $\lambda_j^{-1}$'s larger than
$\frak d_N$;
\item[\hspace*{31pt}(c)]
$\Ttilde$ encircles all the $\lambda_j^{-1}$'s smaller than $\frak
d_N$ and the origin.
\end{enumerate}

\item[(2)]
\begin{enumerate}[(a)]
\item[(a)] $\Upsilon^{(1)}=\Upsilon_*\cup\Ures^{(1)}$ where
\[
\Upsilon_*=\bigl\{\frak d_N+te^{\pm i\pi/3}\dvtx  t\in[0,\rho]\bigr\};
\]
\item[\hspace*{31pt}(b)] $\Ttilde=\Ttilde_*\cup\Ttilde_{\mathrm{res}}$ where
\[
\Ttilde_*=\bigl\{\frak d_N-te^{\pm i\pi/3}\dvtx  t\in[0,\rho]\bigr\}.
\]
\end{enumerate}
\item[(3)]
There exists $K>0$ independent of $N$ such that:
\begin{enumerate}[(a)]
\item[\hspace*{31pt}(a)] $\re ( f_N(z)-f_N(\frak d_N) )\leq-K$ for
all $z\in\Upsilon^{(0)}$;\vspace*{1pt}
\item[\hspace*{31pt}(b)] $\re ( f_N(z)-f_N(\frak d_N) )\leq-K$ for
all $z\in\Ures^{(1)}$;\vspace*{1pt}
\item[\hspace*{31pt}(c)] $\re ( f_N(w)-f_N(\frak d_N) )\geq K$ for
all $w\in\Ttilde_{\mathrm{res}}$.
\end{enumerate}

\item[(4)]
There exists $\d>0$ independent of $N$ such that
\begin{eqnarray*}
\inf \bigl\{\llvert z-w\rrvert\dvtx  z\in\Upsilon^{(0)}, w\in\Ttilde \bigr
\} &\geq & \d,
\\
\inf \bigl\{\llvert z-w\rrvert\dvtx  z\in\Upsilon_{*}, w\in
\Ttilde_{\mathrm{res}} \bigr\} &\geq& \d,
\\
\inf \bigl\{\llvert z-w\rrvert\dvtx  z\in\Ures^{(1)}, w\in
\Ttilde_{*} \bigr\} &\geq & \d,
\\
\inf \bigl\{\llvert z-w\rrvert\dvtx  z\in\Ures^{(1)}, w\in
\Ttilde_{\mathrm{res}} \bigr\} &\geq& \d.
\end{eqnarray*}

\item[(5)]
\begin{enumerate}[(a)]
\item[(a)]
The contours $\Upsilon^{(1)}$ and $\Ttilde$ lie in a bounded subset
of $\C$ independent \hspace*{31pt}of~$N$;
\item[\hspace*{31pt}(b)]
the lengths of $\Upsilon^{(1)}$ and $\Ttilde$ are uniformly bounded
in $N$.
\end{enumerate}
\end{longlist}
\end{proposition}

Note that both the contours $\Upsilon^{(1)}$ and $\Ttilde$ pass
through the
critical point $\frak d_N$.

In order to provide a proof for Proposition~\ref{contoursprop}, we
first establish a few lemmas. We recall that $B(z,\rho)$ for $z\in\C
$, and $\rho> 0$ stands for the open ball of $\C$ with center $z$ and
radius $\rho$.

Recall that $0<\inf_N \lambda_n^{-1}\leq\sup_N\lambda
_1^{-1}<\infty$ by Assumption~\ref{assnu}. By the regularity
assumption, namely
$\liminf_N \min_{j=1}^n \llvert   \frak d - \lambda_j^{-1} \rrvert  >0$, there
exists $\varepsilon>0$ such that $\lambda_j^{-1}\in(0,+\infty
)\setminus B(\frak d,\varepsilon)$ for every $1\leq j\leq n$ and every
$N$ large enough. Denote by ${\mathcal K}$ the compact set
%
%e45 #&#
\begin{equation}
\label{defK} {\mathcal K} = \biggl( \biggl[\inf_N
\frac{1}{\lambda_n},\sup_N\frac{1}{\lambda_1} \biggr] \Bigm\backslash
B(\frak d,\varepsilon) \biggr) \cup\{ 0\}.
\end{equation}
Notice that by construction $\{x\in\R\dvtx   x^{-1}\in\supp(\nu_N)\}
\subset\mathcal K$ for every $N$ large enough, and also that $\{x\in
\R\dvtx   x^{-1}\in\supp(\nu)\}\subset\mathcal K$ because of the weak
convergence $\nu_N\to\nu$.

Recall the definition of $f_N$ (\ref{fN}), and introduce its
asymptotic counterpart,
%
%e46 #&#
\begin{equation}
\label{f} f(z)=-\frak b(z- \frak d)+\log(z)-\gamma\int\log(1-xz)\nu(\d x).
\end{equation}
Notice that whereas $f$ and $f_N$ are defined up to a determination of
the complex logarithm,
%
%e47 #&#
\begin{equation}
\label{defRef} \re f(z)=-\frak b\re(z-\frak d)+\log\llvert z\rrvert -\gamma\int
\log\llvert 1-xz\rrvert \nu (\d x)
\end{equation}
and $\re f_N$ are well defined. The following properties of $\re f$ and
$\re f_N$ around $\frak d$ and $\frak d_N$ will be of constant use in the
sequel.

% WH: la notation $B(z_0,\rho)$ a deja ete introduite.
% \paragraph*{???????????????Notation:} From now, for any $z_0\in\C$ and $\rho>0$, we
%use the notation
% \[
% B(z_0,\rho)=\big\{z\in\C\dvtx    \left\vert z-z_0\right\vert < \rho\big\}.
% \]

%{\red[jam\dvtx  Dans le lemme suivant, j'ai int{\chr"C3}{\chr"A9}gr{
%\chr"C3}{\chr"A9} le lemme~\ref{unifconvre}, le lemme~\ref{controlrho}
%qui {\chr"C3}{\chr"A9}tait dans la section relative {\chr"C3}{\chr"A0}
%$K_N^{(0)}$ mais qui servait ailleurs aussi, et une proposition
%relative {\chr"C3}{\chr"A0} $f$ qu'on utilise mais qui n'{\chr"C3}{
%\chr"A9}tait pas packag{\chr"C3}{\chr"A9}e. n.b. dans la version pr{
%\chr"C3}{\chr"A9}c{\chr"C3}{\chr"A9}dente, la d{\chr"C3}{
%\chr"A9}finition de $f$ n'apparaissait que dans une preuve. Je l'ai
%donc d{\chr"C3}{\chr"A9}finie ci-dessus ainsi que ${\mathcal K}$ dont
%la d{\chr"C3}{\chr"A9}finition intiale {\chr"C3}{\chr"A9}tait un peu
%fausse (ci-dessous)]}

%
%le4.7 #&#
\begin{lemma}\label{lemmaproperties-f-fN}
Let Assumption~\ref{assnu} hold true, and let ${\mathcal K}$ be as in
(\ref{defK}). Then:
\begin{longlist}[(a)]
\item[(a)] The function $\re f_N$ converges locally uniformly to
$\re f$ on $\C\setminus\mathcal K$. Moreover,
%
%e48 #&#
\begin{equation}
\label{reconv} \lim_{N\rightarrow\infty}\re f_N(\frak
d_N)=\re f(\frak d).
\end{equation}
\item[(b)] There exists $\rho_0>0$ and $\Delta=\Delta(\rho
_0)>0$ independent of $N$ such that for every $N$ large enough,
$B(\frak d_N,\rho) \subset\C\setminus{\mathcal K}$ for every $ \rho
\in(0,\rho_0]$, and whatever the analytic representation of $f_N$ on
$B(\frak d_N,\rho)$,
\begin{eqnarray*}
\bigl\llvert f_N(z)-f_N(\frak d_N)-g''_N(
\frak d_N) (z-\frak d_N)^3/6\bigr\rrvert &
\leq& \Delta\llvert z-\frak d_N\rrvert ^4,
\\
\bigl\llvert \re\bigl(f_N(z)-f_N(\frak d_N)
\bigr)-g''_N(\frak d_N)\re
\bigl[(z-\frak d_N)^3\bigr]/6\bigr\rrvert &\leq& \Delta
\llvert z-\frak d_N\rrvert ^4
\end{eqnarray*}
for all $ z\in B(\frak d_N,\rho_0)$.
\item[(c)]
There exists $\rho_0>0$ and $\Delta=\Delta(\rho_0)>0$ such that
$B(\frak d,\rho_0) \subset\C\setminus{\mathcal K}$, and
for all $z\in B(\frak d,\rho_0)$,
\[
\bigl\llvert \re\bigl(f(z)-f(\frak d)\bigr)-g''(
\frak d)\re\bigl[(z-\frak d)^3\bigr]/6\bigr\rrvert \leq\Delta\llvert
z-\frak d\rrvert ^4.
\]
\end{longlist}
\end{lemma}

\begin{pf}
Fix an open ball $B$ of $\C\setminus\mathcal K$. By definition of
$\mathcal K$, one can chose a determination of the logarithm such that $f_N$
is well defined and holomorphic there for $N$ large enough. Indeed,
there exists an analytic determination of the logarithm on every simply
connected domain of $\C\setminus\{0\}$. Use the same determination
for $f$, which is then also well defined and holomorphic on $B$.
By weak convergence of $\nu_N$ to $\nu$, $f_N$ converges pointwise to
$f$ on $B$. Similar to the
proof of Proposition~\ref{gN->g}, the sequence of holomorphic functions
$(f_N)_N$ is uniformly bounded on $B$ and thus has compact closure by the
Montel theorem, which upgrades the pointwise convergence
$f_N\rightarrow f$ to
the uniform one on $B$. The uniform convergence of $\re f_N$ to $\re f$
on $B$
follows since $\llvert  \re f_N(z)-\re f(z)\rrvert  \leq\llvert  f_N(z)-f(z)\rrvert  $ for all $z\in
B$. Now since $\frak d_N\rightarrow\frak d$ and
$\frak d_N,\frak d\in\C\setminus\mathcal K$ for all $N$ large enough
by the
regularity assumption, (\ref{reconv}) follows from the local uniform
convergence $\re f_N\rightarrow\re f$ on $\C\setminus\mathcal K$,
and (a) is proved.

% For $\eta> 0$ small enough, the function $f$ introduced in~\eqref{f}
%is well
%defined and holomorphic on $B(\frak c,\eta)$, see Remark
%\ref{logarithm}. Note that since $\frak b=g(\frak c)$, the definition
%of $f$ provides that $f'(z)=g(z)-g(\frak c)$. Thus, \eqref{critg} and
%a Taylor development for $f$ around $\frak c$ yield for every $z\in B(
%\frak c,\eta)$
%\begin{eqnarray*}
%\left\vert \re f(z)- \re f(\frak c) -g''(\frak c)\re(z-\frak c)^3/6\right\vert  &
%\leq& \left\vert  f(z)- f(\frak c) -g''(\frak c)(z-\frak c)^3/6\right\vert \\
%&\leq& \frac{\left\vert z-\frak c \right\vert ^4}{24}\max_{w\in B(\frak c,\eta)}\left\vert f^{(4)}(w)\right\vert.
%\end{eqnarray*}
%Since $g''(\frak c)>0$, the lemma follows by choosing $\eta$ small
%enough.

It follows from Proposition~\ref{gN->g} that for $\rho_0>0$ small
enough and
every $N$ large enough, we have
$B(\frak d_N,\rho_0)\subset B(\frak d,2\rho_0)\subset\C\setminus
\mathcal K$.
Using the same determination of the $\log$ as previously yields that
$f_N$ is
well defined and holomorphic on $B(\frak d_N,\rho_0)$. Since~(\ref{infocN})
and~(\ref{keyidentity}) yield $f_N'(\frak d_N)=f_N''(\frak d_N)=0$,
$ f_N^{(3)}(\frak d_N)=g_N''(\frak d_N)>0$ and $f_N^{(4)}=g_N^{(3)}$
for all
$N$ large enough, we can perform a Taylor expansion for $f_N$ around
$\frak d_N$ in order to get
\[
\bigl\llvert f_N(z)-f_N(\frak d_N)-g''_N(
\frak d_N) (z-\frak d_N)^3/6\bigr\rrvert \leq
\frac{\llvert  z-\frak d_N\rrvert  ^4}{24}\max_{w\in B(\frak d,2\rho_0)} \bigl\llvert g_N^{(3)}(w)
\bigr\rrvert
\]
provided that $z\in B(\frak d_N,\rho_0)$. Proposition~\ref{gN->g} moreover
provides that $g_N^{(3)}$ converges uniformly on
$B(\frak d,2\rho_0)$ to $g^{(3)}$ which is bounded there. We therefore
get the
existence of $\Delta=\Delta(\rho_0)$ independent of $N$ for which
the first
inequality in part~(b) of the proposition is satisfied. The inequality
for the
real part directly follows, and part~(b) of the proposition is proved,
as is
part~(c) by using similar arguments.
\end{pf}

We now provide a qualitative analysis for the map $\re f$. First, we
study the behavior of $\re f(z)$ as $\llvert  z\rrvert  \to\infty$. To do so, we
introduce the sets
\begin{eqnarray*}
\Omega_- & =& \bigl\{z \in\C\dvtx  \re f(z)<\re f(\frak d) \bigr\},
\\
\Omega_+& =& \bigl\{z \in\C\dvtx  \re f(z)>\re f(\frak d) \bigr\},
\end{eqnarray*}
and prove the following.

%
%le4.8 #&#
\begin{lemma}
\label{behaviorinfinity}
Both $\Omega_+$ and $\Omega_-$ have a unique unbounded connected
component. Moreover, given any $\alpha\in(0,\pi/2)$, there exists
$R>0$ large enough such that
%
%e49 #&#
%e50 #&#
\begin{eqnarray}
\label{omegaR-} \Omega_-^R &=& \biggl\{z\in\C\dvtx  \llvert z\rrvert >R, -
\frac{\pi}{2}+\alpha<\arg (z)<\frac{\pi}{2}-\alpha \biggr\}\subset
\Omega_-,
\\
\label{omegaR+} \Omega_+^R &=& \biggl\{z\in\C\dvtx  \llvert z\rrvert >R,
\frac{\pi}{2}+\alpha<\arg (z)<\frac{3\pi}{2}-\alpha \biggr\}\subset
\Omega_+.
\end{eqnarray}
\end{lemma}

\begin{pf}
Recall expression (\ref{defRef}) of $\re f(z)$ which yields that $\re
f(z)=-\frak b \re(z-\frak d)+O(\log\llvert  z\rrvert  )$ as $\llvert  z\rrvert  \rightarrow\infty$.
Since $\frak b>0$, it follows that for any fixed $\alpha\in(0,\pi
/2)$, there exists $R>0$ large enough such that
%
%e51 #&#
\begin{equation}
\label{rightminus} \Omega_-^R \subset \Omega_-,\qquad
\Omega_+^R \subset \Omega_+.
\end{equation}
Next, we compute for any $A\in\R\setminus\{0\}$,
\[
\frac{\d}{\d t}\re f(t+iA)=-\frak b +\frac{t}{t^2+A^2}+\gamma\int
\frac{(x^{-1}-t)}{(x^{-1}-t)^2+A^2} \nu(\d x).
\]
Since $\frak b>0$ and $\supp(\nu)$ is a compact subset of $(0,+\infty
)$, there exists $A_0>0$ such that for any $A$ satisfying $\llvert  A\rrvert  \geq
A_0$, the map $t\mapsto\frac{\d}{\d t}\re f(t+iA)$ is negative;
namely, $t\mapsto\re f(t+iA)$ is decreasing. Assume there exists
another unbounded connected component of $\Omega_-$, different from
the one containing $\Omega_-^R$. By~(\ref{rightminus}). This
unbounded connected component then lies in $\C\setminus(\Omega
_-^R\cup\Omega_+^R)$, and thus there exists $z_0$ in this component
satisfying $\llvert  \operatorname{Im}(z_0)\rrvert  \geq A_0$. Since the half line $\{\operatorname{Re}(z_0)+t+i\operatorname{Im}(z_0)\dvtx   t\geq0\}$ then belongs to $\Omega_-$ and
eventually hits $\Omega_-^R$, we obtain a contradiction. The same
arguments apply to $\Omega_+$.
%That $\Omega_+$ has a unique unbounded connected component follows by
%using the same line of arguments.
\end{pf}

Next, we describe the behavior of $\re f$ at the neighborhood of $\frak
d$. Taking advantage of Lemma~\ref{lemmaproperties-f-fN}(c), which
encodes that $\re f(z)-\re f(\frak d)$ behaves like $\re[(z-\frak
d)^3]$ around $\frak d$, we describe in the following lemma subdomains
of $\Omega_\pm$ of interest.

%
%le4.9 #&#
\begin{lemma}
\label{trianglef}
There exist $\eta>0$ and $\theta>0$ small enough such that, if
\[
\Delta_k= \biggl\{z\in\C\dvtx  0<\llvert z-\frak d\rrvert <\eta, \biggl
\llvert \arg(z-\frak d)-k\frac{\pi}{3}\biggr\rrvert <\theta \biggr\}
\]
for $-2\le k\le3$, then
\[
\Delta_{2k+1}\subset\Omega_-,\qquad
\Delta_{2k}\subset\Omega_+,\qquad
k\in\{-1,0,1\}. %
\]
\end{lemma}

The regions $\Delta_k$ are shown on Figure~\ref{figptselle-right}.

\begin{pf*}{Proof of Lemma \ref{trianglef}}
Recall Lemma~\ref{lemmaproperties-f-fN}(c), and let $\eta<\rho_0$ as
defined there. Then
\[
\bigl\llvert \re f(z)- \re f(\frak d) -g''(\frak d)
\re\bigl[(z-\frak d)^3\bigr]/6 \bigr\rrvert \le\Delta(
\rho_0)\llvert z - \frak d\rrvert ^4 %
\]
for every $z\in B(\frak d,\eta)$. Notice that $\re[(z-\frak
d)^3]=(-1)^k$ if $z=\frak d + e^{i k\pi/3}$ for consecutive integers $k$.
Since $g''(\frak d)>0$, the lemma follows by choosing $\eta$ small enough.
\end{pf*}

We denote by $\Omega_{2k+1}$ the connected component of $\Omega_-$
which contains $\Delta_{2k+1}$. Similarly, $\Omega_{2k}$ stands for
the connected component of $\Omega_+$ which contains $\Delta_{2k}$.
We now describe these sets by using the maximum principle for
subharmonic functions, in the same spirit as in \cite{delvaux-kuijlaars-10}, Section~6.1 (see also %{\red
\cite{delvaux-kuijlaars-2009}, Section~2.4.2),
%\{il faudrait rajouter \`a la biblio: S. Delvaux and A. B. J.
%Kuijlaars, \emph{A phase transition for nonintersecting Brownian
%motions, and the Painlev\'e II equation},
%Int. Math. Res. Not. 2009 (2009), 3639--3725. \}},
although the setting is more involved here; such a use of the maximum
principle has been communicated to us by Steven Delvaux.

Recall that if $G$ is an open subset of $\C$, a function $u\dvtx G\to\R
\cup\{-\infty\}$ is \textit{subharmonic} if $u$ is upper semicontinuous;
that is, $\{ z\in G,  u(z)< \alpha\}$ is open for every $\alpha
\in\R$, and for every closed disk $\overline{B}(z,\delta)$
contained in $G$, we have the inequality
\[
u(z) \leq \frac{1}{2\pi}\int_0^{2\pi}u
\bigl(z+\delta e^{i\theta}\bigr)\,\d \theta.
\]
A function $u\dvtx G\to\R\cup\{+\infty\}$ is \textit{superharmonic} if $-u$
is subharmonic; in particular, it is lower semicontinuous.
Moreover, if $u\dvtx G\to\C$ is subharmonic, it satisfies a maximum
principle: for any bounded domain (i.e., connected open set) $U\subset
\C$ where $u$ is subharmonic, if for some $\kappa\in\R$ it holds that
\[
\limsup_{z\to\zeta,  z\in U}u(z) \leq \kappa,\qquad\zeta \in\partial U,
\]
then $u\leq\kappa$ on $U$. Similarly, superharmonic functions satisfy
a minimum principle.

The use of the maximum principle for subharmonic functions is made
possible here because of the following observation.

%le4.10 #&#
\begin{lemma}
\label{supersubharmonic}
The function $\re f$ is subharmonic on
$\C\setminus\{x\in\R\dvtx   x^{-1}\in\supp(\nu)\}$ and superharmonic on
$\C\setminus\{0\}$.
\end{lemma}

\begin{pf}
It will be enough to establish the result for the map
\begin{eqnarray}
z&\mapsto&\log\llvert z\rrvert - \gamma\int\log\llvert 1 - xz \rrvert \nu(\d x)
\nonumber\\[-8pt]\\[-8pt]\nonumber
&=& \log\llvert z\rrvert - \gamma\int\log\llvert z - x \rrvert \tau(\d x) - \gamma
\int\log x \nu(\d x),
\nonumber
\end{eqnarray}
where the compactly supported probability measure $\tau$ is the image
of $\nu$
by $x\mapsto x^{-1}$. The assumptions on $\nu$ imply that
$\log x$ is $\nu$-integrable.
Now, it is a standard fact from potential theory that given a positive
Borel measure $\eta$ on $\C$ with compact support, the map
$z\mapsto\int\log\llvert  z-x\rrvert    \eta(\d x)$ is subharmonic on $\mathbb C$ and
harmonic on $\mathbb C \setminus\supp(\eta)$; see, for
example, \cite{book-saff-totik-1997}, Chapter~0.
Consequently, $z\mapsto\log\llvert  z\rrvert  $ is harmonic on $\C\setminus\{0\}$ and
subharmonic on $\C$, and $z\mapsto\gamma\int\log\llvert   z - x \rrvert     \tau
(\d x)$ is
harmonic on $\C\setminus\supp(\tau)$ and subharmonic on $\C$. The result
follows.
\end{pf}

Equipped with Lemma~\ref{supersubharmonic}, we can obtain more
information concerning the connected components of $\Omega_\pm$.

%
%le4.11 #&#
\begin{lemma}\label{keylevelsets}
The following hold true:
\begin{enumerate}[(2)]
\item[(1)]
If $\Omega_*$ is a connected component of $\Omega_+$, then $\Omega
_*$ is open and, if $\Omega_*$ is moreover bounded, there exists $x\in
\supp(\nu)$ such that $x^{-1}\in\Omega_*$.\vadjust{\goodbreak}
\item[(2)]
Let $\Omega_*$ be a connected component of $\Omega_-$ with nonempty interior:
\begin{longlist}[(a)]
\item[(a)]
if $\Omega_*$ is bounded, then $0\in\Omega_*$;
\item[(b)]
if $\Omega_*$ is bounded, then its interior is connected;
\item[(b)]
if $0\notin\Omega_*$, then the interior of $\Omega_*$ is connected.
\end{longlist}
\end{enumerate}
\end{lemma}

\begin{pf}
Let us show (1). We set $\alpha=\re f(\frak d)$.
Since $\re f(z)\to-\infty$ as $\llvert  z\rrvert  \to0$, then $0\notin\{ z\in
\mathbb{C}, \re f > \alpha\}$.
Hence
\[
\{ z\in\mathbb{C}\dvtx  \re f> \alpha\} =\bigl\{ z\in\C\setminus\{0\}\dvtx  \re f>\alpha
\bigr\}. %
\]
However, since $\re f$ is superharmonic on $\C\setminus\{0\}$, $\{
z\in\C\setminus\{ 0\}\dvtx   \re f>\alpha\}$ is an open set on $\C$.
As a consequence, all these connected components are open, hence the
desired result. In particular, $\Omega_*$ is open and $\partial\Omega
_*\subset\partial\Omega_+$; hence $\re f\leq\re f(\frak d)$ on
$\partial\Omega_*$. If $\Omega_*$ is moreover bounded, then we have
$\re f>\re f(\frak d)$ on the bounded domain $\Omega_*$ and $\re f\leq
\re f(\frak d)$ on its boundary. Since subharmonic functions satisfy a
maximum principle, $\re f$ cannot be subharmonic on the whole set
$\Omega_*$, and (1) follows from Lemma~\ref{supersubharmonic}.

We now turn to (2)(a). We argue by contradiction and assume that $\Omega_*$
%\not\subset\R$
is a bounded connected component of $\Omega_-$ which does not contain
the origin. The fact that $\Omega_*$ has a nonempty interior
%The assumption $\Omega_*\not\subset\R$
implies that at least one of the sets $\Omega_*\cap\{\operatorname{Im}(z)>0\}$
or $\Omega_*\cap\{\operatorname{Im}(z)<0\}$ is nonempty. Consider the set
\[
\Omega_*^{\mathrm{sym}}= \{z\in\C\dvtx  \overline z\in\Omega_*\},
\]
and notice it is also a connected component of $\Omega_-$ because of
the symmetry $\re f(z)=\re f(\overline z)$. Without loss of generality,
assume that $\Omega_*\cap\{\operatorname{Im}(z)>0\}\neq\varnothing$ (otherwise
switch the role of $\Omega_*$ and $\Omega_*^{\mathrm{sym}}$ in what follows).
%\textcolor{blue}{[{\mathbf old}] Notice that since $\re f$ is continuous
%on $ \C\setminus\mathcal K $ by Lemma~\ref{supersubharmonic}, the set $
%\Omega_*\cap\{\operatorname{Im}(z)>0\}$ is open and moreover
%\begin{equation}
%\label{valueborder}
%\re f(z)=\re f(\frak d),\qquad z\in\partial\Omega_*\setminus\mathcal
%K.
%\end{equation} }
Since $\re f$ is subharmonic on $\C\setminus\mathcal K$, $\Omega_-$
is open and so are its connected components, in particular $\Omega_*$,
and then $\Omega_*\cap\{\operatorname{Im}(z)>0\}$. Now $\re f$ being
continuous on $ \C\setminus\mathcal K $, by Lemma~\ref
{supersubharmonic}, we have
%
%e52 #&#
\begin{equation}
\label{valueborder} \re f(z)=\re f(\frak d),\qquad z\in\partial\Omega_*\setminus
\mathcal K.
\end{equation}
Let us fix $\varepsilon_0>0$ such that $\Omega_*\cap\{\operatorname{Im}(z)\geq\varepsilon_0\}\neq\varnothing$ and pick $z_0\in\Omega_*$
satisfying $\operatorname{Im}(z_0)\geq\varepsilon_0$ and $\re f(z_0)<\re
f(\frak d)$. Our goal is to construct a bounded domain which contains
$z_0$ but not the origin and where $\re f>\re f(z_0)$ on its boundary.
Indeed, this would lead to a contradiction via the minimum principle
for superharmonic functions since $\re f$ is superharmonic on $\C
\setminus\{0\}$ as stated in Lemma~\ref{supersubharmonic}.

First, notice that if $\operatorname{dist}(\Omega_*,\R)>0$, then $\re f$ is
harmonic on $\Omega_*$, $\re f=\re f(\frak d)$ on $\partial\Omega_*$
and $\re f<\re f(\frak d)$ on $\Omega_*$, which is a bounded domain.
However, this contradicts the minimum principle for (super)harmonic
functions, and thus $\operatorname{dist}(\Omega_*,\R)=0$. Because $\operatorname{dist}(\Omega_*,\R)=0$ and $\Omega_*\cap\{\operatorname{Im}(z)>0\}$ is open
and nonempty, for every $\varepsilon>0$ small enough $\Omega_*\cap\{
\operatorname{Im}(z)=\varepsilon\}=U+i\varepsilon$ where $U$ is a nonempty
open subset of the real line. Thus we can write
\[
\Omega_*\cap\bigl\{\operatorname{Im}(z)=\varepsilon\bigr\} = \bigcup
_{j\in J} \bigl(u_{\min}^{(j)}(
\varepsilon),u_{\max}^{(j)}(\varepsilon) \bigr)+i\varepsilon,
\]
where\vspace*{1pt} $J$ is a countable set satisfying $\operatorname{Card}(J)\geq1$, and the
$u_{\min}^{(j)}(\varepsilon)$'s and $u_{\min}^{(j)}(\varepsilon)$'s
are\vspace*{2pt} real numbers such that any open intervals $(u_{\min
}^{(j_1)}(\varepsilon),u_{\max}^{(j_1)}(\varepsilon))$ and $(u_{\min
}^{(j_2)}(\varepsilon),u_{\max}^{(j_2)}(\varepsilon))$ are disjoint
whenever $j_1\neq j_2$. Notice that by symmetry,
\[
\Omega_*^{\mathrm{sym}}\cap\bigl\{\operatorname{Im}(z) = -\varepsilon\bigr\}= \bigcup
_{j\in J} \bigl(u_{\min}^{(j)}(
\varepsilon),u_{\max
}^{(j)}(\varepsilon) \bigr)-i\varepsilon.
\]
By construction, for every $j\in J$, both $u_{\min}^{(j)}(\varepsilon
)+i\varepsilon$ and $u_{\max}^{(j)}(\varepsilon)+i\varepsilon$
belong to $\partial\Omega_*\setminus\R$. In particular, by (\ref
{valueborder}) and the symmetry $\re f(z)=\re f(\overline z)$,
%
%e53 #&#
\begin{equation}
\label{upbordervalue} \re f \bigl(u_{\min}^{(j)}(\varepsilon)\pm i
\varepsilon \bigr) = \re f \bigl(u_{\max}^{(j)}(\varepsilon)\pm
i\varepsilon \bigr) = \re f(\frak d),\qquad j\in J.
\end{equation}
Since by assumption $0\notin\Omega_*\subset\Omega_-$, there exists
$\delta>0$ such that $B(0,\delta)\cap\Omega_*=\varnothing$
%\textcolor{red}{[c'est certainement vrai mais incomplet ou trop
%rapide]}
otherwise $0\in\partial\Omega_*$, but in this case, the boundary
condition $\re f(z) =\re f(\frak d)$ would be violated near zero as
$\re f(z)\to-\infty$ for $\llvert  z\rrvert  \to0$.
As $\Omega_*$ is moreover bounded by assumption, $\llvert  u_{\min
}^{(j)}(\varepsilon)\rrvert  $ and $\llvert  u^{(j)}_{\max}(\varepsilon)\rrvert  $ stay in a
compact subset of $(0,+\infty)$ independent from $\varepsilon$ and
$j\in J$ as $\varepsilon\to0$.
As a consequence, we can choose $\varepsilon\in(0,\varepsilon_0)$
small enough so that, for every $j\in J$,
%
%e54 #&#
\begin{equation}
\label{assumptionepsilon} \max \biggl( \frac{\varepsilon^2}{u_{\min}^{(j)}(\varepsilon)^2}, \frac{\varepsilon^2}{u_{\max}^{(j)}(\varepsilon)^2} \biggr) <
\min \biggl(\re f(\frak d)-\re f(z_0),\frac{1}2 \biggr).
\end{equation}
If we moreover consider for any $j\in J$ the open rectangle
\[
\mathcal R_j(\varepsilon)= \bigl\{ u+iv\in\C\dvtx  u_{\min
}^{(j)}(
\varepsilon)<u<u_{\max}^{(j)}(\varepsilon), \llvert v\rrvert <
\varepsilon \bigr\},
\]
then we can also assume that $\varepsilon$ is small enough so that
$0\notin\mathcal R_j(\varepsilon)$ for every $j\in J$.

Let $j\in J$ and $\eta\in\R$ be such that $\llvert  \eta\rrvert  \leq\varepsilon
$. Denote by $z_\varepsilon=u_{\min}^{(j)}(\varepsilon)+i\varepsilon
$ and $z_\eta=u_{\min}^{(j)}(\varepsilon)+i\eta$. %
%\textcolor{red}{[me semble que $z_\epsilon$ n'a pas {\chr"C3}{
%\chr"A9}t{\chr"C3}{\chr"A9} encore d{\chr"C3}{\chr"A9}fini]}
Since $\llvert  1-xz_{\eta}\rrvert  \leq\llvert  1-xz_\varepsilon\rrvert  $ for every $x\in\R$, it
follows that
\[
\int\log\llvert 1-x z_\eta\rrvert \nu(\d x) \leq \int\log\llvert 1-x
z_\varepsilon\rrvert \nu(\d x)
\]
and, together with (\ref{upbordervalue}), that
%
%e55 #&#
\begin{equation}
\re f(z_\eta) \geq \re f(z_\varepsilon)+\log\biggl\llvert
\frac{z_\eta
}{z_{\varepsilon}}\biggr\rrvert = \re f(\frak d)+\log\biggl\llvert
\frac
{z_\eta}{z_{\varepsilon}}\biggr\rrvert.
\end{equation}
Next, we have
%
%e56 #&#
\begin{eqnarray}
\label{ineqzeta} \log\biggl\llvert \frac{z_\eta}{z_{\varepsilon}}\biggr\rrvert & =&
\frac
{1}{2}\log \biggl(\frac{u_{\min}^{(j)}(\varepsilon)^2+\eta
^2}{u^{(j)}_{\min}(\varepsilon)^2+\varepsilon^2} \biggr) = \frac
{1}{2}\log
\biggl(1-\frac{\varepsilon^2-\eta^2}{u^{(j)}_{\min
}(\varepsilon)^2+\varepsilon^2} \biggr)
\nonumber\\[-8pt]\\[-8pt]\nonumber
&\geq& \frac{1}{2}\log \biggl(1-\frac{\varepsilon^2}{u^{(j)}_{\min
}(\varepsilon)^2} \biggr) \geq-
\frac{\varepsilon^2}{u^{(j)}_{\min
}(\varepsilon)^2},
\end{eqnarray}
where for the last inequality we use that $\log(1-x)\geq-2x$ for any
$x\in[0,1/2]$. By combining (\ref{assumptionepsilon})--(\ref{ineqzeta}), we have shown that
%
%e57 #&#
\begin{equation}
\label{bordsrect} \re f \bigl(u_{\min}^{(j)}(\varepsilon)+i\eta
\bigr)>\re f(z_0),\qquad\llvert \eta\rrvert \leq\varepsilon,
j\in J.
\end{equation}
The same line of arguments also shows that
%
%e58 #&#
\begin{equation}
\label{bordsrect2} \re f \bigl(u_{\max}^{(j)}(\varepsilon)+i\eta
\bigr)>\re f(z_0),\qquad\llvert \eta\rrvert \leq\varepsilon,
j\in J.
\end{equation}

Now, consider the set
\[
\widetilde\Omega_*= \bigl\{z\in\Omega_*\dvtx  \operatorname{Im}(z)\geq \varepsilon \bigr\}\cup
\bigl\{z\in\Omega^{\mathrm{sym}}_*\dvtx  \operatorname{Im}(z)\leq-\varepsilon \bigr\}\cup
\biggl(\bigcup_{j\in J}\mathcal R_j(
\varepsilon) \biggr),
\]
and\vspace*{1pt} notice it is a bounded open set containing $z_0$
[since $\operatorname{Im}(z_0)\geq\varepsilon_0>\varepsilon$], but which may not be
connected, and which does not contain the origin.\vspace*{2pt} Let $\widetilde
\Omega_*(z_0)$ be the connected component of $\widetilde\Omega_*$
which contains $z_0$. Since $0\notin\widetilde\Omega_*(z_0)$, $\re
f$ is superharmonic on the bounded domain $ \widetilde\Omega_*(z_0)$.
It follows from (\ref{valueborder}), (\ref{bordsrect}), (\ref
{bordsrect2}) and the symmetry $\re f(z)=\re f(\overline z)$ that $\re
f>\re f(z_0)$ on $\partial\widetilde\Omega_*(z_0)$. This yields a
contradiction with the minimum principle for superharmonic functions,
and (2)(a) follows.

We now turn to (2)(b) and again argue by contradiction. Let $\Omega_*$
%\not\subset\R$
be connected component of $\Omega_-$ such that its interior $\operatorname{int}(\Omega_*)$ is not connected. Notice that since $\re f$ is
continuous on $\C\setminus\mathcal K$, we have $\operatorname{int}(\Omega
_*)\setminus\mathcal K= \Omega_*\setminus\mathcal K$, and in
particular (\ref{valueborder}) yields
\[
\re f(z)=\re f(\frak d),\qquad z\in\partial \operatorname{int}(\Omega _*)\setminus
\mathcal K.
\]
If $\Omega_*$ is bounded, then by (2)(a) we have $0\in\Omega_*$, and
moreover, since $\re f(z)\to-\infty$ as $z\to0$, $0\in\operatorname{int}(\Omega_*)$. Let $\Omega_*'$ be a connected component of $\operatorname{int}(\Omega_*)$ which does not contain the origin. It is then a
bounded domain on which $\re f<\re f(\frak d)$ and $\re f=\re f(\frak
d)$ on $\partial\Omega_*'\setminus\mathcal K$. By picking $z_0\in
\Omega_*'\cap\{\operatorname{Im}(z)>0\}$ and by performing the same
construction as in the proof of (2)(a), but replacing $\Omega_*$ by
$\Omega_*'$, we obtain a bounded domain $\widetilde\Omega_*'(z_0)$
containing $z_0$ in its interior, on which $\re f$ is superharmonic,
and such that $\re f>\re f(z_0)$ on its boundary. The minimum principle
for superharmonic functions shows that this is impossible, and (2)(b) follows.

To prove (2)(c), assume now that $0\notin\Omega_*$, so that $\Omega_*$
is necessarily unbounded by~(2)(a). By using that $\operatorname{int}(\Omega
_*)\setminus\mathcal K= \Omega_*\setminus\mathcal K$ where $\mathcal
K$ is a compact set, that $\Omega_-$ has a unique unbounded connected
component by Lemma~\ref{behaviorinfinity}, and that by assumption
$\operatorname{int}(\Omega_*)$ is not connected, it follows that at least one
connected component of $\operatorname{int}(\Omega_*)$, say $\Omega_*'$, is
bounded. Since by assumption $0\notin\Omega_*'$, the same argument as
in the proof of (2)(b) yields a contradiction, and (2)(c) is proved.
\end{pf}

Recall that the $\Delta_k$'s are defined in Lemma~\ref{trianglef},
that $\Delta_{-1}, \Delta_1$ and $\Delta_3$ are in $\Omega_-$ and
$\Delta_{-2},\Delta_0$ and $\Delta_2$ are in $\Omega_+$. Recall
also that the $\Omega_k$'s are the associated connected components
containing the $\Delta_k$'s.
We use the previous lemmas to describe the sets $\Omega_k$'s.

%
%le4.12 #&#
\begin{lemma} \label{studylevelsets}
The following hold true:
\begin{longlist}[(2)]
\item[(1)] The sets $\Omega_1$ and $\Omega_{-1}$ are equal, with a
connected interior and unbounded. In particular,
for every $0<\alpha<\pi/2$ there exists $R>0$ such that
%
%e59 #&#
\begin{equation}
\label{omega1R} \biggl\{z\in\C\dvtx  \llvert z\rrvert >R, -\frac{\pi}{2}+\alpha<
\arg(z)<\frac{\pi
}{2}-\alpha \biggr\}\subset\Omega_1.
\end{equation}
\item[(2)] The sets $\Omega_2$ and $\Omega_{-2}$ are equal, open,
connected and unbounded. In particular there exists $R>0$ such that
%
%e60 #&#
\begin{equation}
\label{omega2R} \biggl\{z\in\C\dvtx  \llvert z\rrvert >R, \frac{\pi}{2}+\alpha<
\arg(z)<\frac{3\pi
}{2}-\alpha \biggr\}\subset\Omega_2.
\end{equation}
%
%\textcolor{blue}{ \item{\mathbf[old]}
% We have $\Omega_{1}=\Omega_{-1}$, the interior of $\Omega_1$ is
%connected, and for every $0<\alpha<\pi/2$ there exists $R>0$ such that
%\begin{equation}
%\label{omega1R}
%\left\{z\in\C\dvtx  \quad\left\vert z\right\vert >R,\quad-\frac{\pi}{2}+\alpha<\arg(z)<\frac{
%\pi}{2}-\alpha\right\}\subset\Omega_1.
%\end{equation}
%\item{\mathbf[old]}
% We have $\Omega_{2}=\Omega_{-2}$, the interior of $\Omega_2$ is
%connected, and there exists $R>0$ such that
%\begin{equation}
%\label{omega2R}
%\left\{z\in\C\dvtx  \quad\left\vert z\right\vert >R,\quad\frac{\pi}{2}+\alpha<\arg(z)<\frac{3
%\pi}{2}-\alpha\right\}\subset\Omega_2.
%\end{equation}
%}
%
\item[(3)]
The interior of $\Omega_3$ is connected, and there exists $\delta>0$
such that $B(0,\delta)\subset\Omega_3$.
\end{longlist}
\end{lemma}

\begin{pf} We first prove (2). Since $\Omega_2$ is by definition a
connected subset of $\Omega_+$, Lemma~\ref{keylevelsets}(1) yields
that it is open. Next, we show by contradiction that $\Omega_2$ is
unbounded. If $\Omega_2$ is bounded, then Lemma~\ref{keylevelsets}(1)
shows there exists $x\in\supp(\nu)$ such that $x^{-1}\in\Omega_2$.
If $x^{-1}<\frak d$ (resp., $x^{-1}>\frak d$), then it follows from the
symmetry $\re f(\overline z)=\re f(z)$ that $\Omega_2$ completely
surrounds $\Omega_3$ (resp., $\Omega_1$); see, for instance,
Figure~\ref{figptselle-right}.
Moreover, Lemma~\ref{trianglef} implies that $\Omega_3$ (resp.,
$\Omega_1$) has nonempty interior.
As a consequence, $\Omega_3$ (resp., $\Omega_1$)
%\not\subset\R$ (resp. $\Omega_1\not\subset\R$)}
is a bounded connected component of $\Omega_-$ which does not contain
the origin, and Lemma~\ref{keylevelsets}(2)(a) shows this is impossible.
The symmetry $\re f(\overline z)=\re f(z)$ moreover provides that
$\Omega_{-2}$ is also unbounded, and~(2) follows from the inclusion
(\ref{omegaR+}) and the fact that $\Omega_+$ has a unique unbounded
connected component; see Lemma~\ref{behaviorinfinity}.

%
%f7 #&#
\begin{figure}

\includegraphics{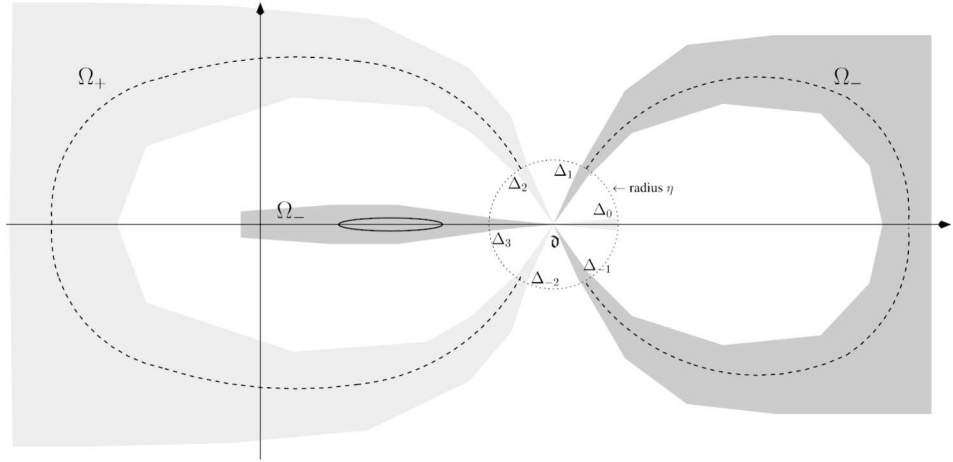}

\caption{Preparation of the saddle point analysis for a right edge.
% Regions $\Delta_k$ defined in Lemma~\ref{trianglef} shown in shaded
%grey.
The dotted path at the right is $\Upsilon^{(1)}_{\mathrm{res}}(N_0)$. The
dotted path at
the left is its counterpart for $\Ttilde$. The\vspace*{1pt} closed contour at the
left of
$\frak d$ is~$\Upsilon^{(0)}$.}\label{figptselle-right}
\end{figure}

We now prove (1). Since $\Omega_2$ is unbounded and symmetric around the
real axis, then $\Omega_1$ does not contain the origin, and it follows
from Lemma~\ref{keylevelsets}(2)(a), (2)(c) that $\Omega_1$ is unbounded
and has a connected interior. Then (1) follows from symmetry $\re
f(\overline z)=\re f(z)$, the inclusion (\ref{omegaR-}) and the fact
that $\Omega_-$ has a unique unbounded connected component; cf. Lemma
\ref{behaviorinfinity}.

Finally, since $\Omega_3$ is bounded as a byproduct of Lemma~\ref
{studylevelsets}(2), it has a connected interior [Lemma~\ref
{keylevelsets}(2)(b)] and contains the origin [Lem-\break ma~\ref
{keylevelsets}(2)(a)]. Moreover, since $\re f(z)\to-\infty$ as $z\to
0$, (3) follows.
\end{pf}

We are finally in position to prove Proposition~\ref{contoursprop}.

\begin{pf*}{Proof of Proposition~\ref{contoursprop}}
Given any $\rho>0$ small enough, it follows from the convergence of
$\frak d_N$
to $\frak d$ that for all $N_0$ large enough, the points
$\frak d_{N_0}+\rho e^{i\pi/3}$ and $\frak d_{N_0}+\rho e^{-i\pi/3}$
belong to
$\Delta_1$ and $\Delta_{-1}$, respectively. Thus both points belong to
$\Omega_1$ by
Lemma~\ref{studylevelsets}(1). As a consequence, we can complete the path
$\{\frak d_{N_0}+te^{\pm i\pi/3}\dvtx   t\in[0,\rho]\}$
into a (closed) contour with a
path $\Upsilon^{(1)}_{\mathrm{res}}(N_0)$ lying in the interior of $\Omega_1$;
see Figure~\ref{figptselle-right}. Since $\Upsilon^{(1)}_{\mathrm{res}}(N_0)$
lies in the interior of $\Omega_1$, the convergence
$\frak d_N\rightarrow\frak d$
moreover yields that we can perform the same construction for all
$N\geq N_0$
with $\Upsilon^{(1)}_{\mathrm{res}}(N)$ in a closed tubular neighborhood
$\mathcal T\subset\Omega_1$ of
$\Upsilon^{(1)}_{\mathrm{res}}(N_0)$. By Lemma~\ref{studylevelsets}(1) again,
we can moreover
choose $\Upsilon^{(1)}_{\mathrm{res}}(N_0)$ in a way that it has finite length
and only
crosses the real axis at a real number lying on the right of $\mathcal
K$. By
construction, this yields that the set $\mathcal T$ is compact and that the
$\Upsilon^{(1)}_{\mathrm{res}}(N)$'s can be chosen with a uniformly bounded
length as long
as $N\geq N_0$. Since $\Omega_1\subset\Omega_-$, there exists $K>0$
such that
$\re f(z)\leq\re f(\frak d)-3K$ on $\mathcal T$. Since moreover $\re
f_N$ uniformly
converges to $\re f$ on $\mathcal T$ and $\re f_N(\frak d_N)\rightarrow
\re f(\frak d
)$, according to Lemma~\ref{lemmaproperties-f-fN}(a), we can choose $N_0$
large enough such that $\re f_N\leq\re f+K$ on $\mathcal T$ and $\re
f(\frak d)\leq\re
f_N(\frak d_N)+K$. This finally yields that $\re(f_N(z)-f_N(\frak
d_N))\leq
-K$ for all $z\in\mathcal T$ and proves the existence of a contour
$\Upsilon^{(1)}$ satisfying
the requirements of Proposition~\ref{contoursprop}, except for point (4).
Similarly, the same conclusion for $\Ttilde$ follows from the same
lines, but
by using $\Omega_2$ instead of $\Omega_1$ and Lemma~\ref{studylevelsets}(2).

As a consequence of Lemma~\ref{studylevelsets}(3), there exists a
contour in
the interior of $\Omega_3$ surrounding $\{x\in\mathcal K\dvtx  0<x<\frak
d\}$ but staying in
$\{z\in\C\dvtx  \re z>0\}$ and which intersects exactly twice the real
axis in
$\R\setminus\mathcal K$ with finite length; see Figure~\ref
{figptselle-right}. Using again
Lemma~\ref{lemmaproperties-f-fN}(a), the existence of $\Upsilon
^{(0)}$ with the
properties provided in the statement of Proposition~\ref{contoursprop}
follows.

Finally, item (4) of Proposition~\ref{contoursprop} is clearly
satisfied by construction since the sets $\Omega_-$ and $\Omega_+$
are disjoint, and the proof of the proposition is therefore complete.
\end{pf*}

%The next step is to show that $ \K_N^{(0)}$ will actually not
%contribute in the large $N$ limit.

%s4.5 #&#
\subsection{Asymptotic analysis for the right edges and proof of Theorem \texorpdfstring{\protect\ref{thfluctuations-TW}\textup{(b)}}{3(b)}}
\label{secasmptotic-analysis-right}

Recall that $\widetilde\K_N=\K_N^{(0)} + \K_N^{(1)}$. We now
analyze the asymptotic behavior of $\K_N^{(0)}$ in the next section
and then investigate $\K_N^{(1)}$ in Section~\ref{KN1}.

%s4.5.1 #&#
\subsubsection{Asymptotic analysis for $\K_N^{(0)}$}\label{K0playsnorole}

Recall definition (\ref{kernelsplitdef}) of the kernel $\K_N^{(0)}$
and its associated contours $\Gamma^{(0)}$ and $\Theta$; cf.
Figure~\ref{fignewcontours}. The aim of this section is to establish
the following statement, which asserts that $ \K_N^{(0)}$ will have no
impact on the asymptotic analysis in the large $N$ limit.

%
%pr4.13 #&#
\begin{proposition}
\label{splitprop}
Let Assumptions~\ref{assgauss}~and~\ref{assnu} hold true. Then for
every $\varepsilon>0$ small enough,
%
%e61 #&#
%e62 #&#
\begin{eqnarray}
\label{splitA} \lim_{N\rightarrow\infty}{\bigl\llVert \mathbf 1_{(s, \varepsilon
N^{2/3}\delta_N)}
\K_N^{(0)}\mathbf 1_{(s, \varepsilon N^{2/3}
\delta
_N)}\bigr\rrVert
}_2&=&0,
\\
\label{splitB} \lim_{N\rightarrow\infty}\Tr \bigl(\mathbf 1_{(s, \varepsilon
N^{2/3}\delta_N)}{
\K}^{(0)}_N\mathbf 1_{(s, \varepsilon
N^{2/3}\delta
_N)} \bigr) & =&0.
\end{eqnarray}
%
%In particular,
%$$
%\mathrm{det} \left(1- \bs1_{(s, \epsilon N^{2/3}\delta_N)}\K_N^{(0)}
%\bs1_{(s, \epsilon N^{2/3} \delta_N)}\right)
%= \mathrm{det} \left(1- \K_N^{(0)}\right)_{L^2(s, \epsilon N^{2/3}
%\delta_N)}\xrightarrow[N\to\infty]{} 1.
%$$
\end{proposition}

\begin{longlist}
\item[\textit{Notation}.] If a contour $\Gamma$ is parametrized by $\gamma\dvtx I\rightarrow\Gamma$ for some interval $I\subset\R$, then for every
map $h\dvtx \Gamma\rightarrow\C$, we set
\[
\int_{\Gamma}h(z)\llvert \d z\rrvert =\int_Ih
\circ\gamma(t) \bigl\llvert \gamma'(t)\bigr\rrvert \,\d t
\]
when it does make sense. In particular, ${\oint}_\Gamma\llvert  \,\d z\rrvert  $ is the
length of the contour~$\Gamma$.
\end{longlist}

\begin{pf*}{Proof of Proposition \ref{splitprop}} %[Proof of Proposition~\ref{splitprop}]
Recall that by definition of $\K_N^{(0)}(x,y)$ [see~(\ref{kernelsplitdef})], we have
%
%e63 #&#
\begin{eqnarray}
\K_{N}^{(0)}(x,y)
&=& \frac{N^{1/3}}{(2i\pi)^2\delta_N}\nonumber
\\
&&{}\times \oint_{\Gamma
^{(0)}}\d z\oint_{\Theta} \d w
\frac{1}{w-z}
\\
&&\hspace*{58pt}{}\times e^{- N^{1/3}x(z-\frak
d_N)/\delta_N +N^{1/3}y(w-\frak d_N)/\delta_N+N f_N(z)-N f_N(w)},
\nonumber
\end{eqnarray}
where $\Theta$ and $\Gamma^{(0)}$ are as in Figure~\ref{fignewcontours}.
%
% is a contour encircling the origin and $\Gamma^{(0)}$ is a contour
%encircling the $1/\lambda_i$'s which are smaller than $\frak c_{N}$.
%Observe now that we can deform the contour $\Theta$ so that it
%encloses the origin and $\Gamma^{(0)}$, and passes through $\frak
%c_N$. Indeed, by using the residue theorem, we see that the residue
%picked at $w=z$ for every $z\in\Gamma^{(0)}$ yields the extra term
%\[
%\frac{N^{1/3}}{2i\pi\sigma_N}\oint_{\Gamma^{(0)}}e^{N^{1/3}(y-x)(z-
%\frak c_N)/\sigma_N}\d z
%\]
%which identically vanishes since the integrand is analytic.
We now deform the contours $\Theta$ and $\Gamma^{(0)}$ so that
$\Theta=\Ttilde$ and $\Gamma^{(0)}=\Upsilon^{(0)}$ where $\Ttilde$
and $\Upsilon^{(0)}$ are given by Proposition~\ref{contoursprop}.
%Thus we can and do specify $\Gamma^{(0)}$ and $\Theta$ to be as in
%Proposition~\ref{contoursprop} with $\rho<\rho_0$, where $\rho_0$
%comes from Lemma~\ref{lemmaproperties-f-fN}-(b).
As a consequence of Proposition~\ref{contoursprop}(4), we have the
upper bound
\begin{eqnarray}
\label{K01}
\qquad\bigl\llvert \K_{N}^{(0)}(x,y)\bigr\rrvert
&\leq&\frac{N^{1/3}}{\d(2\pi)^2
\delta_N}\oint_{\Upsilon^{(0)}}e^{-N^{1/3}x\re(z-\frak d_N)/\delta
_N+N\re(f_N(z)-f_N(\frak d_N))}\llvert \d z\rrvert
\nonumber\\[-8pt]\\[-8pt]\nonumber
&&{} \times\oint_{\Ttilde
}e^{N^{1/3}y\re(w-\frak d_N)/\delta_N-N\re(f_N(w)-f_N(\frak
d_N))}\llvert \d w\rrvert.
\nonumber
\end{eqnarray}
Recall that $\Upsilon^{(0)}$ does not depend on $N$.
By Proposition~\ref{contoursprop}(5)(b), the contour $\Ttilde$ lies in a
compact set.
Hence there exists $L>0$ independent of $N$ such that $\llvert  \re(z-\frak c_N)\rrvert  \leq L$ for $z\in\Upsilon^{(0)}$ or $z\in\Ttilde$. Together
with Proposition~\ref{contoursprop}(3)(a), we obtain that for all $x\geq s$,
%
%e64 #&#
\begin{eqnarray}
\label{K02}
&& \oint_{\Upsilon^{(0)}}e^{-N^{1/3}x\re(z-\frak d_N)/\delta_N+N\re
(f_N(z)-f_N(\frak d_N))}\llvert \d z\rrvert
\nonumber\\[-8pt]\\[-8pt]
&&\qquad \leq e^{-NK+\sklfrac{N^{1/3}}{\delta_N}  ( L (x-s)+L \llvert  s\rrvert
)}\oint_{\Upsilon^{(0)}}\llvert \d z\rrvert.
\nonumber
\end{eqnarray}
Similarly, by splitting $\Ttilde$ into $\Ttilde_{\mathrm{res}}$ and $\Ttilde
_*$, we get from Proposition~\ref{contoursprop}(3)(c) for every $y\geq s$
%
%e65 #&#
\begin{eqnarray}
\label{K03}
&& \oint_{\Ttilde}e^{N^{1/3}y\re(w-\frak d_N)/\delta_N-N\re
(f_N(w)-f_N(\frak d_N))}\llvert \d w\rrvert\nonumber
\\
&&\qquad \leq e^{N^{1/3}L (y-s)/\delta_N+N^{1/3}L \llvert  s\rrvert  /\delta_N}
\\
&&\quad\qquad{}\times \biggl(e^{-NK}\int_{\Ttilde_{\mathrm{res}}}
\llvert \d w\rrvert +\int_{\Ttilde_*} e^{-N\re
(f_N(w)-f_N(\frak d_N))}\llvert \d w
\rrvert \biggr).
\nonumber
\end{eqnarray}
The definition of $\Ttilde_*$ and Lemma~\ref{lemmaproperties-f-fN}(b)
then yield
%
%e66 #&#
\begin{eqnarray}
\label{K04} \int_{\Ttilde_*} e^{-N\re(f_N(w)-f_N(\frak d_N))}\llvert \d w
\rrvert &\leq& \int_{\Ttilde_*} e^{-Ng_N''(\frak d_N)\re(w-\frak d_N)^3+N\Delta\llvert  w-\frak d_N\rrvert  ^4}\llvert \d w\rrvert
\nonumber
\\
&\leq& \int_{\Ttilde_*} e^{-Ng_N''(\frak d_N)\re(w-\frak
d_N)^3+N\rho\Delta\llvert  w-\frak d_N\rrvert  ^3}\llvert \d w\rrvert
\\
& =& 2\int_0^\rho e^{-Nt^3(g_N''(\frak d_N) -\rho\Delta) }\,\d t \leq 2
\rho\nonumber
\end{eqnarray}
provided that $\rho$ is chosen small enough so that $g_N''(\frak d_N)
-\rho\Delta>0$.

By combining (\ref{K01})--(\ref{K04}), we thus obtained that there
exist constants $C_0,C_1>0$ independent of $N$ such that for every
$x,y\geq s$ and every $N$ large enough,
%
%e67 #&#
\begin{equation}
\label{estimateK0} \bigl\llvert \K_{N}^{(0)}(x,y)\bigr\rrvert
\leq C_0 e^{-C_1N+\sklfrac{N^{1/3}}{\delta_N} 2L(x+y)}.
\end{equation}
Since by (\ref{defHS-kernel}),
\[
{\bigl\llVert \mathbf 1_{(s, \varepsilon N^{2/3}\delta_N)}\K _N^{(0)}
\mathbf  1_{(s, \varepsilon N^{2/3} \delta_N)}\bigr\rrVert }_2= \biggl(\int
_{s}^{\varepsilon N^{2/3}\delta_N}\int_{s}^{\varepsilon N^{2/3}\delta
_N}
\K_N^{(0)}(x,y)^2\,\d x\,\d y \biggr)^{1/2},
\]
we obtain from (\ref{estimateK0}) the rough estimate
\[
{\bigl\llVert \mathbf 1_{(s, \varepsilon N^{2/3}\delta_N)}\K _N^{(0)}
\mathbf  1_{(s, \varepsilon N^{2/3} \delta_N)}\bigr\rrVert }_2\leq C_0
\bigl(\varepsilon N^{2/3}\delta_N-s\bigr)e^{-N(C_1-4\varepsilon L)}
\]
from which (\ref{splitA}) follows, provided that we choose
$\varepsilon$ small enough.
Similarly, by~(\ref{defTrace-kernel})
\[
\Tr \bigl(\mathbf 1_{(s, \varepsilon N^{2/3}\delta_N)}{\K }^{(0)}_N
\mathbf  1_{(s, \varepsilon N^{2/3}\delta_N)} \bigr)=\int_{s}^{\varepsilon
N^{2/3}\delta_N}
\K_N^{(0)}(x,x)\,\d x,
\]
and (\ref{estimateK0}) yields the estimate
\begin{eqnarray*}
\bigl\llvert \Tr \bigl(\mathbf 1_{(s, \varepsilon N^{2/3}\delta_N)}{\K }^{(0)}_N
\mathbf  1_{(s, \varepsilon N^{2/3}\delta_N)} \bigr)\bigr\rrvert &\leq& \int_{s}^{\varepsilon N^{2/3}\delta_N}
\bigl\llvert \K_N^{(0)}(x,x)\bigr\rrvert \,\d x
\\
&\leq& C_0\bigl(\varepsilon N^{2/3}\sigma_N-s
\bigr)e^{-N(C_1-4\varepsilon L)},
\end{eqnarray*}
which proves (\ref{splitB}) as soon as $\varepsilon$ is small enough.
Proof of Proposition~\ref{splitprop} is therefore complete.
\end{pf*}

%Next, we investigate the convergence of $\K_N^{(1)}$ toward $\K_{\Aii}$
%and complete the proof of Theorem~\ref{thfluctuations-TW}.

%s4.5.2 #&#
\subsubsection{Asymptotic analysis for $ \K_N^{(1)}$ and proof of Theorem\texorpdfstring{ \protect\ref{thfluctuations-TW}\textup{(b)}}{3(b)}}
\label{KN1}

We now investigate the convergence of $\K_N^{(1)}$ toward $\K_{\Aii}$
and thereafter complete the proof of Theorem~\ref{thfluctuations-TW}(b).
%establish the following statement, where we recall that the operator $
%\K_\Aii$ is the integral operator associated with the Airy kernel $\K_
%\Ai(x,y)$ introduced in \eqref{Kai}.

%
%pr4.14 #&#
\begin{proposition}
\label{keyTW} For every $\varepsilon>0$ small enough, we have
%
%e68 #&#
%e69 #&#
\begin{eqnarray}
\label{HSAiry} \lim_{N\rightarrow\infty}{\bigl\llVert \mathbf 1_{(s,\varepsilon
N^{2/3}\delta_N)}
\bigl(\K_N^{(1)}-\K_\Aii\bigr)
\mathbf 1_{(s,\varepsilon N^{2/3}
\delta_N)}\bigr\rrVert }_2&=&0,
\\
\label{TrAiry} \lim_{N\rightarrow\infty}\Tr \bigl(\mathbf 1_{(s,\varepsilon
N^{2/3}\delta_N)}
\bigl(\K_N^{(1)}-\K_{\Aii}\bigr)
\mathbf 1_{(s,\varepsilon N^{2/3}
\delta_N)} \bigr)&=& 0.
\end{eqnarray}
\end{proposition}

First, we represent the Airy kernel as a double complex integral. To do
so, we introduce for some $\delta>0$, which will be specified later,
the contours
%
%
%e70 #&#
%e71 #&#
\begin{eqnarray}
\label{gammainftydef} \Gamma^\infty&=& \bigl\{{ \frak d_N +}
\delta e^{i\pi\theta}\dvtx  \theta\in[-\pi/3,\pi/3] \bigr\}\cup \bigl\{\frak
d_N+te^{\pm i\pi
/3}\dvtx  t\in[\delta,\infty) \bigr\},
\\
\label{thetainftydef} \Theta^\infty&=& \bigl\{{ \frak d_N +}
\delta e^{i\pi\theta}\dvtx  \theta\in[2\pi/3,4\pi/3] \bigr\}\cup \bigl\{\frak
d_N-te^{\pm
i\pi/3}\dvtx  t\in[\delta,\infty) \bigr\},
\end{eqnarray}
and prove the following.

%
%le4.15 #&#
\begin{lemma} For every $\delta>0$ and $x,y\in\R$, we have
\label{Airydoubleintegral}
\begin{eqnarray*}
\K_\Aii(x,y) &=&\frac{N^{1/3}}{(2i\pi)^2\delta_N}
\\
&&{} \oint_{\Gamma
^\infty}\,\d z
\oint_{\Theta^\infty}\,\d w\frac
{1}{w-z}e^{-N^{1/3}\sklfrac{x(z-\frak d_N)}{\delta_N}+\sklfrac{N}6
g_N''(\frak d_N)(z-\frak d_N)^3}
\\
&&\hspace*{56pt}{}\times e^{N^{1/3}\sklfrac{y(w-\frak d_N)}{\delta_N} -\sklfrac{N}6
g_N''(\frak d_N)(w-\frak d_N)^3}.
\end{eqnarray*}
\end{lemma}

\begin{pf}
First, it easily follows from the differential equation
satisfied by the Airy function, namely ${\Ai}''(x)=x\Ai(x)$, and an
integration by part that
%
%e72 #&#
\begin{equation}
\label{Airyintegral} \K_\Aii(x,y)=\int_0^\infty
\Ai(x+u)\Ai(y+u)\,\d u.
\end{equation}
The Airy function admits the following complex integral representation
(see, e.g., \cite{book-olver-1974}, page~53)
%
%e73 #&#
\begin{equation}
\label{Airysimplecontour} \Ai(x)=-\frac{1}{2i\pi}\oint_{\Xi}e^{-xz+z^3/3}
\,\d z=\frac
{1}{2i\pi}\oint_{\Xi'}e^{xw-w^3/3}\,\d w,
\end{equation}
where $\Xi$ and $\Xi'$ are disjoint unbounded contours, and $\Xi$
goes from $e^{i\pi/3}\infty$ to $e^{-i\pi/3}\infty$ whereas $\Xi'$
goes from $e^{-2i\pi/3}\infty$ to $e^{2i\pi/3}\infty$. By plugging
(\ref{Airysimplecontour}) into (\ref{Airyintegral}) and by using the
Fubini theorem, we obtain
%
%e74 #&#
\begin{eqnarray}\label{Airykerneldoublecontours}
\qquad\K_\Aii(x,y) & =&-\frac{1}{(2i\pi)^2}\oint_{\Xi}\d z
\oint_{\Xi'} \d w\,e^{-xz+z^3/3+yw-w^3/3}\int_0^\infty
e^{u(w-z)}\,\d u
\nonumber\\[-8pt]\\[-8pt]\nonumber
& =& \frac{1}{(2i\pi)^2}\oint_{\Xi}\d z\oint_{\Xi'}\d w
\frac
{1}{w-z}e^{-xz+z^3/3+yw-w^3/3},
\end{eqnarray}
since $\re(w-z)<0$ for all $z\in\Xi$ and $w\in\Xi'$. Lemma~\ref
{Airydoubleintegral} then follows after the changes of variables
$z\mapsto N^{1/3}(z-\frak d_N)/\delta_N$ and $w\mapsto N^{1/3}(w-\frak
d_N)/\delta_N$, the mere definition $\delta_N^3=2/g_N''(\frak d_N)$
and an appropriate deformation of the contours.
\end{pf}

We now turn to the proof of Proposition~\ref{keyTW}.

\begin{pf*}{Proof of Proposition~\ref{keyTW}}
Recall that
%
%e75 #&#
\begin{eqnarray}\label{K1kernelvariation}
\K_{N}^{(1)}(x,y)
&=& \frac{N^{1/3}}{(2i\pi)^2\delta_N}\nonumber
\\
&&{}\times \oint_{\Gamma^{(1)}} \d z \oint_{\Theta} \d w
\frac{1}{w-z}
\\
&&\hspace*{58pt}{}\times e^{- N^{1/3}\sklfrac{x(z-\frak
d_N)}{\delta_N}+N^{1/3}\sklfrac{y(w-\frak d_N)}{\delta_N}+N f_N(z)-N
f_N(w)}.\nonumber
\end{eqnarray}
The key step in the analysis is to deform the contours $\Gamma^{(1)}$
and $\Theta$ into $\Upsilon^{(1)}$ and $\Ttilde$ of Proposition~\ref
{contoursprop}, but since the later intersect
in $\frak d_N$, we need to slightly modify them.

Let $\rho_0$ be fixed so that Lemma~\ref{lemmaproperties-f-fN} holds
true, fix $\rho\le\rho_0$ and recall the definitions of
%
%e76 #&#
\begin{equation}
\label{rappel-contours} \Upsilon^{(1)}=\Upsilon_*\cup\Upsilon^{(1)}_{\mathrm{res}}
\quad\mbox {and}\quad\Ttilde=\Ttilde_* \cup\Ttilde_{\mathrm{res}}
\end{equation}
as provided by Proposition~\ref{contoursprop}. Since $\Upsilon_*\cap
\Ttilde_*=\{\frak d_N\}$, we deform them to make them disjoint. Set
%
%e77 #&#
\begin{equation}
\label{defdelta} \delta=N^{-1/3},
\end{equation}
and from now until the end of the proof, denote (with a slight abuse of
notation)
%
%e78 #&#
%e79 #&#
%e80 #&#
\begin{eqnarray}
\Upsilon_* &=& \biggl\{{ \frak d_N+ }\delta e^{i\theta}\dvtx  \theta
\in \biggl[-\frac
{\pi}3,\frac{\pi}3 \biggr] \biggr\}\cup \bigl\{
\frak d_N+te^{\pm
i\sklfrac{\pi}3}\dvtx  t\in[\delta,\rho] \bigr\}
\label{defupsilon-star-modified}
\\
&:=& \Upsilon_{*,1}\cup\Upsilon_{*,2}, \label{defupsilon-star-split}
\\
\Ttilde_* &=& \biggl\{{ \frak d_N+ } \delta e^{i\theta}\dvtx  \theta
\in \biggl[\frac
{2\pi}3,\frac{4\pi}3 \biggr] \biggr\}\cup \bigl\{
\frak d_N-te^{\pm
i\sklfrac{\pi}3}\dvtx  t\in[\delta,\rho] \bigr\}.
\label{deftheta-modified}
\end{eqnarray}
Notice in particular that this deformation provides now the control
\[
\min \bigl\{\llvert w-z\rrvert\dvtx  z\in\Upsilon_*, w\in\Ttilde_* \bigr\}\geq
\delta. %
\]
Now, let $\Gamma^{(1)}= \Upsilon^{(1)}$ and $\Theta=\Ttilde$. We
can also express the Airy contours $\Gamma^\infty$ and $\Theta
^\infty$ as
\begin{eqnarray*}
\Gamma^\infty&=& \Upsilon_* \cup\Gamma^\infty_{\mathrm{res}}
\qquad\mbox {with } \Gamma^\infty_{\mathrm{res}} = \bigl\{ \frak
d_N+te^{\pm i\pi/3}\dvtx  t\in[\rho,\infty) \bigr\},
\\
\Theta^\infty_{\mathrm{res}} &=& \Ttilde_* \cup\Theta^\infty_{\mathrm{res}}
\qquad \mbox{with } \Theta^\infty_{\mathrm{res}} = \bigl\{\frak
d_N-te^{\pm i\pi
/3}\dvtx  t\in[\rho,\infty) \bigr\}.
\end{eqnarray*}

%{\blue
%we specify the contours $\Gamma^{(1)}$ and $\Theta$ to be as in
%Proposition~\ref{contoursprop} with $\rho<\rho_0$ (see Lemma
%\ref{lemmaproperties-f-fN}-(b) for the definition of $\rho_0$), except
%that we deform the paths $\Gamma_*$ and $\Theta_*$ in the following way
%\begin{eqnarray*}
%\Gamma_* & =&\{\delta e^{i\pi\theta}\dvtx  \theta\in[-\pi/3,\pi/3]\}\cup\{
%\frak c_N+te^{\pm i\pi/3}\dvtx   t\in[\delta,\rho]\}\\
%\Theta_*& =&\{\delta e^{i\pi\theta}\dvtx  \theta\in[2\pi/3,4\pi/3]\}\cup\{
%\frak c_N-te^{\pm i\pi/3}\dvtx   t\in[\delta,\rho]\},
%\end{eqnarray*}
%where we set
%\begin{equation}
%\delta=N^{-1/3}.
%\end{equation}
%Thus we can express $\Gamma^\infty$ and $\Theta^\infty$, respectively
%introduced in \eqref{gammainftydef} and \eqref{thetainftydef}, as $
%\Gamma_*\cup\Gamma_{\mathrm{res}}^\infty$ and $\Theta_*\cup\Theta^\infty_{\mathrm{res}}$
%with
%\[
%\Gamma_{\mathrm{res}}^\infty=\{\frak c_N+te^{\pm i\pi/3}\dvtx   t\in[\rho,\infty)\},
%\qquad\Theta_{\mathrm{res}}^\infty=\{\frak c_N-te^{\pm i\pi/3}\dvtx   t\in[\rho,
%\infty)\}.
%\]
%}

It follows from Proposition~\ref{contoursprop}(4) and the definition
of the contours that there exists $\d'$ such that for any
\begin{eqnarray*}
\bigl(\Xi,\Xi'\bigr) &\in& \bigl\{(\Upsilon_*,\widetilde
\Theta_{\mathrm{res}}),\bigl(\Upsilon ^{(1)}_{\mathrm{res}},\widetilde
\Theta_*\bigr),\bigl(\Upsilon^{(1)}_{\mathrm{res}},\widetilde
\Theta_{\mathrm{res}}\bigr),
\\
&&\hspace*{15pt} \bigl(\Gamma_*,\Theta_{\mathrm{res}}^\infty
\bigr), \bigl(\Gamma_{\mathrm{res}}^\infty,\widetilde\Theta_*\bigr),
\bigl(\Gamma_{\mathrm{res}}^\infty,\Theta_{\mathrm{res}}^\infty
\bigr) \bigr\},
\end{eqnarray*}
we have
\[
\min \bigl\{\llvert w-z\rrvert\dvtx  z\in\Xi, w\in\Xi' \bigr\}\geq
\d'. %
\]
%
%. Moreover, it follows form the definitions of $\Gamma_*$ and $
%\Theta_*$ that $\min\{\left\vert w-z\right\vert\dvtx  z\in\Gamma_*, w\in\Theta_*\}\geq
%\delta$.
As a consequence, by using (\ref{K1kernelvariation}), (\ref
{defdelta}), Lemma~\ref{Airydoubleintegral} and by splitting contours
into their different components, we obtain that
\begin{eqnarray}
\label{splitEi}
&& \bigl\llvert \K_N^{(1)}(x,y)-
\K_\Aii(x,y)\bigr\rrvert
\nonumber\\[-8pt]\\[-8pt]\nonumber
&&\qquad \leq\frac{N^{2/3}}{(2\pi
)^2\delta_N}E_0+\frac{N^{1/3}}{\d'(2\pi)^2\delta_N} (E_1+E_2+E_3+E_4+E_5+E_6),
\nonumber
\end{eqnarray}
where, setting for convenience
\begin{eqnarray*}
F_N(x,z) & =&e^{- N^{1/3}\sklfrac{x(z-\frak d_N)}{\delta_N} +N
(f_N(z)-f_N(\frak d_N))},
\\
F_\Aii(x,z) & =& e^{- N^{1/3}\sklfrac{x(z-\frak d_N)}{\delta_N} +\sklfrac{N}6 g_N''(\frak d_N)(z-\frak d_N)^3},
\\
G_N(y,w) & =&e^{N^{1/3}\sklfrac{y(w-\frak d_N)}{\delta_N}
-N(f_N(w)-f_N(\frak d_N))},
\\
G_\Aii(y,w) & =& e^{N^{1/3}\sklfrac{y(w-\frak d_N)}{\delta_N} -\sklfrac{N}6
g_N''(\frak d_N)(w-\frak d_N)^3},
\end{eqnarray*}
we introduce
%
%e81 #&#
%e82 #&#
%e83 #&#
%e84 #&#
%e85 #&#
%e86 #&#
%e87 #&#
\begin{eqnarray}
\label{E0def} E_0 & =&\int_{\Upsilon_*}\llvert \d z
\rrvert \int_{\Ttilde_*}\llvert \d w\rrvert \bigl\llvert
F_N(x,z)G_N(y,w)-F_\Aii(x,z)G_\Aii(y,w)
\bigr\rrvert,
\\
\label{Eidef1} E_1 & =& \biggl(\int_{\Upsilon_*}\bigl
\llvert F_N(x,z)\bigr\rrvert \llvert \d z\rrvert \biggr) \biggl(\int
_{\Ttilde_{\mathrm{res}}}\bigl\llvert G_N(y,w)\bigr\rrvert \llvert \d
w\rrvert \biggr),
\\
E_2 & =& \biggl(\int_{\Upsilon_{\mathrm{res}}}\bigl\llvert
F_N(x,z)\bigr\rrvert \llvert \d z\rrvert \biggr) \biggl(\int
_{\Ttilde_{*}}\bigl\llvert G_N(y,w)\bigr\rrvert \llvert \d
w\rrvert \biggr),
\\
E_3 & =& \biggl(\int_{\Upsilon_{\mathrm{res}}}\bigl\llvert
F_N(x,z)\bigr\rrvert \llvert \d z\rrvert \biggr) \biggl(\int
_{\Ttilde_{\mathrm{res}}}\bigl\llvert G_N(y,w)\bigr\rrvert \llvert \d
w\rrvert \biggr),
\\
E_4 & =& \biggl(\int_{\Upsilon_{*}}\bigl\llvert
F_\Aii(x,z) \bigr\rrvert \llvert \d z\rrvert \biggr) \biggl(\int
_{\Theta_{\mathrm{res}}^\infty}\bigl\llvert G_\Aii(y,w) \bigr\rrvert \llvert \d w\rrvert \biggr),
\\
E_5 & =& \biggl(\int_{\Gamma_{\mathrm{res}}^\infty}\bigl\llvert
F_\Aii(x,z) \bigr\rrvert \llvert \d z\rrvert \biggr) \biggl(\int
_{\Ttilde_{*}}\bigl\llvert G_\Aii(y,w) \bigr\rrvert \llvert \d w\rrvert \biggr),
\\
\label{Eidef2} E_6 & =& \biggl(\int_{\Gamma_{\mathrm{res}}^\infty}\bigl
\llvert F_\Aii(x,z) \bigr\rrvert \llvert \d z\rrvert \biggr) \biggl(
\int_{\Theta_{\mathrm{res}}^\infty}\bigl\llvert G_\Aii(y,w)\bigr\rrvert
\llvert \d w\rrvert \biggr).
\end{eqnarray}

\begin{longlist}
\item[\textit{Convention}.] In the rest of the proof, $C,C_0,C_1,\ldots$
stand for positive constants which are independent on $N$ or $x,y$, but
which may change from one line to an other.
\end{longlist}

\begin{longlist}
\item[\textit{Step} 1: \textit{Estimates for $E_0$}.] We rely on the following
elementary inequality:
%
%e88 #&#
\begin{eqnarray}
\label{ineqdiffkernels} \bigl\llvert e^u-e^v\bigr\rrvert &=&
e^{\re(v)}\bigl\llvert e^{(u-v)}-1\bigr\rrvert
\nonumber\\[-8pt]\\[-8pt]\nonumber
&\leq& e^{\re(v)}\sum_{k\geq1}
\frac{\llvert  u-v\rrvert  ^k}{k!} \leq \llvert u-v\rrvert e^{\re(v)+\llvert  u-v\rrvert  },
\end{eqnarray}
which holds for every $u,v\in\C$.
By combining this inequality for
\begin{eqnarray*}
u& =&N\bigl(f_N(z)-f_N(\frak d_N)\bigr)-N
\bigl(f_N(w)-f_N(\frak d_N)\bigr),
\\
v & =&\frac{Ng_N''(\frak d_N)}6 \bigl\{ (z-\frak d_N)^3-(w-
\frak d_N)^3 \bigr\}
\end{eqnarray*}
together with Lemma~\ref{lemmaproperties-f-fN}(b), we obtain
\begin{eqnarray*}
\hspace*{-3pt}&& \bigl\llvert F_N(x,z)G_N(y,w)-F_\Aii(x,z)G_\Aii(y,w)
\bigr\rrvert
\\
\hspace*{-3pt}&&\qquad\leq\Delta N \bigl(\llvert z-\frak d_N\rrvert
^4+\llvert w-\frak d_N\rrvert ^4\bigr)
e^{-N^{1/3}\sklfrac{x\re(z-\frak d_N)}{\delta_N}+N^{1/3}\sklfrac{y\re
(w-\frak d_N)}{\delta_N}}
\\
\hspace*{-3pt}&&\quad\qquad{} \times e^{\sklfrac{Ng_N''(\frak d_N)}6 \re(z-\frak d_N)^3+N\Delta
\llvert  z-\frak d_N\rrvert  ^4-\sklfrac{Ng_N''(\frak d_N)} 6 \re(w-\frak d_N)^3+N\Delta
\llvert  w-\frak d_N\rrvert  ^4},
\end{eqnarray*}
provided that $z,w\in B(\frak d_N,\rho)$. This yields with (\ref{E0def})
%
%
%e89 #&#
\begin{eqnarray}
\label{ineqE0} E_0 &\leq& \Delta\int_{\Upsilon_*}N
\llvert z-\frak d_N\rrvert ^4e^{-N^{1/3}\sklfrac{x\re(z-\frak d_N)}{\delta_N}}\nonumber
\\
&&\hspace*{24pt}{}\times e^{\sklfrac {Ng_N''(\frak d_N)}6 \re(z-\frak d_N)^3+N\Delta\llvert  z-\frak d_N\rrvert  ^4}
\llvert \d z\rrvert\nonumber
\\
&&{} \times\int_{\Ttilde_*}e^{N^{1/3}\sklfrac{x\re(w-\frak
d_N)}{\delta_N}}e^{-\sklfrac{Ng_N''(\frak d_N)}6 \re(w-\frak
d_N)^3+N\Delta\llvert  w-\frak d_N\rrvert  ^4}\llvert \d w\rrvert
\nonumber\\[-8pt]\\[-8pt]
&&{} + \Delta\int_{\Upsilon_*}e^{-N^{1/3}\sklfrac{x\re(z-\frak
d_N)}{\delta_N} }e^{\sklfrac {Ng_N''(\frak d_N)}6 \re(z-\frak
d_N)^3+N\Delta\llvert  z-\frak d_N\rrvert  ^4}
\llvert \d z\rrvert
\nonumber
\\
&&{} \times\int_{\Ttilde_*}N\llvert w-\frak d_N\rrvert
^4e^{N^{1/3}\sklfrac{x\re
(w-\frak d_N)}{\delta_N}}\nonumber
\\
&&\hspace*{28pt}{}\times e^{-\sklfrac{Ng_N''(\frak d_N)}6 \re(w-\frak
d_N)^3+N\Delta\llvert  w-\frak d_N\rrvert  ^4}\llvert \d w\rrvert.\nonumber
\end{eqnarray}

We first handle the integrals over the contour $\Upsilon_*=\Upsilon
_{*,1} \cup\Upsilon_{*,2}$ [see (\ref{defupsilon-star-split})] and
consider separately the two different portions of the contour. First,
let $z\in\Upsilon_{*,1}$, and recall that $x\geq s$ by
assumption. Since
\[
\frac{\delta}2 \leq \re(z-\frak d_N) \leq \llvert z-\frak
d_N\rrvert \leq \delta \quad\mbox{and} \quad\delta=
\frac{1}{N^{1/3}}, %
\]
we have $\llvert  z-\frak d_N\rrvert  ^4 =N^{-4/3}$ and the estimates
\begin{eqnarray*}
%& \left\vert z-\frak d_N\right\vert ^4 &=& \quad\frac1{N^{4/3}},\\
e^{-N^{1/3}\sklfrac{x\re(z-\frak d_N)}{\delta_N}} &\leq& e^{ -\sklvfrac {x-s}{2\delta_N}+\sklfrac{\llvert  s\rrvert  }{\delta_N}},
\\
e^{Ng_N''(\frak d_N)\re(z-\frak d_N)^3+N\Delta\llvert  z-\frak d_N\rrvert  ^4}&\leq& e^{ g_N''(\frak d_N)+\sklfrac{\Delta}{N^{1/3}}}.
\end{eqnarray*}
This immediately yields
%
%e90 #&#
\begin{eqnarray}
\label{ineq-portion1} && \int_{\Upsilon_{*,1}} N\llvert z-\frak
d_N\rrvert ^4e^{-N^{1/3}\sklfrac{x\re
(z-\frak d_N)}{\delta_N}}e^{\sklfrac {Ng_N''(\frak d_N)}6 \re(z-\frak
d_N)^3+N\Delta\llvert  z-\frak d_N\rrvert  ^4}\llvert \d
z\rrvert\nonumber
\\
&&\qquad\leq \frac{1}{N^{1/3}}e^{-\sklvafrac{x-s}{2\delta_N}+ \sklfrac{\llvert  s\rrvert  }{\delta_N}} e^{g_N''(\frak d_N)+\sklfrac{\Delta}{N^{1/3}}}
\biggl(\frac{2\pi
}{3N^{1/3}} \biggr)
\\
&&\qquad  \leq\frac{C}{N^{2/3}} e^{-\vafrac{x-s}{2\delta_N}},%{N^{2/3}}
\nonumber
\end{eqnarray}
where $2\pi/3N^{1/3}$ accounts for the length of ${\Upsilon_{*,1}}$. Similarly
%
%e91 #&#
\begin{eqnarray}
\label{ineq-bis-portion1}
&& \int_{\Upsilon_{*,1}} e^{-N^{1/3}\sklfrac{x\re(z-\frak d_N)}{\delta
_N}}e^{\sklfrac {Ng_N''(\frak d_N)}6 \re(z-\frak d_N)^3+N\Delta\llvert  z-\frak d_N\rrvert  ^4}
\llvert \d z\rrvert
\nonumber\\[-8pt]\\[-8pt]\nonumber
&&\qquad  \leq \frac{C}{N^{1/3}} e^{-\vafrac{x-s}{2\delta_N}}. %\frac{C'e^{-\vafrac{x-s}{2\delta_N}}}{N^{1/3}}
\end{eqnarray}

Consider now the situation where $z\in\Upsilon_{*,2}$. In this case,
\[
\re(z-\frak d_N)=\frac{t} 2,\qquad\re(z-\frak
d_N)^3=-t^3,\qquad \llvert z-\frak
d_N\rrvert ^4= t^4, %
\]
with $N^{-1/3}\leq t\leq\rho$ and thus
\begin{eqnarray*}
e^{-N^{1/3}\sklfrac{x\re(z-\frak d_N)}{\delta_N}}&\leq& e^{-tN^{1/3} \sklvafrac{x-s}{2\delta_N}+N^{1/3}\sklvafrac{\llvert  s\rrvert  t}{2\delta_N}}
\\
&\le &
e^{-\sklvafrac{x-s}{2\delta_N}+N^{1/3}\sklvafrac{\llvert  s\rrvert  t}{2\delta_N}},
\\
e^{Ng_N''(\frak d_N)\re(z-\frak d_N)^3+N\Delta\llvert  z-\frak d_N\rrvert  ^4}&\leq& e^{-N(g_N''(\frak d_N)-\rho\Delta)t^3}.
\end{eqnarray*}
Assuming that we choose $\rho$ small enough so that $g''(\frak d) -
\rho\Delta>0$ and recalling that $g''_N(\frak d_N)\to g''(\frak d)$,
this provides for every $N$ large enough the inequalities
%
%e92 #&#
\begin{eqnarray}
\label{ineq-portion2} \qquad&& \int_{\Upsilon_{*,2}} N\llvert z-\frak
d_N\rrvert ^4e^{-N^{1/3}\sklfrac{x\re
(z-\frak d_N)}{\delta_N}}e^{\sklfrac {Ng_N''(\frak d_N)}6 \re(z-\frak
d_N)^3+N\Delta\llvert  z-\frak d_N\rrvert  ^4}\llvert \d
z\rrvert\nonumber
\\
&&\qquad\leq2 e^{-\vafrac{x-s}{2\delta_N}} \int_{N^{-1/3}}^\rho
Nt^4e^{N^{1/3}\sklvafrac{\llvert  s\rrvert  t}{2\delta
_N}-N(g_N''(\frak d_N)-\rho\Delta)t^3}\,\d t
\nonumber\\[-8pt]\\[-8pt]\nonumber
&&\qquad\leq\frac{2}{N^{2/3}}e^{-\vafrac{x-s}{2\delta_N}} \int_1^\infty
u^4e^{\sklvafrac{\llvert  s\rrvert  u}{2\delta_N}-(g_N''(\frak d_N)-\rho
\Delta)u^3}\,\d u
\\
&&\qquad \le \frac{C}{N^{2/3}} e^{-\vafrac{x-s}{2\delta_N}}.\nonumber
% \leq& \qquad\frac{Ce^{-(x-s)/(2\delta_N)}}{N^{2/3}}
\end{eqnarray}
%
%\begin{eqnarray}
%\label{ineq-portion2}
% \int_{\Upsilon_{*,2}} N\left\vert z-\frak d_N\right\vert ^4e^{-N^{1/3}x\re(z-\frak d_N)/
%\delta_N}e^{Ng_N''(\frak d_N)\re(z-\frak d_N)^3/6+N\Delta\left\vert z-\frak
%d_N\right\vert ^4}\left\vert \d z\right\vert  \\
%\leq\quad\quad2 e^{-\vafrac{x-s}{2\delta_N}}
%\int_{N^{-1/3}}^\rho Nt^4e^{N^{1/3}\sklfrac{\left\vert s\right\vert t}{2\delta_N}-N(g_N''(
%\frak d_N)-\rho\Delta)t^3}\d t \\
%\leq\quad2 \frac{e^{-\vafrac{x-s}{2\delta_N}}}{N^{2/3}}
%\int_1^\infty u^4e^{\frac{\left\vert s\right\vert u}{2\delta_N}-(g_N''(\frak d_N)-\rho
%\Delta)u^3}\d u \quad\le\quad
%C \frac{e^{-\vafrac{x-s}{2\delta_N}}}{N^{2/3}}
%% \leq& \qquad\frac{Ce^{-(x-s)/(2\delta_N)}}{N^{2/3}}
%\nonumber
%\end{eqnarray}
Similarly,
%
%e93 #&#
\begin{eqnarray}
\label{ineq-bis-portion2}
&& \int_{\Upsilon_{*,2}} e^{-N^{1/3}\sklfrac{x\re(z-\frak d_N)}{\delta
_N}}e^{N \sklfrac{g_N''(\frak d_N)}6 \re(z-\frak d_N)^3+N\Delta\llvert  z-\frak d_N\rrvert  ^4}
\llvert \d z\rrvert
\nonumber\\[-8pt]\\[-8pt]\nonumber
&&\qquad \leq %\quad\quad2 e^{-\vafrac{x-s}{2\delta_N}}
%\int_{N^{-1/3}}^\rho e^{N^{1/3}\sklfrac{\left\vert s\right\vert t}{2\delta_N}-N(g_N''(\frak
%d_N)-\rho\Delta)t^3}\d t \\
%\leq\quad2 \frac{e^{-\vafrac{x-s}{2\delta_N}}}{N^{2/3}}
%\int_1^\infty u^4e^{\left\vert s\right\vert u/(2\delta_N)-(g_N''(\frak d_N)-\rho\Delta)u^3}
%\d u \quad\le\quad
\frac{C}{N^{1/3}} e^{-\vafrac{x-s}{2\delta_N}}.%\frac{e^{-
%\frac{x-s}{2\delta_N}}}{N^{1/3}}
\end{eqnarray}
Gathering (\ref{ineq-portion1})--(\ref{ineq-bis-portion2}), we
finally obtain estimates over the whole contour $\Upsilon_*$,
%
%e94 #&#
\begin{eqnarray}
\label{estimate-E0-1}
&& \int_{\Upsilon_*} N\llvert z-\frak d_N
\rrvert ^4e^{-N^{1/3}\sklfrac{x\re(z-\frak
d_N)}{\delta_N}}e^{\sklfrac {Ng_N''(\frak d_N)}6 \re(z-\frak
d_N)^3+N\Delta\llvert  z-\frak d_N\rrvert  ^4}\llvert \d z\rrvert\hspace*{-20pt}\nonumber
\\
&&\qquad \le
\frac{ C}{N^{2/3}} e^{-\vafrac{x-s}{2\delta_N}},
\nonumber\\[-8pt]\\[-8pt]\nonumber
&& \int_{\Upsilon_*} e^{-N^{1/3}\sklfrac{x\re(z-\frak d_N)}{\delta
_N}}e^{\sklfrac {Ng_N''(\frak d_N)}6 \re(z-\frak d_N)^3+N\Delta\llvert  z-\frak d_N\rrvert  ^4}\llvert \d z
\rrvert
\\
&&\qquad  \le \frac{C}{N^{1/3}} e^{-\vafrac{x-s}{2\delta_N}}.\nonumber
\end{eqnarray}
The same line of arguments also yields equivalent estimates for the
integrals over~$\Ttilde^*$. Namely,
%
%e95 #&#
\begin{eqnarray}
\label{estimate-E0-4}
&& \int_{\Ttilde_*}N\llvert w-\frak d_N
\rrvert ^4e^{N^{1/3}\sklfrac{y\re(w-\frak
d_N)}{\delta_N}}e^{-\sklfrac{Ng_N''(\frak d_N)}6 \re(w-\frak
d_N)^3+N\Delta\llvert  w-\frak d_N\rrvert  ^4}\llvert \d w\rrvert\nonumber\hspace*{-18pt}
\\
&&\qquad  \leq
\frac{C}{N^{2/3}} e^{-\vafrac{y-s}{2\delta_N}},
\nonumber\\[-8pt]\\[-8pt]\nonumber
&& \int_{\Ttilde_*}e^{N^{1/3}\sklfrac{y\re(w-\frak d_N)}{\delta
_N}}e^{-\sklfrac{Ng_N''(\frak d_N)}6\re(w-\frak d_N)^3+N\Delta\llvert  w-\frak d_N\rrvert  ^4}\llvert \d w
\rrvert
\\
&&\qquad  \leq \frac{C}{N^{1/3}} e^{-\vafrac{y-s}{2\delta_N}}.\nonumber
\end{eqnarray}
Combining (\ref{estimate-E0-1})--(\ref{estimate-E0-4}), we have shown that
%
%e96 #&#
\begin{equation}
\label{conclusionE0} E_0\leq\frac{C}{N}e^{-\vafrac{x+y-2s}{2\delta_N}}.
\end{equation}
\end{longlist}

\begin{longlist}
\item[\textit{Step} 2: \textit{Estimates for the remaining $E_i$'s}.]
Using the same estimates as in step~1, we can prove that
%
%e97 #&#
%e98 #&#
\begin{eqnarray}
\label{ineqreste1}
\int_{\Upsilon_{*}}\bigl\llvert F_N(x,z)
\bigr\rrvert \llvert \d z\rrvert &\leq& \frac
{C}{N^{1/3}} e^{-\vafrac{x-s}{2\delta_N}},
\nonumber\\[-8pt]\\[-8pt]\nonumber
\int_{\Upsilon
_{*}}\bigl\llvert F_\Aii(x,z) \bigr
\rrvert \llvert \d z\rrvert &\leq&\frac{C}{N^{1/3}} e^{-\vafrac{x-s}{2\delta_N}},
\\
\int_{\Ttilde_*}\bigl\llvert G_N(y,w)\bigr\rrvert
\llvert \d w\rrvert &\leq& \frac{C}{N^{1/3}} e^{-\vafrac{y-s}{2\delta_N}},
\nonumber\\[-8pt]\\[-8pt]\nonumber
\int _{\Ttilde_*}\bigl\llvert G_\Aii (y,w)\bigr\rrvert \llvert \d w\rrvert &\leq& \frac{C}{N^{1/3}} e^{-\vafrac{y-s}{2\delta_N}}.
\end{eqnarray}

The definitions of the paths and Proposition~\ref{contoursprop} yield
that there exists $L>0$ independent of $N$ such that
\[
\bigl\llvert \re(z-\frak c_N)\bigr\rrvert \leq L,\qquad z\in
\Upsilon_* \cup\Ttilde_* \cup\Upsilon_{\mathrm{res}}\cup\Ttilde_{\mathrm{res}}.
\]
This estimate, together with Proposition~\ref{contoursprop}(3)(b), (3)(c)
and (5)(c) yields that for every $x,y\geq s$,
%
%e99 #&#
%e100 #&#
\begin{eqnarray}
\int_{\Upsilon_{\mathrm{res}}}\bigl\llvert F_N(x,z)\bigr\rrvert
\llvert \d z\rrvert &\leq& C e^{-NK+N^{1/3}L \sklvfrac{x-s}{\delta_N}+N^{1/3}\sklfrac{L \llvert  s\rrvert  }{\delta
_N}},
\\
\int_{\Ttilde_{\mathrm{res}}}\bigl\llvert G_N(y,w)\bigr\rrvert
\llvert \d w\rrvert &\leq& C e^{-NK+N^{1/3}L\sklvafrac{y-s}{\delta_N}+N^{1/3}\sklfrac{L \llvert  s\rrvert  }{\delta_N}}.\label{ineqreste2}
\end{eqnarray}
Combining (\ref{ineqreste1})--(\ref{ineqreste2}), we readily obtain
\[
E_1 + E_2 + E_3 \le C e^{ - C_1 N + C_2 N^{1/3} \sklvfrac{x+y}{\delta
_N}}.
\]
We now handle
\[
\int_{\Gamma_{\mathrm{res}}^\infty}\bigl\llvert F_\Aii(x,z) \bigr\rrvert
\llvert \d z\rrvert \quad\mbox {and}\quad \int_{\Theta_{\mathrm{res}}^\infty}\bigl
\llvert G_\Aii(y,w) \bigr\rrvert \llvert \d w\rrvert. %
\]
We have
\begin{eqnarray*}
\int_{\Gamma_{\mathrm{res}}^\infty}\bigl\llvert F_\Aii(x,z) \bigr\rrvert
\llvert \d z\rrvert & =& \int_{\Theta_{\mathrm{res}}^\infty}\bigl\llvert
G_\Aii(y,w) \bigr\rrvert \llvert \d w\rrvert
\\
& =& 2\int_\rho^\infty e^{-\sklafrac{N^{1/3} x t}{2\delta_N} - \sklfrac
{Ng''_N(\frak d_N)}{6} t^3}\,\d t
\\
&\le& 2
\int_\rho^\infty e^{\sklafrac{N^{1/3} \llvert  s\rrvert   t}{2\delta_N} -
\sklfrac{Ng''_N(\frak d_N)}{6} t^3}\,\d t.
\end{eqnarray*}
Let now $N$ large enough so that
\[
3 \frac{g''_N(\frak d_N) N}6 \rho^2 - \frac{N^{1/3}\llvert  s\rrvert  }{2\delta_N} \geq\rho
\]
(beware that such a condition only depends on $s$). Then
\begin{eqnarray*}
&& 2\int_\rho^\infty e^{\sklafrac{N^{1/3} \llvert  s\rrvert   t}{2\delta_N} - \sklfrac
{Ng''_N(\frak d_N)}{6} t^3}\,\d t
\\
&&\qquad \le
\frac{2}\rho\int_\rho^\infty \biggl(3
\frac{g''_N(\frak d_N)
N}6 t^2- \frac{N^{1/3}\llvert  s\rrvert  }{2\delta_N} \biggr) e^{\sklafrac{N^{1/3} \llvert  s\rrvert   t}{2\delta_N} - \sklfrac{Ng''_N(\frak d_N)}{6}
t^3}\,\d t
\\
&&\qquad \le \frac{2}\rho \bigl[ - e^{\sklafrac{N^{1/3} \llvert  s\rrvert   t}{2\delta_N} -
\sklfrac{Ng''_N(\frak d_N)}{6} t^3} \bigr]_\rho^\infty
\\
&&\qquad= \frac{2}\rho e^{\sklafrac{N^{1/3} \llvert  s\rrvert   \rho}{2\delta_N} - \sklfrac
{Ng''_N(\frak d_N)}{6} \rho^3},
\end{eqnarray*}
and we hence obtain the estimate
%
%e101 #&#
\begin{equation}
\int_{\Gamma_{\mathrm{res}}^\infty}\bigl\llvert F_\Aii(x,z) \bigr\rrvert
\llvert \d z\rrvert = \int_{\Theta_{\mathrm{res}}^\infty}\bigl\llvert
G_\Aii(y,w) \bigr\rrvert \llvert \d w\rrvert \le C e^{- C_1 N}.
\end{equation}
We can now easily handle $E_4$, $E_5$ and $E_6$ and finally obtain
%
%e102 #&#
\begin{equation}
\label{conclusionEi} \sum_{k=1}^6E_k
\leq Ce^{-C_1N+C_2N^{1/3}\sklvfrac{x+y}{\delta_N}}.
\end{equation}
\end{longlist}

\begin{longlist}
\item[\textit{Step} 3: \textit{Conclusions}.]
By combining (\ref{splitEi}), (\ref{conclusionE0}) and (\ref
{conclusionEi}), we have shown for every $x,y\geq s$ and $N$ large
enough that
\[
\bigl\llvert \K_N^{(1)}(x,y)-\K_\Aii(x,y)\bigr
\rrvert \leq\frac
{C}{N^{1/3}}e^{-\vafrac{x+y-2s}{2\delta_N}}+C_1e^{-C_2N+C_3N^{1/3}
\sklvfrac{x+y}{\delta_N}}.
\]
As a consequence,
\begin{eqnarray*}
&& \bigl\llvert \Tr \bigl(\mathbf 1_{(s, \varepsilon N^{2/3}\delta_N)}\bigl(\K _N^{(1)}-
\K_\Aii\bigr)\mathbf 1_{(s, \varepsilon N^{2/3}\delta
_N)} \bigr)\bigr\rrvert
\\
&&\qquad \leq
\int_{s}^{\varepsilon N^{2/3}\delta_N}\bigl\llvert \K _N^{(1)}(x,x)-
\K_\Aii(x,x)\bigr\rrvert \,\d x
\\
&&\qquad \leq \frac{\delta_N C}{N^{1/3}}+\bigl(\varepsilon N^{2/3}\delta
_N-s\bigr)C_1e^{-N(C_2-2\varepsilon C_3)},
\end{eqnarray*}
and (\ref{TrAiry}) follows provided $\varepsilon$ is chosen small
enough. Similarly,
\begin{eqnarray*}
&& {\bigl\llVert \mathbf 1_{(s, \varepsilon N^{2/3}\delta_N)}\bigl(\K _N^{(1)}-
\K _\Aii\bigr)\mathbf 1_{(s, \varepsilon N^{2/3}\delta_N)}\bigr\rrVert
}_2^2
\\
&&\qquad= \int_{s}^{\varepsilon N^{2/3}\delta_N}\int_{s}^{\varepsilon N^{2/3}\delta_N}
\bigl(\K_N^{(1)}(x,y)-\K_\Aii (x,y)
\bigr)^2\,\d x\,\d y
\\
&&\qquad\leq \biggl(\frac{\delta_N C}{N^{1/3}} \biggr)^2+\bigl(\varepsilon
N^{2/3}\delta_N-s\bigr)^2C'_1e^{-N(C_2-2\varepsilon C_3)},
\end{eqnarray*}
where $C'_1>0$ is independent on $N$. This yields (\ref{HSAiry}) as
soon as $\varepsilon$ is chosen small enough and thus completes the
proof of Proposition~\ref{keyTW}.\quad\qed
\end{longlist}\noqed
\end{pf*}

We are finally in position to prove Theorem~\ref{thfluctuations-TW}(b).

\begin{pf*}{Proof of Theorem~\ref{thfluctuations-TW}\textup{(b)}}
First, we check that the Airy operator $\K_\Aii$ is trace class and
Hilbert--Schmidt on $L^2(s,\infty)$ for every $s\in\R$. Indeed,
representation~(\ref{Airyintegral}) provides the factorization $\K
_\Aii=\A_s^2$ of operators on $L^2(s,\infty)$, where $\A_s$ is the
integral operator having for kernel $\A_s(x,y)=\Ai(x+y-s)$. The fast
decay as $x\rightarrow+\infty$ of the Airy function (see \cite{book-olver-1974}, page~394)
%
%e103 #&#
\begin{equation}
\label{Airyasymptotic} \Ai(x)\leq\frac{e^{-\sklfrac{2}{3} x^{3/2}}}{2\pi
^{1/2}x^{1/4}},\qquad x>0,
\end{equation}
then shows that both $\A_s$ and $\K_\Aii$ are Hilbert--Schmidt, and
moreover that $\K_\Aii$ is trace class being the product of two
Hilbert--Schmidt operators.

Next, by using again upper bound (\ref{Airyasymptotic}), it follows
that for every $\varepsilon>0$,
\begin{eqnarray*}
\lim_{N\rightarrow\infty}{\llVert \mathbf 1_{(s,\varepsilon N^{2/3}
\delta_N)}
\K_\Aii\mathbf 1_{(s,\varepsilon N^{2/3} \delta
_N)}-\mathbf  1_{(s,\infty)}
\K_\Aii\mathbf 1_{(s,\infty)}\rrVert }_2&=&0,
\\
\lim_{N\rightarrow\infty}\Tr (\mathbf 1_{(s,\varepsilon N^{2/3}
\delta_N)}\K_\Aii
\mathbf 1_{(s,\varepsilon N^{2/3}\delta_N)} )&=&\Tr ( \mathbf 1_{(s,\infty)}
\K_\Aii\mathbf 1_{(s,\infty
)} ).
\end{eqnarray*}
Together with Proposition~\ref{keyTW}, this yields
\begin{eqnarray*}
\lim_{N\rightarrow\infty}{\bigl\llVert \mathbf 1_{(s,\varepsilon N^{2/3}
\delta_N)}
\K_N^{(1)}\mathbf 1_{(s,\varepsilon N^{2/3} \delta
_N)}- \mathbf
1_{(s,\infty)}\K_\Aii\mathbf 1_{(s,\infty)}\bigr\rrVert
}_2&=&0,
\\
\lim_{N\rightarrow\infty}\Tr \bigl(\mathbf 1_{(s,\varepsilon N^{2/3}
\delta_N)}
\K_N^{(1)}\mathbf 1_{(s,\varepsilon N^{2/3}\delta
_N)} \bigr)&=&\Tr (
\mathbf 1_{(s,\infty)}\K_\Aii\mathbf  1_{(s,\infty)} ),
\end{eqnarray*}
and, combined moreover with Proposition~\ref{splitprop} and (\ref
{split}), we obtain
%
%e104 #&#
%e105 #&#
\begin{eqnarray}
\label{final1}
&& \lim_{N\rightarrow\infty}{\llVert \mathbf 1_{(s,\varepsilon
N^{2/3}\delta _N)}\widetilde\K_N\mathbf 1_{(s,\varepsilon N^{2/3}\delta
_N)}-\mathbf 1_{(s,\infty)}\K_\Aii\mathbf 1_{(s,\infty)}\rrVert
}_2 = 0,
\\
&& \lim_{N\rightarrow\infty}\Tr (\mathbf 1_{(s,\varepsilon N^{2/3}
\delta_N)}
\widetilde\K_N\mathbf 1_{(s,\varepsilon N^{2/3}\delta
_N)} )
\nonumber\\[-8pt] \label{final2}   \\[-8pt] \nonumber
&&\qquad  =\Tr (
\mathbf 1_{(s,\infty)} \K_\Aii\mathbf  1_{(s,\infty)} ),
\end{eqnarray}
provided we choose $\varepsilon$ small enough. Finally, it follows
from (\ref{changevariables})--(\ref{resc}), (\ref{final1})--(\ref{final2})
and Proposition~\ref{propdet2}
that for every $s\in\R$,
\[
\lim_{N\rightarrow\infty} \p \bigl( N^{2/3}\delta_N (
\tilde x_{\phi(N)}-\frak b_N )\leq s \bigr) =\det (I-
\K_{\Aii} )_{L^2(s,\infty)}.
\]
Proof of Theorem~\ref{thfluctuations-TW}(b) is therefore complete.
\end{pf*}

In the next section, we provide a proof for Theorem~\ref
{thfluctuations-TW}(a) and thus complete the proof for Theorem~\ref
{thfluctuations-TW}. We shall see that we can recover the setting of
the proof of Theorem~\ref{thfluctuations-TW}(b); the only task left is
to prove the existence of appropriate contours for the saddle point
analysis, which differ from the case of a right edge.

%s4.6 #&#
\subsection{Asymptotic analysis for the left edges and proof of Theorem \texorpdfstring{\protect\ref{thfluctuations-TW}\textup{(a)}}{3(a)}}
\label{variationsleftedges}

This section is devoted to the end of the proof of Theorem~\ref
{thfluctuations-TW}. We precisely recall the setting for the analysis
of a left regular soft edge $\frak a$;
we state and prove the counterparts of Proposition~\ref{contoursprop}
(i.e., the existence of appropriate contours for the asymptotic
analysis), that is, Proposition~\ref{contoursleft+} for the case where
$\frak c>0$ with $\frak a = g(\frak c)$, and Proposition~\ref
{contoursleft-} for the case where $\frak c<0$. The remainder of the
asymptotic analysis is omitted since we show it is essentially the same
than in Section~\ref{secasmptotic-analysis-right}.

Let $\frak a$ be a left regular soft edge; recall the definitions of
$g$, $\frak c$, $(\frak c_N)$ as provided by Proposition~\ref
{propproperties-regular-left}, and set
%
%e106 #&#
\begin{equation}
\label{notations*} \frak a_N=g_N(\frak c_N),
\qquad\sigma_N= \biggl(- \frac
{2}{g''_N(\frak c_N)} \biggr)^{1/3}.
\end{equation}
Recall moreover that
%
%e107 #&#
\begin{eqnarray}\label{infocN*}
g'_N(\frak c_N)&=&0,\qquad
\lim_{N\rightarrow\infty}\frak c_N=\frak c,\qquad\lim
_{N\rightarrow\infty} \frak a_N=\frak a,
\nonumber\\[-8pt]\\[-8pt]\nonumber
\lim
_{N\rightarrow\infty}\sigma_N &=& \biggl( -\frac{2}{g''(\frak c)}
\biggr)^{1/3}.
\end{eqnarray}
In particular, for $N$ large enough, $-g_N''(\frak c_N)$ and $\sigma
_N$ are positive numbers, and $\frak c_N$ and $\frak c$ have the same sign.

%s4.6.1 #&#
\subsubsection{Reduction to the right edge setting}
The definition of the extremal eigenvalue $\tilde x_{\varphi(N)}$ (see
Theorem~\ref{thdef-phi-N}) and Proposition~\ref{asymptoticFredholm}
yield that
for every $\varepsilon>0$ small enough,
%
%e108 #&#
\begin{eqnarray}
\label{changevariablesL}
&& \p \bigl(N^{2/3}\sigma_N (\frak
a_{N}-\tilde x_{\varphi(N)} )\leq s \bigr)
\nonumber\\[-8pt]\\[-8pt]\nonumber
&&\qquad = \det (I-\K_N )_{L^2(\frak a_{N}-\varepsilon, \frak a_{N}-s/(N^{2/3}\sigma_N))} +o(1)
\end{eqnarray}
as $N\rightarrow\infty$. We then write
\begin{eqnarray*}
\hspace*{-2pt}&& \det (I-\K_N )_{L^2(\frak a_{N}-\varepsilon, \frak
a_{N}-s/(N^{1/3}\sigma_N))}
=\det (I-\mathbf1_{(s,
N^{2/3}\varepsilon c_N)}{
\widetilde{\K}_N} \mathbf1_{(s,
N^{2/3}\varepsilon c_N)} )_{L^2(s,\infty)},
\end{eqnarray*}
where the scaled operator $\widetilde\K_N$ has for kernel
\[
{\widetilde{\K}_N}(x,y)=-\frac{1}{N^{2/3}\sigma_N}\K_N \biggl(
\frak a_{N}-\frac{x}{N^{2/3}\sigma_N},\frak a_{N}-
\frac{y}{N^{2/3}\sigma
_N} \biggr),
\]
and where $\K_N(x,y)$ was introduced in (\ref{KNversion}) (with
$\frak d_N$ replaced by $\frak c_N$). If we introduce the map
%
%e109 #&#
\begin{equation}
\label{deffN*} f_N^*(z)=\frak a_{N}(z-\frak
c_N)-\log(z)+\frac{1}{N}\sum_{j=1}^n
\log(1-\lambda_j z),
\end{equation}
which differs from $f_N$ defined in (\ref{fN}) by a minus sign and by
the fact that $\frak b_N$ is replaced by $\frak a_N$,
then we have
\begin{eqnarray*}
{\widetilde{\K}_N}(x,y)
&=& -\frac{N^{1/3}}{(2i\pi)^2\sigma_N}
\\[-2pt]
&&\hspace*{6pt}{}\times \oint _{\Gamma} \d z \oint_{\Theta} \d w
\frac{1}{w-z}
\\[-2pt]
&&\hspace*{55pt}{}\times e^{N^{1/3}x(z-\frak c_N)/\sigma_N-N^{1/3}y(w-\frak c_N)/\sigma_N-N
f_N^*(z)+N f_N^*(w)}.
\end{eqnarray*}
Set moreover $\K_N^*(x,y)=\widetilde\K_N(y,x)$. Then it follows by
exchanging $z$ and $w$ in the last integral that
%
%e110 #&#
\begin{eqnarray}\label{kernelleft}
\K_N^*(x,y)
&=& \frac{N^{1/3}}{(2i\pi)^2\sigma_N}\nonumber
\\
&&{}\times \oint_{\Theta} \d z \oint_{\Gamma} \d w
\frac{1}{w-z}
\\
&&\hspace*{49pt}{}\times  e^{-N^{1/3}x(z-\frak
c_N)/\sigma_N+N^{1/3}y(w-\frak c_N)/\sigma_N+N f_N^*(z)-Nf_N^*(w)}.
\nonumber
\end{eqnarray}
Note that, as a consequence of the definition of $f_N^*$ and (\ref
{infocN*}), we have
%
%e111 #&#
\begin{equation}
\label{f*crit} \bigl(f^*_N\bigr)'(\frak
c_N)=\bigl(f^*_N\bigr)''(
\frak c_N)=0,\qquad\bigl(f^*_N\bigr)^{(3)}(\frak
c_N)=-g_N''(\frak
c_N)>0.
\end{equation}
Thus, by comparing (\ref{kernelleft}) with (\ref{kernelright}) and
(\ref{f*crit}) with (\ref{csqkeyidentity}), we recover the setting of
the proof of Theorem~\ref{thfluctuations-TW}(b), except that we
exchange $x$~and~$y$, the role of $\Gamma$ and $\Theta$ as well, and
that we replace $f_N$ by $f_N^*$. Since the Airy kernel is symmetric
[see (\ref{Kai})], it is enough to show that
%
%e112 #&#
%e113 #&#
\begin{eqnarray}
\label{finalleft1} \lim_{N\rightarrow\infty}{\bigl\llVert \mathbf 1_{(s,\varepsilon
N^{2/3}\sigma
_N)}
\bigl(\K_N^*-\K_\Aii\bigr)\mathbf 1_{(s,\varepsilon N^{2/3}\sigma
_N)}\bigr
\rrVert }_2&=&0,
\\
\label{finalleft2} \lim_{N\rightarrow\infty}\Tr \bigl(\mathbf 1_{(s,\varepsilon N^{2/3}
\sigma_N)}
\bigl(\K_N^*-\K_\Aii\bigr)\mathbf 1_{(s,\varepsilon
N^{2/3}\sigma
_N)}
\bigr)&=&0,
\end{eqnarray}
in order to prove (\ref{TWleft}), as explained in the proof of Theorem
\ref{thfluctuations-TW}(b).

In the case of left regular soft edges, the analysis substantially
changes whether~$\frak c$ (cf. Proposition~\ref
{propproperties-regular-left}) is positive or not, and we consider
separately the two cases in the sequel.

%s4.6.2 #&#
\subsubsection{The case where $\frak c$ is positive}\label{cpositive}
We first consider the case where $\frak c>0$, which is always the case,
except if $\frak a$ is the leftmost edge and $\gamma>1$; see
Proposition~\ref{1111}. In particular, $\frak c_N>0$ for all $N$
large enough. We then split $\Gamma$ into two disjoint contours,
$\Gamma^{(0)}$ and $\Gamma^{(1)}$, in the following way: $\Gamma
^{(0)}$ encloses the $\lambda_j^{-1}$'s which are larger that $\frak
c_N$, while $\Gamma^{(1)}$ encloses the $\lambda_j^{-1}$'s which are
smaller that $\frak c_N$. Proposition~\ref{1111}(e), applied to
the measure $\nu_N$, shows that the
set $\{ j, 1\le j\le n\dvtx    \lambda_j^{-1} < \frak c_N \}$ is not
empty, and thus the contour $\Gamma^{(1)}$ is always well defined.
If $\frak c_N$ is actually larger than all the $\lambda_j^{-1}$'s, as
it is
the case when dealing with the smallest eigenvalue when $\gamma<1$,
then set $\Gamma^{(1)}=\Gamma$, $\K_N^{(1)}= \K^*_N$; any later statement
involving $\Gamma^{(0)}$ will be considered empty. Otherwise,
$\Gamma^{(0)}$ is well defined, and we introduce for $\alpha\in\{
0,1\}$ the kernels
\begin{eqnarray*}
\K_{N}^{(\alpha)}(x,y)
&=& \frac{N^{1/3}}{(2i\pi)^2\sigma_N}
\\
&&{}\times \oint _{\Theta} \d z \oint_{\Gamma^{(\alpha)}} \d w
\frac{1}{w-z}
\\
&&\hspace*{59pt}{}\times e^{-
N^{1/3}(z-\frak c_N)x/\sigma_N+N^{1/3}(w-\frak c_N)y/\sigma_N+N
f_N^*(z)-N f_N^*(w)}
\end{eqnarray*}
so that $\K^*_N(x,y)=\K_N^{(0)}(x,y)+\K_N^{(1)}(x,y)$. We similarly
have for the associated operators that $\K_N^*=\K_N^{(0)}+\K
_N^{(1)}$. Observe moreover that we can deform $\Theta$ in $\K
_N^{(1)}(x,y)$ so that it encloses the origin and $\Gamma^{(1)}$ since
the residue we pick at $z=w$ vanishes.

In order to establish (\ref{finalleft1}) and (\ref{finalleft2}), it
is then enough to prove that
%
%e114 #&#
%e115 #&#
\begin{eqnarray}
\label{leftsplit} \lim_{N\rightarrow\infty}{\bigl\llVert \mathbf 1_{(s,\varepsilon
N^{2/3}\sigma
_N)}
\K_N^{(0)}\mathbf 1_{(s,\varepsilon N^{2/3}\sigma_N)}\bigr\rrVert
}_2&=&0,
\\
\lim_{N\rightarrow\infty}\Tr \bigl(\mathbf 1_{(s,\varepsilon N^{2/3}
\sigma_N)}
\K_N^{(0)}\mathbf 1_{(s,\varepsilon N^{2/3}\sigma
_N)} \bigr)&=&0
\end{eqnarray}
and
%
%e116 #&#
%e117 #&#
\begin{eqnarray}
\lim_{N\rightarrow\infty}{\bigl\llVert \mathbf 1_{(s,\varepsilon
N^{2/3}\sigma
_N)}\bigl(
\K_N^{(1)}-\K_\Aii\bigr)\mathbf 1_{(s,\varepsilon N^{2/3}\sigma
_N)}
\bigr\rrVert }_2&=&0,
\\
\label{keyleft} \lim_{N\rightarrow\infty}\Tr \bigl(\mathbf 1_{(s,\varepsilon N^{2/3}
\sigma_N)}
\bigl(\K_N^{(1)}-\K_\Aii\bigr)
\mathbf 1_{(s,\varepsilon
N^{2/3}\sigma
_N)} \bigr)&=&0.
\end{eqnarray}
The exact same estimates as in the proof of the Propositions~\ref
{splitprop}~and~\ref{keyTW} show that (\ref{leftsplit})--(\ref
{keyleft}) hold true, provided we can show the existence of appropriate
contours similar to Proposition~\ref{contoursprop}. More precisely,
it is enough to establish the next proposition in order to prove
Theorem~\ref{thfluctuations-TW}(a), in the case where $\frak c>0$.

%
%pr4.16 #&#
\begin{proposition}
\label{contoursleft+}
For every $\rho>0$ small enough, there exists a contour $\Upsilon^{(0)}$
independent of $N$ and two contours $\Upsilon^{(1)}=\Upsilon
^{(1)}(N)$ and
$\widetilde\Theta=\widetilde\Theta(N)$,
which satisfy for every $N$ large enough the following:
\begin{longlist}[(2)]
\item[(1)]
\begin{enumerate}[(a)]
\item[(a)]
$\Upsilon^{(0)}$ encircles the $\lambda_j^{-1}$'s larger than $\frak c_N$;
\item[\hspace*{31pt}(b)]
$\Upsilon^{(1)}$ encircles the $\lambda_j^{-1}$'s smaller than $\frak c_N$;
\item[\hspace*{31pt}(c)]
$\widetilde\Theta$ encircles the $\lambda_j^{-1}$'s smaller than
$\frak c_N$ and the origin.
\end{enumerate}

\item[(2)]
\begin{enumerate}[(a)]
\item[(a)] $\Upsilon^{(1)}=\Upsilon_*\cup\Upsilon
^{(1)}_{\mathrm{res}}$ where
\[
\Upsilon_*=\bigl\{\frak c_N-te^{\pm i\pi/3}\dvtx  t\in[0,\rho]\bigr\};
\]
\item[\hspace*{31pt}(b)] $\widetilde\Theta=\widetilde\Theta_*\cup
\widetilde\Theta_{\mathrm{res}}$ where
\[
\widetilde\Theta_*=\bigl\{\frak c_N+te^{\pm i\pi/3}\dvtx  t\in[0,\rho]
\bigr\}.
\]
\end{enumerate}

\item[(3)]
There exists $K>0$ independent of $N$ such that:
\begin{enumerate}[(a)]
\item[\hspace*{31pt}(a)] $\re ( f_N(w)-f_N(\frak c_N) )\geq K$ for
all $w\in\Upsilon^{(0)}$;
\item[\hspace*{31pt}(b)] $\re ( f_N(w)-f_N(\frak c_N) )\geq K$ for
all $w\in\Upsilon^{(1)}_{\mathrm{res}}$;\vspace*{1pt}
\item[\hspace*{31pt}(c)] $\re ( f_N(z)-f_N(\frak c_N) )\leq-K$ for
all $z\in\widetilde\Theta_{\mathrm{res}}$.
\end{enumerate}

\item[(4)]
There exists $\d>0$ independent of $N$ such that
\begin{eqnarray*}
\inf \bigl\{\llvert z-w\rrvert\dvtx  z\in\Upsilon^{(0)}, w\in\widetilde
\Theta \bigr\} &\geq& \d,
\\
\inf \bigl\{\llvert z-w\rrvert\dvtx  z\in\Upsilon_{*}, w\in\widetilde
\Theta _{\mathrm{res}} \bigr\} &\geq& \d,
\\
\inf \bigl\{\llvert z-w\rrvert\dvtx  z\in\Upsilon^{(1)}_{\mathrm{res}}, w
\in\widetilde \Theta_{*} \bigr\} &\geq& \d,
\\
\inf \bigl\{\llvert z-w\rrvert\dvtx  z\in\Upsilon^{(1)}_{\mathrm{res}}, w
\in\widetilde \Theta_{\mathrm{res}} \bigr\} &\geq& \d.
\end{eqnarray*}
\item[(5)]
\begin{enumerate}[(a)]
\item[(a)]
The contours $\Upsilon^{(1)}$ and $\widetilde\Theta$ lie in a
compact subset of $\C$, independent \hspace*{31pt}of~$N$.
\item[\hspace*{31pt}(b)]
The lengths of $\Upsilon^{(1)}$ and $\widetilde\Theta$ are uniformly
bounded in $N$.
\end{enumerate}
\end{longlist}
\end{proposition}

\begin{figure}

\includegraphics{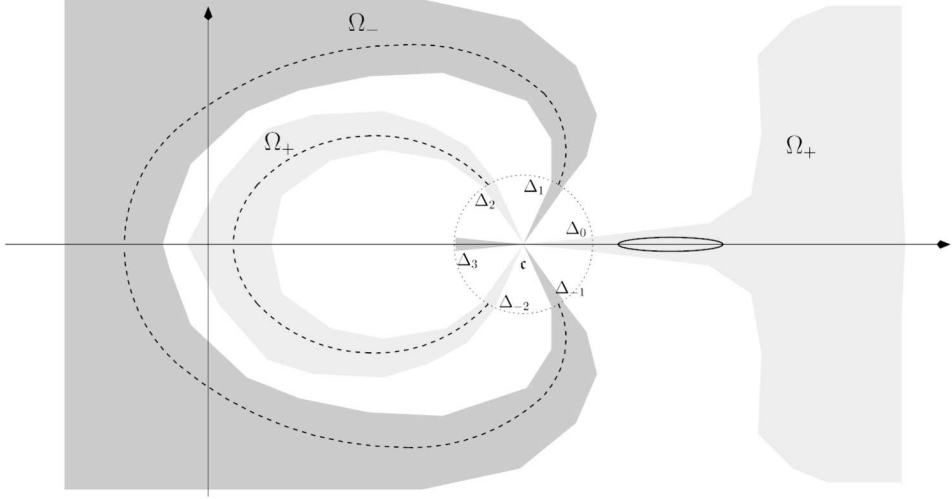}

\caption{Preparation of the saddle point analysis for a left edge with
$\frak c > 0$. The path $\Upsilon_{\mathrm{res}}^{(1)}$ is close to the inner
dotted path
at the left of $\frak c$. The path $\widetilde\Theta_{\mathrm{res}}$ is close
to the outer
dotted path\vspace*{1pt} at the left of $\frak c$. The contour at the right of
$\frak c$ is $\Upsilon^{(0)}$.}
\label{figptselle-c>0}
\end{figure}

Although the proof uses the same type of arguments as in the proof of
Proposition~\ref{contoursprop}, the analytical setting
is not identical. Thus, although we shall provide fewer details than in the
proof of Proposition~\ref{contoursprop}, we shall emphasize the required
changes. Figure~\ref{figptselle-c>0} may help as a visual support for the
argument.

%
%f8 #&#

\begin{pf*}{Proof of Proposition \ref{contoursleft+}}
The regularity assumption yields $\varepsilon>0$ such that $\lambda
_j^{-1}\in(0,+\infty)\setminus B(\frak c,\varepsilon)$ for every
$1\leq j\leq n$ and every $N$ large enough. We then introduce the
compact set $\mathcal K$ defined by
%
%e118 #&#
\begin{equation}
\label{defKleft} {\mathcal K} = \biggl( \biggl[\inf_N
\frac{1}{\lambda_n},\sup_N\frac{1}{\lambda_1} \biggr] \Bigm\backslash
B(\frak c,\varepsilon) \biggr) \cup\{ 0\}
\end{equation}
and notice that by construction $\{x\in\R\dvtx   x^{-1}\in\supp(\nu
_N)\}\subset\mathcal K$ for every $N$ large enough, and also that $\{
x\in\R\dvtx   x^{-1}\in\supp(\nu)\}\subset\mathcal K$. If we
introduce the map
\[
f^*(z)=\frak a(z- \frak c)-\log(z)+\gamma\int\log(1-xz)\nu(\d x),
\]
then, given any simply connected subset of $\C\setminus\mathcal K$,
we can choose a determination of the logarithm such that
both the maps $f^*_N$ and $f^*$ are well defined and holomorphic there
for every $N$ large enough. Notice that the definition of $\re f^*$
does not depend on the determination of the logarithm. Moreover, the
proof of Lemma~\ref{lemmaproperties-f-fN}(a) shows that $\re f_N^*$
converges locally uniformly on $\C\setminus\mathcal K$ toward $\re
f^*$, and moreover $\re f_N^*(\frak c_N)\rightarrow\re f^*(\frak c)$
as $N\rightarrow\infty$.

Next, we perform a qualitative analysis for $\re f^*$ and introduce the sets
\begin{eqnarray*}
\Omega_-& =& \bigl\{z \in\C\dvtx  \re f^*(z)<\re f^*(\frak c) \bigr\},
\\
\Omega_+& =& \bigl\{z \in\C\dvtx  \re f^*(z)>\re f^*(\frak c) \bigr\}.
\end{eqnarray*}
Since $\frak a>0$, the asymptotic behavior $\re f^*(z)=\frak a\re(z-
\frak c)+O(\log\llvert  z\rrvert  )$ as \mbox{$z\to\infty$} shows that for every $\alpha
\in(0,\pi/2)$ there exisits $R>0$ large enough such that
%
%e119 #&#
\begin{equation}
\label{leftplus} \biggl\{z\in\C\dvtx  \llvert z\rrvert >R, -\frac{\pi}{2}+\alpha<
\arg(z)<\frac{\pi
}{2}-\alpha \biggr\}\subset\Omega_+
\end{equation}
and
%
%e120 #&#
\begin{equation}
\label{leftminus} \biggl\{z\in\C\dvtx  \llvert z\rrvert >R, \frac{\pi}{2}+\alpha<
\arg(z)<\frac{3\pi
}{2}-\alpha \biggr\}\subset\Omega_-.
\end{equation}
Notice that the role of $\Omega_+$ and $\Omega_-$ has been exchanged
compared to the setting of a right edge. Moreover, the arguments of the
proof of Lemma~\ref{behaviorinfinity} show that both $\Omega_+$ and
$\Omega_-$ have a unique unbounded connected component.

As for the behavior of $\re f^*$ around $\frak c$, because $\frak
a=g(\frak c)$, it follows from the definition of $f^*$ that
$(f^*)'(z)=g(\frak c)-g(z)$. Thus, by Proposition~\ref
{propproperties-regular-left}, we have
$(f^*)'(\frak c)=(f^*)''(\frak c)=0$ and $ (f^*)^{(3)}(\frak
c)=-g''(\frak c)>0$.
As a consequence, the same proofs as those of Lemmas~\ref{trianglef}~and~\ref{lemmaproperties-f-fN}(b),~(c) show there exist $\eta>0$ and
$0<\theta<\pi/2$ small enough such that
\[
\Delta_{2k+1}\subset\Omega_-,\qquad\Delta_{2k}\subset\Omega
_+,\qquad k\in\{-1,0,1\},
\]
where we introduce, as in Section~\ref{contourssection},
\[
\Delta_k= \biggl\{z\in\C\dvtx  0<\llvert z-\frak c \rrvert <\eta, \biggl
\llvert \arg(z-\frak c)-k\frac{\pi}{3}\biggr\rrvert <\theta \biggr\}.
\]
Notice that the role of $\Omega_-$ and $\Omega_+$ is the same as in
the right edge setting. We then denote by $\Omega_{2k+1}$, the
connected component of $\Omega_-$ which contains $\Delta_{2k+1}$, and
similarly by $\Omega_{2k}$, the connected component of $\Omega_+$
which contains $\Delta_{2k}$.

The proof of Lemma~\ref{supersubharmonic} yields that $\re f^*$ is
subharmonic in $\C\setminus\{0\}$ and is superharmonic in $\C
\setminus\{x\in\R\dvtx   x^{-1}\in\supp(\nu)\}$. As a consequence, it
follows from the proof of Lemma~\ref{keylevelsets} that we obtain a
similar statement as in Lemma~\ref{keylevelsets} for $\re f^*$ after
having exchanged the role of $\Omega_+$ and $\Omega_-$ (to
furthermore convince the reader, notice that $\re f^*(z)-\frak a\re
(z-\frak c)=-\re f(z)-\frak b\re(z-\frak d)$ and that both the maps
$z\mapsto\frak a\re(z-\frak c)$ and $z\mapsto\frak b\re(z-\frak d)$
are harmonic). Namely:
\begin{enumerate}[(2)]
\item[(1)]
If $\Omega_*$ is a connected component of $\Omega_-$, then $\Omega
_*$ is open, and if $\Omega_*$ is moreover bounded, there exists $x\in
\supp(\nu)$ such that $x^{-1}\in\Omega_*$.
\item[(2)]
Let $\Omega_*$ be a connected component of $\Omega_+$ with nonempty interior:
\begin{longlist}[(a)]
\item[(a)]
if $\Omega_*$ is bounded, then $0\in\Omega_*$;
\item[(b)]
if $\Omega_*$ is bounded, then its interior is connected;
\item[(c)]
if $0\notin\Omega_*$, then the interior of $\Omega_*$ is connected.
\end{longlist}
\end{enumerate}
Equipped with the previous observations we are now in position to
provide the counterpart of Lemma~\ref{studylevelsets} in the present
setting, namely to prove that the following statements hold true:
\begin{longlist}[(A)]
\item[(A)] we have $\Omega_1=\Omega_{-1}$, the interior of $\Omega
_1$ is connected, and for every $0<\alpha<\pi/2$ there exists $R>0$
such that
\[
\biggl\{z\in\C\dvtx  \llvert z\rrvert >R, \frac{\pi}{2}+\alpha<\arg(z)<
\frac{3\pi
}{2}-\alpha \biggr\}\subset\Omega_1;
\]
\item[(B)]
the interior of $\Omega_0$ is connected, and for every $0<\alpha<\pi
/2$, there exists $R>0$ such that
\[
\biggl\{z\in\C\dvtx  \llvert z\rrvert >R, -\frac{\pi}{2}+\alpha<\arg(z)<
\frac{\pi
}{2}-\alpha \biggr\}\subset\Omega_0;
\]
\item[(C)]
we have $\Omega_2=\Omega_{-2}$, the interior of $\Omega_2$ is
connected, and there exists $\delta>0$ such that $B(0,\delta)\subset
\Omega_2 $.
\end{longlist}

Let us first prove (A). Since by definition, $\Omega_1$ is a connected
component of $\Omega_-$, its interior is connected by (1). Let us
prove by contradiction that $\Omega_1$ is unbounded, from which (A)
will follow by using the symmetry $\re f^*(\overline z)=\re f^*(z)$,
inclusion~(\ref{leftminus}) and that $\Omega_-$ has a unique
unbounded connected component. Assume $\Omega_1$ is bounded. Then (1)
yields the existence of $x\in\supp(\nu)$ such that $x^{-1}\in\Omega
_1$. If $x^{-1}<\frak c$ (resp., $x^{-1}>\frak c$), it then follows from
the symmetry $\re f^*(\overline z)=\re f^*(z)$ that $\Omega_1$
surrounds $\Omega_2$ (resp., $\Omega_0$) so that $\Omega_2*$
(resp., $\Omega_0$)
is a bounded connected component of $\Omega_+$ which does not contain
the origin. Notice that by (\ref{leftplus}), $\Omega_2$ (resp.,
$\Omega_0$) has a nonempty interior.
This yields, with (2)(a), a contradiction, and our claim follows.
Since we just proved that $\Omega_1$ is unbounded, the origin does not
belong to $\Omega_0$. As a consequence, (2)(a) and (2)(c) yield,
respectively, that $\Omega_0$ is unbounded and has a connected
interior. Using moreover inclusion (\ref{leftplus}) and that $\Omega
_+$ has a unique unbounded connected component, (B) follows.

As a byproduct of (A), $\Omega_2$ is bounded. Thus $\Omega_2$
contains the origin by (2)(a) and has a connected interior by (2)(b). By
using the symmetry $\re f^*(\overline z)=\re f^*(z)$ and that $\re
f^*(z)\to+\infty$ as $z\to0$, (C) is proved.

Finally,\vspace*{1pt} as a consequence of (A), (B) and (C), the existence of the
contour $\Upsilon^{(0)}$, respectively, $\Upsilon^{(1)}$, respectively, $\widetilde
\Theta$, in Proposition~\ref{contoursleft+} is proved by choosing
$\Upsilon^{(0)}$ in the interior of $\Omega_0$ encircling $\{x\in
\mathcal K\dvtx  x>\frak c\}$ and intersecting the real axis exactly twice
in $\R\setminus\mathcal K$ with finite length, respectively, by
completing $\{\frak c_N-te^{\pm i\pi/3}\dvtx   t\in[0,\rho]\}$ for $\rho
$ small enough and $N$ large enough so that both the points $\frak
c_N-\rho e^{ i\pi/3}$ and $\frak c_N-\rho e^{ -i\pi/3}$ lie in
$\Omega_2$ into a closed contour with a path lying in the interior of
$\Omega_2$ but staying in $\{z\in\C\dvtx    \re(z)>0\}$ and
intersecting the real line exactly once at the left of $\mathcal K$
with finite length, respectively, by completing $\{\frak c_N+te^{\pm
i\pi/3}\dvtx   t\in[0,\rho]\}$ for $\rho$ small enough and $N$ large
enough so that both the points $\frak c_N+\rho e^{ i\pi/3}$ and $\frak
c_N+\rho e^{ -i\pi/3}$ belong to $\Omega_1$ into a closed contour
with a path lying in the interior of $\Omega_1$ and crossing the real
axis exactly once at the left of the origin with finite length, and
then by using the local uniform convergence of $\re f_N^*\rightarrow
\re f^*$ on $\C\setminus\mathcal K$; see the proof of Proposition
\ref{contoursprop} for the details.
\end{pf*}

%s4.6.3 #&#
\subsubsection{The case where $\frak c$ is negative}
\label{cnegativesection}
Here we consider the case where $\frak c$ is negative, which only
happens if we are looking at the leftmost edge $\frak a$ when $\gamma
>1$, and thus $\frak c_N<0$ for all $N$ large enough. We recall that
\begin{eqnarray*}
\K_N^*(x,y)
&=& \frac{N^{1/3}}{(2i\pi)^2\sigma_N}
\\
&&{}\times \oint_{\Theta} \d z \oint _{\Gamma} \d w
\frac{1}{w-z}
\\
&&\hspace*{49pt}{}\times  e^{-N^{1/3}x(z-\frak c_N)/\sigma
_N+N^{1/3}y(w-\frak c_N)/\sigma_N+N f_N^*(z)-Nf_N^*(w)}.
\end{eqnarray*}
Note that the $\lambda_j^{-1}$'s are zeros for $e^{f_N^*}$, and that
$0$ is a zero for $e^{-f_N^*}$. Thus, since the residue picked at $w=z$
vanishes, we can deform $\Theta$ and $\Gamma$ in a way that $\Gamma$
encircles $\Theta$ and all the $\lambda_j^{-1}$'s, whereas $\Theta$
encircles the origin and possibly some~$\lambda_j^{-1}$'s.

\begin{figure}

\includegraphics{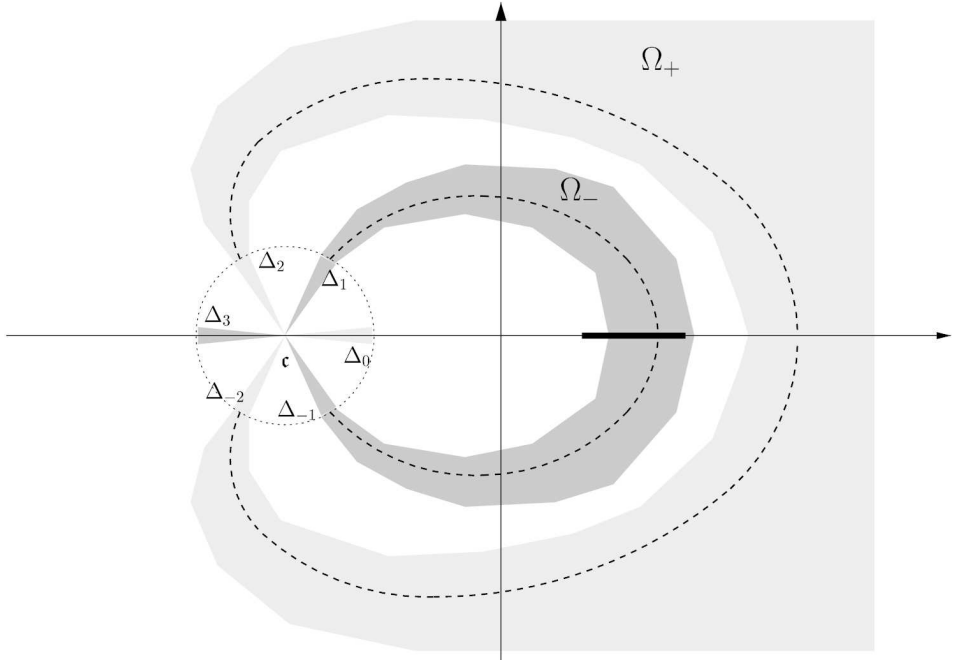}

\caption{Preparation of the saddle point analysis for a left edge with
$\frak c < 0$. The path $\widetilde\Theta_{\mathrm{res}}$ is close to the
inner dotted
path. The path $\widetilde\Upsilon_{\mathrm{res}}$ is close to the outer
dotted path.
The thick segment represents the support of the image of $\nu$ by
the map $x\mapsto x^{-1}$.}
\label{figleft-c<0}
\end{figure}

It is enough to establish the next proposition in order to obtain (\ref
{finalleft1}) and (\ref{finalleft2}) in the case where $\frak c<0$,
and thus to complete the proof of Theorem~\ref{thfluctuations-TW}(a),
since the same estimates as in the proof of Proposition~\ref{keyTW}
can be used after setting $\K^{(1)}_N=\K^*_N$ and $\Gamma
^{(1)}=\Gamma$.
The reader may refer to Figure~\ref{figleft-c<0} to better visualize the
results of the next proposition as well as the proof argument.

%
%pr4.17 #&#
\begin{proposition}
\label{contoursleft-}
For every $\rho>0$ small enough, there exist contours $\Upsilon
=\Upsilon(N)$ and $\widetilde\Theta=\widetilde\Theta(N)$ which
satisfy for every $N$ large enough the following:
\begin{longlist}[(2)]
\item[(1)]
\begin{enumerate}[(a)]
\item[(a)]
$\Upsilon$ encircles $\widetilde\Theta$, the origin and all the
$\lambda_j^{-1}$'s;
\item[\hspace*{31pt}(b)]
$\widetilde\Theta$ encircles the origin (and possibly some $\lambda
_j^{-1}$'s).
\end{enumerate}

\item[(2)]
\begin{enumerate}[(a)]
\item[(a)] $\Upsilon=\Upsilon_*\cup\Upsilon_{\mathrm{res}}$ where
\[
\Upsilon_*=\bigl\{\frak c_N-te^{\pm i\pi/3}\dvtx  t\in[0,\rho]\bigr\};
\]
\item[\hspace*{31pt}(b)] $\widetilde\Theta=\widetilde\Theta_*\cup
\widetilde\Theta_{\mathrm{res}}$ where
\[
\widetilde\Theta_*=\bigl\{\frak c_N+te^{\pm i\pi/3}\dvtx  t\in[0,\rho]
\bigr\}.
\]
\end{enumerate}
\item[(3)]
There exists $K>0$ independent of $N$ such that:
\begin{enumerate}[(a)]
\item[\hspace*{31pt}(a)] $\re ( f_N(w)-f_N(\frak c_N) )\geq K$ for
all $w\in\Upsilon_{\mathrm{res}}$;
\item[\hspace*{31pt}(b)] $\re ( f_N(z)-f_N(\frak c_N) )\leq-K$ for
all $z\in\widetilde\Theta_{\mathrm{res}}$.
\end{enumerate}

\item[(4)]
There exists $\d>0$ independent of $N$ such that
\begin{eqnarray*}
\inf \bigl\{\llvert z-w\rrvert\dvtx  z\in\Upsilon_{*}, w\in\widetilde
\Theta _{\mathrm{res}} \bigr\} &\geq& \d,
\\
\inf \bigl\{\llvert z-w\rrvert\dvtx  z\in\Upsilon_{\mathrm{res}}, w\in\widetilde
\Theta _{*} \bigr\} &\geq& \d,
\\
\inf \bigl\{\llvert z-w\rrvert\dvtx  z\in\Upsilon_{\mathrm{res}}, w\in\widetilde
\Theta _{\mathrm{res}} \bigr\}&\geq& \d.
\end{eqnarray*}
\item[(5)]
\begin{enumerate}[(a)]
\item[(a)]
$\Upsilon$ and $\widetilde\Theta$ lie in a bounded subset of $\C$
independently of $N$;
\item[\hspace*{31pt}(b)]
the lengths of $\Upsilon$ and $\widetilde\Theta$ are uniformly
bounded in $N$.
\end{enumerate}
\end{longlist}
\end{proposition}

%
%f9 #&#

\begin{pf} We use the notation, definitions and properties used in the
proof of Proposition~\ref{contoursleft+}, except for $\mathcal K$ that
we define by
\[
{\mathcal K} = \biggl[\inf_N \frac{1}{\lambda_n},\sup
_N\frac
{1}{\lambda_1} \biggr] \cup\{ 0\}.
\]
Clearly $\{x\in\R\dvtx   x^{-1}\in\supp(\nu_N)\}\subset\mathcal K$
for every $N$ and moreover $\{x\in\R\dvtx   x^{-1}\in\supp(\nu)\}
\subset\mathcal K$.
We now prove that the following facts hold true:
\begin{longlist}[(A)]
\item[(A)] we have $\Omega_1=\Omega_{-1}$, the interior of $\Omega
_1$ is connected and there exists $x_0\in\supp(\nu)$ and $\delta>0$
such that $B(x_0^{-1},\delta)\subset\Omega_1$;
\item[(B)]
we have $\Omega_2=\Omega_{-2}$, the interior of $\Omega_2$ is
connected and for every $0<\alpha<\pi/2$, there exists $R>0$ such that
\[
\biggl\{z\in\C\dvtx  \llvert z\rrvert >R, -\frac{\pi}{2}+\alpha<\arg(z)<
\frac{\pi
}{2}-\alpha \biggr\}\subset\Omega_2.
\]
\end{longlist}

The proof will mainly use properties (1) and (2)(a)/(b)/(c) from the
proof of Proposition~\ref{contoursleft+}. Let us show (A). First,
$\Omega_1$ has a connected interior by (1). Let us show by
contradiction that $\Omega_1$ is bounded. If $\Omega_1$ is unbounded,
then by using the symmetry $\re f^*(\overline z)=\re f^*(z)$, inclusion
(\ref{leftminus}) and the uniqueness of the unbounded connected
component of $\Omega_-$, it follows that $\Omega_2$ is bounded
without containing the origin, which contradicts (2)(a). Thus $\Omega_1$
is bounded and has to contain some $x_0^{-1}$ with $x_0\in\supp(\nu
)$ as a consequence of (1). Moreover, since $\re f^*$ is upper
semicontinuous on an open neighborhood of $x_0^{-1}$ (because it is
subharmonic on $\C\setminus\{0\}$), there exists $\delta>0$ such
that $B(x_0^{-1},\delta)\subset\Omega_1$. As a consequence, together
with the symmetry $\re f^*(\overline z)=\re f^*(z)$, (A) is proved.

Next, since $\Omega_1$ thus surrounds the origin, then $\Omega_2$ has
to be unbounded by (2)(a) and has a connected interior by (2)(c). Finally,
(B) follows from the symmetry $\re f^*(\overline z)=\re f^*(z)$,
inclusion (\ref{leftplus}) and the uniqueness of the unbounded
connected component of $\Omega_+$.

To construct $\Upsilon$ satisfying the conditions of Proposition~\ref
{contoursleft-}, by (B) we can complete $\{\frak c_N-te^{\pm i\pi
/3}\dvtx   t\in[0,\rho]\}$, for $N$ large enough and $\rho$ small enough
so that both the points $\frak c_N-\rho e^{ i\pi/3}$ and $\frak
c_N-\rho e^{ -i\pi/3}$ lie in $\Omega_2$, into a closed contour with
a path lying in the interior of $\Omega_2$ and intersecting the real
line exactly once at the right of $\mathcal K$ with finite length, and
then use the local uniform convergence of $\re f_N^*$ to $ \re f^*$ on
$\C\setminus\mathcal K$; see the proof of Proposition~\ref
{contoursprop} for the details.

To construct $\widetilde\Theta$, we need to proceed more carefully
since $\Omega_1$ actually crosses~$\mathcal K$, and $\re f^*_N$ may
not converge uniformly to $\re f^*$ there. For $N$ large enough and
$\rho$ small enough so that the points $\frak c_N+\rho e^{ i\pi/3}$
and $\frak c_N+\rho e^{ -i\pi/3}$ lie in $\Omega_1$, by~(A) we can
complete $\{\frak c_N+te^{\pm i\pi/3}\dvtx   t\in[0,\rho]\}$ into a
closed contour with a path $\Xi$ lying in the interior of $\Omega_1$
and crossing the real axis exactly once at $x_0^{-1}$ with finite
length. Since $B(x_0^{-1},\delta)\subset\Omega_1$ we can moreover
assume that $\Xi$ crosses the real axis perpendicularly, namely that
there exists $\eta_1>0$ small enough such that the segment $\{
x_0^{-1}+i\eta\dvtx   \llvert  \eta\rrvert  \leq\eta_1\}$ is contained in $\Xi$. Since
$\Omega_1\subset\Omega_-$, there exists $K>0$ independent on $N$
such that
%
%e121 #&#
\begin{equation}
\label{ineqXi} \re f^*(z)-\re f^*(\frak c)\leq-4K,\qquad z\in\Xi.
\end{equation}

Notice that the map $z\mapsto\int\log\llvert  1-xz\rrvert  \nu(\d x)$ is upper
semicontinuous on $\C$ since it is subharmonic; see the proof of Lemma
\ref{supersubharmonic}. As a consequence, if $\int\log\llvert  1-x
x_0^{-1}\rrvert  \nu(\d x)=-\infty$, then there exists $\eta_0\in(0,\eta
_1)$ small enough so that
\begin{eqnarray}
\label{etainftycond}
&& \gamma\int\log\bigl\llvert 1-x \bigl(x_0^{-1}+i
\eta_0\bigr)\bigr\rrvert \nu(\d x)
\nonumber\\[-8pt]\\[-8pt]\nonumber
&&\qquad \leq-2K-\sup_N \bigl(\frak a_N
\bigl(x_0^{-1}-\frak c_N\bigr)-\re f^*(\frak
c_N) \bigr)-\log(x_0).
\nonumber
\end{eqnarray}
If instead $\int\log\llvert  1-x/x_0\rrvert  \nu(\d x)>-\infty$, then by upper
semicontinuity there exists $\eta_0\in(0,\eta_1)$ small enough to that
%
%e122 #&#
\begin{equation}
\label{eta>cond} \qquad\gamma\int\log\bigl\llvert 1-x\bigl(x_0^{-1}+i
\eta_0\bigr)\bigr\rrvert \nu(\d x)\leq\gamma\int \log\llvert
1-x/x_0\rrvert \nu(\d x)+K.
\end{equation}

Let $\eta_0$ be defined as above, and consider a compact tubular
neighborhood $\mathcal T$ of $\Xi\setminus\{x_0^{-1}+i\eta\dvtx   \llvert  \eta
\rrvert  < \eta_0\}$ small enough so that $\mathcal T$ lies in $\C\setminus
\mathcal K$ and $\re f^*-\re f^*(\frak c)\leq-3K$ there (the latter is
possible since $\re f^*$ is upper semicontinuous on $\C\setminus\{0\}
$). Notice that by construction the interior of $\mathcal T$ contains
both the points $\frak c_N+\rho e^{ i\pi/3}$ and $\frak c_N+\rho e^{
-i\pi/3}$ for every $N$ large enough, and the points $x_0^{-1}+i\eta
_0$ and $x_0^{-1}-i\eta_0$ as well. Using the local uniform
convergence of $\re f_N^*$ to $\re f^*$ on $\C\setminus\mathcal K$
and the convergence $\re f^*_N(\frak c_N)\to\re f(\frak c)$, we can
show as in the proof of Proposition~\ref{contoursprop} that for every
$N$ large enough, we have
\[
\re f_N^*(z)-\re f_N^*(\frak c_N)\leq-K
\]
for every $z\in\mathcal T$. As a consequence, for every $N$ large
enough, we can construct the path $\widetilde\Theta_{\mathrm{res}}$ in the
following way: it goes from $\frak c_N+\rho e^{ -i\pi/3}$ to
$x_0^{-1}+i\eta_0$ staying in $\mathcal T$, then follows the segment
$\{x_0^{-1}+i\eta\dvtx   0\leq\eta\leq\eta_0\}$, and is finally
completed by symmetry with respect to the real axis. As for what is
happening on $ \{x_0^{-1}+i\eta\dvtx   \llvert  \eta\rrvert  < \eta_0\}$, since a priori
$\re f_N^*$ does not converge uniformly there toward $\re f^*$, we need
an extra argument to complete the proof of Proposition~\ref
{contoursleft-}. Namely, we need to show that for every $N$ large
enough, uniformly in $\llvert  \eta\rrvert  <\eta_0$,
%
%e123 #&#
\begin{equation}
\label{toprovec<0} \re f^*_N\bigl(x_0^{-1}+i
\eta\bigr)-\re f_N^*(\frak c_N)\leq-K.
\end{equation}

Let us set for convenience $z_\eta=x_0^{-1}+i\eta$ for any $\llvert  \eta
\rrvert  \leq\eta_0$. First, since the map $x\mapsto\log\llvert  1-x z_{\eta_0}\rrvert  $
is bounded and continuous on any compact subset of $\R$, the weak
convergence $\nu_N\to\nu$ and the convergence $n/N\to\gamma$ yield
that for any $N$ large enough,
%
%e124 #&#
\begin{equation}
\label{weakconvlog} \frac{n}N\int\log\llvert 1-xz_{\eta_0}\rrvert
\nu_N(\d x)\leq\gamma\int\log \llvert 1-xz_{\eta_0}\rrvert \nu(
\d x) +K.
\end{equation}

If we assume $\int\log\llvert  1-x / x_0\rrvert  \nu(\d x)=-\infty$, then for every
$N$ large enough, uniformly in $\llvert  \eta\rrvert  <\eta_0$,
\begin{eqnarray*}
\label{ineqkeyc<00} && \re f^*_N(z_\eta)-\re
f_N^*(\frak c_N)
\\
&&\qquad\leq\sup_N \bigl(\frak a_N
\bigl(x_0^{-1}-\frak c_N\bigr)-\re f^*(\frak
c_N) \bigr)-\log\llvert z_\eta\rrvert +\frac{n}N
\int\log\llvert 1-xz_\eta\rrvert \nu_N(\d x)
\nonumber
\\
&&\qquad\leq\sup_N \bigl(\frak a_N
\bigl(x_0^{-1}-\frak c_N\bigr)-\re f^*(\frak
c_N) \bigr)+\log(x_0)+\frac{n}N\int\log
\llvert 1-xz_{\eta_0}\rrvert \nu_N(\d x)
\nonumber
\\
&&\qquad\leq-K,
\nonumber
\end{eqnarray*}
%
%{\red[j'ai du mal {\chr"C3}{\chr"A0} reconstituer ces in{\chr"C3}{
%\chr"A9}galit{\chr"C3}{\chr"A9}s..]}
where for the last inequality we use (\ref{weakconvlog}) and (\ref
{etainftycond}).

Now, assume instead that $\int\log\llvert  1-x /x_0\rrvert  \nu(\d x)>-\infty$. By
using the convergences $\frak a_N\to\frak a$, $\frak c_N\to\frak c$
and $\re f_N^*(\frak c_N)\to\re f^*(\frak c)$, we obtain for every $N$
large enough (and independently on $\eta$)
\[
\frak a_N\re(z_\eta-\frak c_N)-\re
f_N^*(\frak c_N)\leq\frak a\re (z_\eta-\frak
c)-\re f^*(\frak c)+K.
\]
Combined with inequalities (\ref{ineqXi}), (\ref{weakconvlog}) and
(\ref{eta>cond}), we obtain that for every $N$ large enough and
uniformly in $\llvert  \eta\rrvert  <\eta_0$,
\begin{eqnarray*}
\label{ineqkeyc<0} && \re f^*_N(z_\eta)-\re
f_N^*(\frak c_N)
\\
&&\qquad\leq K + \re f^*(z_\eta)-\re f^*(\frak c)
\\
&&\quad\qquad{}+\frac{n}N
\int\log \llvert 1-xz_\eta\rrvert \nu_N(\d x)-\gamma\int
\log\llvert 1-xz_\eta\rrvert \nu(\d x)
\nonumber
\\
&&\qquad\leq-3K+\frac{n}N\int\log\llvert 1-xz_\eta\rrvert
\nu_N(\d x)-\gamma \int\log\llvert 1-xz_\eta\rrvert \nu(\d
x)
\nonumber
\\
&&\qquad\leq-3K+\frac{n}N\int\log\llvert 1-xz_{\eta_0}\rrvert
\nu_N(\d x)-\gamma\int\log\llvert 1-x/x_0\rrvert \nu(\d
x)
\nonumber
\\
&&\qquad\leq-2K+ \gamma\int\log\llvert 1-xz_{\eta_0}\rrvert \nu(\d x)-
\gamma\int\log\llvert 1-x/x_0\rrvert \nu(\d x)
\nonumber
\\
&&\qquad\leq-K,
\nonumber
\end{eqnarray*}
and this completes the proof of Proposition~\ref{contoursleft-}.
\end{pf}

%s5 #&#
\section{Proof of Theorem \texorpdfstring{\protect\ref{thasymptindep}}{4}: Asymptotic independence}\label{secasymptotic-independence}

Our strategy to prove Theorem~\ref{thasymptindep} builds on an
approach used by Bornemann \cite{bornemann-2010}. Indeed, the
asymptotic independence for the smallest and largest eigenvalues of an
$N\times N$ GUE random matrix is established in \cite{bornemann-2010}
by showing that the trace class norm of the \mbox{off-}diagonal entries of a
two-by-two operator valued matrix goes to zero as $N\rightarrow\infty
$. Here we obtain that proving the asymptotic joint independence of
several extremal eigenvalues leads to considering a larger operator
valued matrix. Moreover, we show that it is actually sufficient to
establish that the Hilbert--Schmidt norms of the off-diagonal entries
go to zero as $N\rightarrow\infty$, instead of the trace class norms.
The former can be provided by an asymptotic analysis for double complex
integrals as we performed in the previous section.

More generally, our method can be applied to several other
determinantal point processes for which a contour integral
representation for the kernel and its asymptotic analysis are known,
for example, the eigenvalues of an additive perturbation of a GUE
matrix \cite{capitaine-peche-2014-preprint}.

\begin{longlist}[]
\item[\textit{Conventions}.]
In this section, we fix two finite sets $I$ and $J$ of indices, and real
numbers $(s_i)_{i\in I}$ and $(t_j)_{j\in J}$ as well.
Assume that $(\frak a_i = g(\frak c_i))_{i\in I}$ are regular left soft edges
and $(\frak b_j = g(\frak d_j))_{j\in J}$ are regular right edges. We denote
by $\frak c_{i,N}$ and $\frak d_{j,N}$ the sequences associated, respectively,
with $\frak a_i$ and $\frak b_j$ as specified by Proposition~\ref{gN->g}(c).
We moreover set
\[
\frak a_{i,N}=g_N(\frak c_{i,N}),\qquad\frak
b_{j,N}=g_N(\frak d_{j,N})
\]
and
\[
\sigma_{i,N}= \biggl(\frac{2}{g''_N(\frak c_{i,N})} \biggr)^{1/3},\qquad
\delta_{j,N}= \biggl(\frac{2}{g''_N(\frak d_{j,N})} \biggr)^{1/3},
\]
where $g_N$ has been introduced in (\ref{gN}). Similarly, $\varphi_i(N)$
[resp., $\phi_j(N)$] denotes the sequence associated with
$\frak a_{i,N}$ (resp., $\frak b_{j,N}$) as in Theorem~\ref
{thdef-phi-N}; see also Propositions~\ref{propproperties-regular-left}
and~\ref{propproperties-regular-right}.

Finally, we shall consider that the free parameter $q$ introduced in
the statement of Proposition~\ref{propdet-representation} is zero when
dealing with
the kernel $\K_N(x,y)$; see Remark~\ref{freedomq}.

Our starting point is the following proposition.
\end{longlist}

%
%pr5.1 #&#
\begin{proposition}
\label{multiasymptoticFredholm}
Consider the setting of Theorem~\ref{thasymptindep}. Then, for every
\mbox{$\varepsilon>0$} small enough and for every sequences $(\eta
_{i,N})_N$, $(\chi_{j,N})_N$ of positive numbers growing with $N$ to
infinity, it holds that
\begin{eqnarray*}
&& \p \bigl( \eta_{i,N} (\frak a_{i,N}-x_{\varphi_i(N)} )
\leq s_i, \chi_{j,N} (x_{\phi_j(N)}-\frak
b_{j,N} )\leq t_j, i\in I, j\in J \bigr)
\\
&&\qquad= \det ( I -\K_N )_{L^2 ((\bigcup_{i\in
I}A_i)
\cup (\bigcup_{j\in J}B_j) )} + o(1)
\end{eqnarray*}
as $N\rightarrow\infty$, where
\[
A_i= \biggl(\frak a_{i,N}-\varepsilon, \frak
a_{i,N}-\frac{s_i}{\eta
_{i,N}} \biggr),\qquad B_j= \biggl(
\frak b_{j,N}+\frac{t_j}{\chi
_{j,N}}, \frak b_{j,N}+\varepsilon
\biggr).
\]
\end{proposition}

The proof is omitted, being very similar to that of Proposition~\ref
{asymptoticFredholm}. Now, if we specify $\eta_{i,N}=N^{2/3}\sigma
_{i,N}$ and $\chi_{j,N}=N^{2/3}\delta_{j,N}$, then Proposition~\ref
{multiasymptoticFredholm} reads
%
%e125 #&#
\begin{eqnarray}
\label{AIfredrep} && \p \bigl( N^{2/3}\sigma_{i,N} (\frak
a_{i,N}-x_{\varphi_i(N)} )\leq s_i,\nonumber
\\
&&\hspace*{11pt}  N^{2/3}
\delta_{j,N} (x_{\phi_j(N)}-\frak b_{j,N} )\leq
t_j, i\in I, j\in J \bigr)
\\
&&\qquad= \det ( I -\K_N )_{L^2 ((\bigcup_{i\in
I}A_i)
\cup (\bigcup_{j\in J}B_j) )} + o(1),\nonumber
\end{eqnarray}
where
\[
A_i= \biggl(\frak a_{i,N}-\varepsilon, \frak
a_{i,N}-\frac
{s_i}{N^{2/3}\sigma_{i,N}} \biggr),\qquad B_j= \biggl(
\frak b_{j,N}+\frac{t_j}{N^{2/3}\delta_{j,N}}, \frak b_{j,N}+\varepsilon
\biggr).
\]

For every $i\in I$ and $j\in J$, we introduce the maps
%
%e126 #&#
%e127 #&#
\begin{eqnarray}
\label{fj*} f^*_{i,N}(z) & =&\frak a_{i,N}(z-\frak
c_{i,N})-\log(z)+\frac
{1}{N}\sum_{k=1}^n
\log(1-\lambda_k z),
\\
f_{j,N}(z) & =&-\frak b_{j,N}(z-\frak d_{j,N})+
\log(z)-\frac
{1}{N}\sum_{k=1}^n
\log(1-\lambda_k z)
\end{eqnarray}
and the multiplication operators $\E_i^*$ and $\E_j$ acting on
$L^2(A_i)$ and $L^2(B_j)$, respectively, by
\begin{eqnarray*}
\E_i^*h(x) & =&e^{Nf_{i,N}^*(\frak c_{i,N})+Nx\frak
c_{i,N}}h(x),\qquad h\in
L^2(A_i),
\\
\E_jh(x) & =&e^{-Nf_{j,N}(\frak d_{j,N})+Nx\frak d_{j,N}}h(x),\qquad h\in
L^2(B_j).
\end{eqnarray*}

The next proposition is the key to obtain Theorem~\ref{thasymptindep}.

%
%pr5.2 #&#
\begin{proposition}
\label{AIkey}
For every $\varepsilon$ small enough, the following holds true:
\begin{longlist}[(a)]
\item[(a)]
for every $(i,j)\in J\times J$ such that $i\neq j$, we have
%
%e128 #&#
\begin{equation}
\lim_{N\rightarrow\infty}{\bigl\llVert \mathbf 1_{B_i}
\E_i \K_N \E _j^{-1}
\mathbf 1_{B_j}\bigr\rrVert }_2=0;
\end{equation}

\item[(b)]
for every $(i,j)\in I\times I$ such that $i\neq j$, we have
%
%e129 #&#
\begin{equation}
\lim_{N\rightarrow\infty}{\bigl\llVert \mathbf 1_{A_i}
\E_i^* \K_N \bigl(\E _j^*
\bigr)^{-1}\mathbf 1_{A_j}\bigr\rrVert }_2=0;
\end{equation}

\item[(c)] for every $(i,j)\in I\times J$, we have
%
%e130 #&#
\begin{equation}
\lim_{N\rightarrow\infty}{\bigl\llVert \mathbf 1_{A_i}
\E_i^* \K_N \E _j^{-1}
\mathbf 1_{B_j}\bigr\rrVert }_2=0
\end{equation}
and
%
%e131 #&#
\begin{equation}
\lim_{N\rightarrow\infty}{\bigl\llVert \mathbf 1_{B_j}
\E_j \K_N \bigl(\E _i^*
\bigr)^{-1}\mathbf 1_{A_j}\bigr\rrVert }_2=0.
\end{equation}
\end{longlist}
\end{proposition}

Before proving Proposition~\ref{AIkey}, let us show how does it lead
to the
asymptotic joint independence of the extremal eigenvalues:

\begin{pf*}{Proof of Theorem~\ref{thasymptindep}}
Our purpose is to show that for large $N$, the determinant at the right-hand
side of~(\ref{AIfredrep}) converges to a product of Fredholm
determinants involving
the Airy kernel.
Assume that $N$ is large enough so that all the $A_i$'s and $B_j$'s are disjoint
sets. Then, as shown in \cite{bor-mc-10}
(see also~\cite{book-gohberg-2000}, Chapter~6), the Fredholm determinant
$\det ( I -\K_N )_{L^2 ((\bigcup_{i\in I}A_i) \cup
(\bigcup_{j\in J}B_j) )}$
admits the operator matrix representation
%
%e132 #&#
\begin{eqnarray}
\label{AIstep1} && \det ( I -\K_N )_{L^2 ((\bigcup_{i\in I}A_i) \cup
(\bigcup_{j\in J}B_j) )}
\nonumber\\[-2pt]\\[-14pt]\nonumber
&&\qquad= \det\lleft(I-\lleft[\matrix{ \bigl[ \K_{II}^{i,j}
\bigr]_{(i,j) \in I \times I} & \bigl[ \K_{IJ}^{i,j}
\bigr]_{(i,j) \in I \times J}
\cr
\bigl[ \K_{JI}^{i,j}
\bigr]_{(i,j) \in J \times I} & \bigl[ \K_{JJ}^{i,j}
\bigr]_{(i,j) \in J \times J} } \rright] \rright)_{ (\bigoplus_{i\in I}L^2(A_i) ) \oplus
 (\bigoplus_{j\in J} L^2(B_j) )},
\end{eqnarray}
where $\K_{II}^{i,j}\dvtx  L^2(A_j) \to L^2(A_i)$ denotes the integral operator
\[
\K_{II}^{i,j} h (x) = \int_{A_j}
\K_N(x,y) h(y)\,\d y, \qquad x \in A_i,
\]
and similarly the operators $\K_{IJ}^{i,j}\dvtx  L^2(B_j) \to L^2(A_i)$,
$\K_{JI}^{i,j}\dvtx  L^2(A_j) \to L^2(B_i)$ and
$\K_{JJ}^{i,j}\dvtx  L^2(B_j) \to L^2(B_i)$ are defined by restricting $\K
_N$ on appropriate subspaces of $L^2(\R)$. Consider now the diagonal operator
\[
\E= \biggl(\bigoplus_{i\in I} \E_i^*
\biggr)\oplus \biggl(\bigoplus_{j\in J}
\E_j \biggr)
\]
acting on $ (\bigoplus_{i\in I}L^2(A_i)  )\oplus
(\bigoplus_{j\in J} L^2(B_j) )$. Since the $A_i$'s and $B_j$'s are
compact sets and $\K_N$ is locally trace class, identity (\ref
{switchidentity}) then yields
\begin{eqnarray}
\label{AIstep2} && {\fontsize{7.5}{8}\selectfont{\operatorname{det}\lleft(I-\lleft[\matrix{ \bigl[
\K_{II}^{i,j} \bigr]_{(i,j)\in I\times I} & \bigl[
\K_{IJ}^{i,j} \bigr]_{(i,j)\in I\times J}
\cr
\bigl[
\K_{JI}^{i,j} \bigr]_{(i,j)\in J\times I} & \bigl[
\K_{JJ}^{i,j} \bigr]_{(i,j)\in J\times J} } \rright] \rright)_{ (\bigoplus_{i\in I}L^2(A_i)  )\oplus
 (\bigoplus_{j\in J} L^2(B_j) )}}}\nonumber
\\
&&{\fontsize{7.5}{8}\selectfont{\qquad= \operatorname{det}\lleft(I-\E\lleft[\matrix{ \bigl[
\K_{II}^{i,j} \bigr]_{(i,j)\in I\times I} & \bigl[
\K_{IJ}^{i,j} \bigr]_{(i,j)\in I\times J}
\cr
\bigl[
\K_{JI}^{i,j} \bigr]_{(i,j)\in J\times I} & \bigl[
\K_{JJ}^{i,j} \bigr]_{(i,j)\in J\times J} } \rright]
\E^{-1} \rright)_{ (\bigoplus_{i\in I}L^2(A_i)
)\oplus
 (\bigoplus_{j\in J} L^2(B_j) )}}}
\nonumber\\[-8pt]\\[-8pt]\nonumber
&&{\fontsize{7.5}{8}\selectfont{\qquad= \operatorname{det}\lleft(I- \lleft[\matrix{ \bigl[
\E_i^* \K_{II}^{i,j}\bigl(\E_j^*
\bigr)^{-1} \bigr]_{(i,j)\in I\times I} & \bigl[ \E_i^*
\K_{IJ}^{i,j} \E_j^{-1}
\bigr]_{(i,j)\in I\times J}
\cr
\bigl[ \E_i\K_{JI}^{i,j}
\bigl(\E_j^*\bigr)^{-1} \bigr]_{(i,j)\in J\times I} & \bigl[
\E_i\K_{JJ}^{i,j} \E_j^{-1}
\bigr]_{(i,j)\in J\times J} } \rright] \rright)_{ (\bigoplus_{i\in I}L^2(A_i)  )\oplus
 (\bigoplus_{j\in J} L^2(B_j) )}}}
\nonumber
\\
&&{\fontsize{7.5}{8}\selectfont{\qquad= \operatorname{det}\lleft(I- \lleft[\matrix{ \bigl[
\mathbf 1_{A_i}\E_i^* \K_N \bigl(
\E_j^*\bigr)^{-1}\mathbf  1_{A_j}
\bigr]_{(i,j)\in I\times I}& \bigl[\mathbf 1_{A_i}\E_i^*
\K_N \E_j^{-1}\mathbf 1_{B_j}
\bigr]_{(i,j)\in
I\times J}
\cr
\bigl[\mathbf 1_{B_i}\E_i
\K_N \bigl(\E_j^*\bigr)^{-1}
\mathbf 1_{A_j} \bigr]_{(i,j)\in J\times I} & \bigl[\mathbf 1_{B_i}
\E_i \K_N \E_i^{-1}
\mathbf 1_{B_j} \bigr]_{(i,j)\in
J\times J} } \rright] \rright)_{L^2(\R)^{\oplus(\llvert  I\rrvert  +\llvert  J\rrvert  )}},}}\nonumber
\end{eqnarray}
where $\llvert  I\rrvert  $ and $\llvert  J\rrvert  $ stand for the cardinalities of $I$ and $J$, respectively.
By using definition (\ref{det2}) of $\det_2$, it follows
from~(\ref{AIstep1}) and~(\ref{AIstep2}) that
\begin{eqnarray}
\label{bidulechouette} &&{\fontsize{7.5}{8}\selectfont{\det ( I -\K_N )_{L^2 ((\bigcup_{i\in I}A_i) \cup
(\bigcup_{j\in J}B_j) )}}}\nonumber
\\
&&{\fontsize{7.5}{8}\selectfont{\qquad= \prod_{i\in I}e^{\Tr(\mathbf 1_{A_i}\E_i^* \K_N (\E
_i^*)^{-1}\mathbf 1_{A_i})} \prod
_{j\in J}e^{\Tr(\mathbf
1_{B_j}\E_j \K_N
\E_j^{-1}\mathbf 1_{B_j})}}}
\nonumber\\[-8pt]\\[-8pt]\nonumber
&&{\fontsize{7.5}{8}\selectfont{\quad\qquad{} \times{\det}_2\lleft(I- \lleft[\matrix{ \bigl[
\mathbf 1_{A_i}\E_i^* \K_N \bigl(
\E_j^*\bigr)^{-1}\mathbf  1_{A_j}
\bigr]_{(i,j)\in I\times I}& \bigl[\mathbf 1_{A_i}\E_i^*
\K_N \E_j^{-1}\mathbf 1_{B_j}
\bigr]_{(i,j)\in
I\times J}
\cr
\bigl[\mathbf 1_{B_i}\E_i
\K_N \bigl(\E_j^*\bigr)^{-1}
\mathbf 1_{A_j} \bigr]_{(i,j)\in J\times I} & \bigl[\mathbf 1_{B_i}
\E_i \K_N \E_i^{-1}
\mathbf 1_{B_j} \bigr]_{(i,j)\in
J\times J} } \rright] \rright)_{L^2(\R)^{\oplus(\llvert  I\rrvert  +\llvert  J\rrvert  )}}.}}
\nonumber
\end{eqnarray}
Let us inspect the diagonal elements of the matrix valued operator in
the Fredholm determinant at the right-hand
side of the previous identity. In Section~\ref{SectionTW}, we have
precisely shown
that for every $i\in I$ and $j\in J$,
\begin{eqnarray*}
\lim_{N\rightarrow\infty}{\bigl\llVert \mathbf 1_{A_i}
\E_i^* \K_N \bigl(\E _i^*
\bigr)^{-1}\mathbf 1_{A_i}-\mathbf 1_{(s_i,\infty)}
\K_\Aii \mathbf 1_{(s_i,\infty
)}\bigr\rrVert }_2 &=&0,
\\
\lim_{N\rightarrow\infty}{\bigl\llVert \mathbf 1_{B_j}
\E_j \K_N \E _j^{-1}
\mathbf 1_{B_j}-\mathbf 1_{(t_j,\infty)}\K_\Aii
\mathbf 1_{(t_j,\infty
)}\bigr\rrVert }_2 &=&0
\end{eqnarray*}
and
\begin{eqnarray*}
\lim_{N\rightarrow\infty}\Tr\bigl(\mathbf 1_{A_i}
\E_i^* \K_N \bigl(\E _i^*
\bigr)^{-1}\mathbf 1_{A_i}\bigr)&=&\Tr(
\mathbf 1_{(s_i,\infty)}\K _\Aii\mathbf  1_{(s_i,\infty)}),
\\
\lim_{N\rightarrow\infty}\Tr\bigl(\mathbf 1_{B_j}
\E_j \K_N \E _j^{-1}\mathbf
1_{B_j}\bigr)&=&\Tr(\mathbf 1_{(t_j,\infty)}\K_\Aii
\mathbf  1_{(t_j,\infty)}).
\end{eqnarray*}
Proposition~\ref{AIkey} then yields that the Hilbert--Schmidt norms of
the off
diagonal entries of the matrix valued operator in the Fredholm
determinant at the right-hand side
of~(\ref{bidulechouette}) converge to zero. Recalling that $\det_2$ is
continuous with respect to the Hilbert--Schmidt norm, we obtain from
(\ref{bidulechouette}) that
\begin{eqnarray*}
&& \lim_{N\rightarrow\infty}\det ( I -\K_N )_{L^2
((\bigcup_{i\in I}A_i) \cup (\bigcup_{j\in J}B_j) )}
\\[-1pt]
&&\qquad= \prod_{i\in I}e^{\Tr(\mathbf 1_{(s_i,\infty)}\K_\Aii
\mathbf
1_{(s_i,\infty)})} {
\det}_2 (I-\mathbf 1_{(s_i,\infty)}\K _\Aii\mathbf
1_{(s_i,\infty)} )_{L^2(\R)}
\\[-1pt]
&&\quad\qquad{} \times\prod_{j\in J} e^{\Tr(\mathbf
1_{(t_j,\infty
)}\K_\Aii\mathbf 1_{(t_j,\infty)})}
{\det}_2 (I-\mathbf  1_{(t_j,\infty
)}\K_\Aii
\mathbf 1_{(t_j,\infty)} )_{L^2(\R)}
\\[-1pt]
&&\qquad= \prod_{i\in I}{\det} (I-\K_\Aii
)_{L^2(s_i,\infty
)} \prod_{j\in J}{\det} (I-
\K_\Aii )_{L^2(t_j,\infty)},
\end{eqnarray*}
and Theorem~\ref{thasymptindep} is proved.
\end{pf*}

Now we turn to the proof of Proposition~\ref{AIkey}.

To do so, we shall deform the contours $\Gamma$ and $\Theta$ in the
integral representation of~$\K_N$ to appropriate contours for the
asymptotic analysis, as provided by the propositions~\ref
{contoursprop},~\ref{contoursleft+}~and~\ref{contoursleft-}. The
problem is that since $\Theta$ and $\Gamma$ will be associated to
different critical points $\frak c_N$'s or $\frak d_N$'s, the
possibility that they intersect holds true.
This raises a problem related to the presence of the factor
$(w-z)^{-1}$ in the integral representation of $\K_N$. This problem
can be
avoided by using the following alternative expression of the kernel $\K
_N$, that
was established in \cite{bleher-kuijlaars-05}; since the proof is
short, we provide it for the sake of completeness.

%
%le5.3 #&#
\begin{lemma}
\label{CNz,w}
For every $x\neq y$ we have
\begin{eqnarray}
\K_N(x,y)
&=& \frac{N}{(2i\pi)^2(x-y)}
\nonumber\\[-9pt]\\[-9pt]\nonumber
&&{} \times \oint_\Gamma\d z\oint_\Theta \d w\,
e^{-Nxz+Nyw}C_N(z,w) \biggl(\frac{z}{w}
\biggr)^N\prod_{i=1}^n \biggl(
\frac{1-\lambda_iw}{1-\lambda_iz} \biggr),\hspace*{-20pt}
\nonumber
\end{eqnarray}
where
%
%e133 #&#
\begin{equation}
C_N(z,w)=\frac{1}{zw}-\frac{1}{N}\sum
_{j=1}^n\frac{\lambda
_j^2}{(1-\lambda_j z)(1-\lambda_j w)}.
\end{equation}
\end{lemma}

\begin{pf}
Starting from (\ref{kernel}) with $q=0$ and following \cite{bleher-kuijlaars-05}, Section~3.3, we obtain by integrations by parts
\begin{eqnarray*}
x\K_N(x,y)&=&\frac{1}{(2i\pi)^2}\oint_\Gamma\d z
\oint_\Theta\d w\, e^{-Nxz+Nyw}\frac{1}{w-z} \biggl(
\frac{z}{w} \biggr)^N
\\
&&\hspace*{70pt}{}\times \prod_{i=1}^n
\biggl(\frac{1-\lambda_iw}{1-\lambda_iz} \biggr)\Biggl(\frac{1}{w-z}+\frac{N}{z}-\sum
_{j=1}^n\frac{\lambda
_j}{1-\lambda_j z} \Biggr)
\end{eqnarray*}
and
\begin{eqnarray*}
y\K_N(x,y)&=&\frac{1}{(2i\pi)^2}\oint_\Gamma\d z
\oint_\Theta\d w\, e^{-Nxz+Nyw}\frac{1}{w-z} \biggl(
\frac{z}{w} \biggr)^N
\\[-1pt]
&&\hspace*{70pt}{}\times \prod_{i=1}^n
\biggl(\frac{1-\lambda_iw}{1-\lambda_iz} \biggr)
\Biggl(\frac{1}{w-z}+\frac{N}{w}-\sum
_{j=1}^n\frac{\lambda
_j}{1-\lambda_j w} \Biggr).
\end{eqnarray*}
This provides
\begin{eqnarray*}
&& (x-y)\K_N(x,y)
\\
&&\qquad  =\frac{N}{(2i\pi)^2}\oint_\Gamma\d z
\oint_\Theta\d w\, e^{-Nxz+Nyw} \biggl(\frac{z}{w}
\biggr)^N
\\[-1pt]
&&\hspace*{104pt}{}\times \prod_{i=1}^n
\biggl(\frac{1-\lambda_iw}{1-\lambda_iz} \biggr)
\Biggl(\frac{1}{zw}-\frac{1}{N}\sum
_{j=1}^n\frac{\lambda
_j^2}{(1-\lambda_j z)(1-\lambda_j w)} \Biggr),
\end{eqnarray*}
and Lemma~\ref{CNz,w} follows.
\end{pf}

Equipped with Lemma~\ref{CNz,w}, we are now in position to prove
Proposition~\ref{AIkey}.

\begin{pf*}{Proof of Proposition~\ref{AIkey}}
Since the sets of indices $I$ and $J$ are finite by assumption, the
regularity condition provides $\varepsilon>0$ such that $\lambda
_j^{-1}\in(0,+\infty)\setminus\mathcal B$ for every $1\leq j\leq n$
and every $N$ large enough, where
\[
\mathcal B=\bigcup_{i\in I, j\in J} \bigl(B(\frak
c_i,\varepsilon )\cup B(\frak d_j,\varepsilon)
\bigr).
\]
We then set
\[
\mathcal K= \biggl( \biggl[\inf_N \frac{1}{\lambda_n},\sup
_N\frac
{1}{\lambda_1} \biggr] \Bigm\backslash\mathcal B \biggr)
\cup\{ 0\},
\]
so that $\{x\in\R\dvtx   x^{-1}\in\supp(\nu_N)\}\subset\mathcal K$
for every $N$ large enough and moreover $\{x\in\R\dvtx   x^{-1}\in\supp
(\nu)\}\subset\mathcal K$.

We start by proving (a). To do so, we essentially use the estimates
from the Section~\ref{KN1}. For any $(i,j)\in J\times J$ such that
$i\neq j$, we have
\begin{eqnarray}
\label{AIbb1}
&& {\bigl\llVert \mathbf 1_{B_i}\E_i
\K_N \E_j^{-1}\mathbf  1_{B_j}
\bigr\rrVert }_2^2
\nonumber\\[-8pt]\\[-8pt]\nonumber
&&\qquad =\int_{B_i}\int_{B_j}
\bigl(e^{-Nf_{i,N}(\frak d_{i,N})+Nx\frak d_{i,N}} \K_N(x,y) e^{Nf_{j,N}(\frak d_{j,N})-Ny\frak d_{j,N}}
\bigr)^2\,\d x\,\d y.\hspace*{-20pt}
\nonumber
\end{eqnarray}
By using Lemma~\ref{CNz,w} and performing the changes of variables
$x\mapsto N^{2/3}\delta_{i,N}(x-\frak b_{i,N})$ and $y\mapsto
N^{2/3}\delta_{j,N}(y-\frak b_{j,N}) $, we obtain
\begin{eqnarray}
\label{AIbb2} &&\int_{B_i}\int_{B_j}
\bigl(e^{-Nf_{i,N}(\frak d_{i,N})+Nx\frak
d_{i,N}}\K_N(x,y)e^{Nf_{j,N}(\frak d_{j,N})-Ny\frak d_{j,N}} \bigr)^2\,
\d x\,\d y
\nonumber\\[-8pt]\\[-8pt]\nonumber
&&\qquad =\frac{1}{\delta_{i,N}\delta_{j,N}}\int_{t_i}^{N^{2/3}\delta
_{i,N}\varepsilon}\int
_{t_j}^{N^{2/3}\delta_{j,N}\varepsilon
}\widetilde\K_N^{(\frak b_i,\frak b_j)}(x,y)^2\,\d x\,\d y,
\nonumber
\end{eqnarray}
where
%
%e134 #&#
\begin{eqnarray}
\label{defKNstrabisme}
&& \widetilde\K_N^{(\frak b_i,\frak b_j)}(x,y)\nonumber
\\
&&\qquad = \frac{N^{1/3}}{(2i\pi)^2(\frak b_{i,N}-\frak
b_{j,N}+x/(N^{2/3}\delta_{i,N})-y/(N^{2/3}\delta_{j,N}))}
\nonumber\\[-8pt]\\[-8pt]\nonumber
&&\quad\qquad{} \times \oint _\Gamma\d z\oint_\Theta\d w\,
C_N(z,w) e^{-N^{1/3}x\sklfrac{(z-\frak d_{i,N})}{\delta
_{i,N}}+N(f_{i,N}(z)-f_{i,N}(\frak d_{i,N}))}
\\
&&\hspace*{83pt}{}\times e^{N^{1/3}y\sklfrac{(w-\frak
d_{j,N})}{\delta_{j,N}}-N(f_{j,N}(w)-f_{j,N}(\frak d_{j,N}))}.
\nonumber
\end{eqnarray}
The main point here is that since $i\neq j$, there exists $C>0$
independent of $N$, $x$ and $y$ such that
%
%e135 #&#
\begin{equation}
\label{maintermAI} \biggl\llvert \frac{N^{1/3}}{(2i\pi)^2(\frak b_{i,N}-\frak
b_{j,N}+x/(N^{2/3}\delta_{i,N})-y/(N^{2/3}\delta_{j,N}))}\biggr\rrvert \leq
CN^{1/3}.
\end{equation}
Then, as in Sections~\ref{kernelofWishart}~and~\ref{contourssection},
we replace the contour
$\Gamma$ by $\Upsilon^{(0)}\cup\Upsilon^{(1)}$ where the contours
$\Upsilon^{(0)}$ and $\Upsilon^{(1)}$ are specified
by\vspace*{1pt} Proposition~\ref{contoursprop} with $\frak d_N=\frak d_{i,N}$. (If
$\Upsilon^{(0)}$ does not exist, we just
deform $\Gamma$ to $\Upsilon^{(1)}$.) Similarly, we deform the
contour $\Theta$ and
replace it with the contour $\Ttilde$ specified by
Proposition~\ref{contoursprop} with $\frak d_N=\frak d_{j,N}$. We then
deform the contours $\Upsilon^{(1)}$ and
$\Ttilde$ around the saddle points similar to Section~\ref{KN1}.
More precisely,
\[
\Upsilon^{(1)}=\Upsilon_*\cup\Upsilon^{(1)}_{\mathrm{res}}
\quad\mbox {and}\quad \Ttilde=\Ttilde_* \cup\Ttilde_{\mathrm{res}},
\]
where we introduce
\begin{eqnarray*}
\Upsilon_* & =&\bigl\{ \frak d_{i,N}+N^{-1/3}
e^{i\pi\theta}\dvtx  \theta\in [-\pi/3,\pi/3]\bigr\}\cup\bigl\{\frak
d_{i,N}+te^{\pm i\pi/3}\dvtx  t\in \bigl[N^{-1/3},\rho\bigr]\bigr
\},
\\
\Ttilde_*& =&\bigl\{\frak d_{j,N}+N^{-1/3} e^{i\pi\theta}\dvtx
\theta\in [2\pi/3,4\pi/3]\bigr\}\cup\bigl\{\frak d_{j,N}-te^{\pm i\pi/3}\dvtx
t\in \bigl[N^{-1/3},\rho\bigr]\bigr\},
\end{eqnarray*}
with $\rho$ chosen small enough so that Lemma~\ref{lemmaproperties-f-fN}(b)
applies for both $f_{i,N}$ and $f_{j,N}$. In addition,
Proposition~\ref{contoursprop} provides $K>0$ independent of $N$ such that
%
%e136 #&#
%e137 #&#
%e138 #&#
\begin{eqnarray}
\label{fastdecay1} \re\bigl(f_{i,N}(z)-f_{i,N}(\frak
d_{i,N})\bigr)&\leq& -K, \qquad z\in\Upsilon ^{(0)},
\\
\re\bigl(f_{i,N}(z)-f_{i,N}(\frak d_{i,N})\bigr)&
\leq& -K, \qquad z\in\Upsilon ^{(1)}_{\mathrm{res}}
\\
\label{fastdecay2} \re\bigl(f_{j,N}(w)-f_{j,N}(\frak
d_{j,N})\bigr)&\geq& K, \qquad w\in\Ttilde_{\mathrm{res}}.
\end{eqnarray}
Note that the contours $\Upsilon^{(0)}$ and $\Ttilde$ may now
intersect, and the contours $\Upsilon^{(1)}$ and $\Ttilde$ as well,
since the contours are
associated with different edges. This raises no problem since
$C_N(z,w)$ is
analytic on $\C\setminus\mathcal K$. More precisely, since by
construction the contours $\Upsilon^{(0)}$,
$\Upsilon^{(1)}$ and $\Ttilde$ lie inside a compact subset of $\C
\setminus\mathcal K$
which does not dependent on $N$, there exists $C'>0$ independent of $N$ such
that
%
%e139 #&#
\begin{equation}
\label{boundCNz,w} \bigl\llvert C_N(z,w)\bigr\rrvert \leq
C',\qquad z\in\Upsilon^{(0)}\cup\Upsilon ^{(1)},
\qquad w\in\Ttilde.
\end{equation}

Next, Lemma~\ref{lemmaproperties-f-fN}(b) yields
%
%e140 #&#
%e141 #&#
\begin{eqnarray}
\re\bigl(f_{i,N}(z)-f_{i,N}(\frak d_{i,N})\bigr)&
\leq& g_N''(\frak d_{i,N})\re
(z-\frak d_{i,N})^3/6+\Delta\llvert z-\frak
d_{i,N}\rrvert ^4,\nonumber
\\
\eqntext{z\in\Upsilon _*}
\\
\re\bigl(f_{j,N}(w)-f_{j,N}(\frak d_{j,N})\bigr)&
\geq& g_N''(\frak d_{j,N})\re
(w-\frak d_{j,N})^3/6-\Delta\llvert w-\frak
d_{j,N}\rrvert ^4, \nonumber
\\
\eqntext{w\in\Ttilde_*,}
\end{eqnarray}
where $\Delta>0$ is independent of $N$. We moreover assume we choose
$\rho$ small enough so that
%
%e142 #&#
\begin{equation}
\label{controlrhobis} g_N''(\frak
d_{i,N})-\rho\Delta>0,\qquad g_N''(
\frak d_{j,N})-\rho \Delta>0,
\end{equation}
for all $N$ large enough. Then, by using the same estimates as in
Sections~\ref{K0playsnorole} and~\ref{KN1}, we obtain for every
$x,y\geq s$ and $N$ large enough,
\begin{eqnarray*}
&& \int_{\Upsilon^{(0)}} e^{-N^{1/3}\sklfrac{x\re(z-\frak d_{i,N})}{
\delta_{i,N}}+N\re(f_{i,N}(z)-f_{i,N}(\frak d_{i,N}))}\llvert \d z\rrvert
\\
&&\qquad \leq
C_1e^{-C_2N+C_3N^{1/3}\sklfrac{x}{\delta_{i,N}}},
\\
&& \int_{\Upsilon^{(1)}} e^{-N^{1/3}\sklfrac{x\re(z-\frak d_{i,N})}{
\delta_{i,N}}+N\re(f_{i,N}(z)-f_{i,N}(\frak d_{i,N}))}\llvert \d z\rrvert
\\
&&\qquad \leq
\frac{C}{N^{1/3}}e^{-\vafrac{x-s}{2\delta
_{i,N}}}+C_1e^{-C_2N+C_3N^{1/3}\sklfrac{x}{\delta_{i,N}}},
\\
&& \int_{\widetilde\Theta} e^{N^{1/3}\sklfrac{y\re(w-\frak d_{j,N})}{
\delta_{j,N}}-N\re(f_{j,N}(w)-f_{j,N}(\frak d_{j,N}))}\llvert \d w\rrvert
\\
&&\qquad \leq
\frac{C}{N^{1/3}}e^{-\vafrac{y-s}{2\delta
_{j,N}}}+C_1e^{-C_2N+C_3N^{1/3}\sklfrac{y}{\delta_{j,N}}},
\end{eqnarray*}
for some $C,C_1,C_2,C_3>0$ independent on $N$ and $x,y$. Combined with
(\ref{maintermAI}) and~(\ref{boundCNz,w}), it follows from (\ref
{defKNstrabisme}) that
%
%e143 #&#
\begin{eqnarray}\label{ineqfinalAI}
\bigl\llvert \widetilde\K_N^{(\frak b_i,\frak b_j)}(x,y)
\bigr\rrvert &\leq&\frac
{C'}{N^{1/3}}e^{-(x-s)/(2\delta_{i,N})-(y-s)/(2\delta
_{j,N})}
\nonumber\\[-8pt]\\[-8pt]\nonumber
&&{} +C_1'e^{-C_2'N+C_3'N^{1/3}(\sklfrac{x}{\delta_{i,N}}+\sklfrac{y}{\delta_{j,N}})},
\end{eqnarray}
where $C',C_1',C_2',C_3'>0$ are independent on $N$ and $x,y$. Finally,
by mimicking the step~3 of the proof of Proposition~\ref{keyTW}, we obtain
\[
\lim_{N\rightarrow\infty}{\bigl\llVert \mathbf 1_{B_i}
\E_i \K_N \E _j^{-1}
\mathbf 1_{B_j}\bigr\rrVert }_2^2=0,
\]
as soon as $\varepsilon$ is small enough. We thus have proved~(a).

Concerning points (b) and (c), we proceed similarly to point (a) and use
Lemma~\ref{CNz,w} and the changes of variables
$x\mapsto N^{2/3}\sigma_{i,N}(\frak a_{i,N}-x)$ and
$y\mapsto N^{2/3}\delta_{j,N}(y-\frak b_{j,N}) $ in order to obtain
\[
{\bigl\llVert \mathbf 1_{A_i}\E_i^* \K_N
\E_j^{-1}\mathbf  1_{B_j}\bigr\rrVert
}_2^2 =\frac{1}{\sigma_{i,N}\delta_{j,N}}\int_{s_i}^{N^{2/3}\sigma
_{i,N}\varepsilon}
\int_{t_j}^{N^{2/3}\delta_{j,N}\varepsilon
}\widetilde\K_N^{(\frak a_i,\frak b_j)}(x,y)^2\,\d x\,\d y,
\]
where
%
%e144 #&#
\begin{eqnarray}
&& \widetilde\K_N^{(\frak a_i,\frak b_j)}(x,y)\nonumber
\\
&&\qquad = \frac{N^{1/3}}{(2i\pi)^2(\frak a_{i,N}-\frak
b_{j,N}-x/(N^{2/3}\sigma_{i,N})-y/(N^{2/3}\delta_{j,N}))}
\nonumber\\[-8pt]\\[-8pt]\nonumber
&&\quad\qquad{}\times \oint _\Gamma\d z\oint_\Theta\d w\,
C_N(z,w)
e^{N^{1/3}x(z-\frak c_{i,N})/ \sigma
_{i,N}-N(f_{i,N}^*(z)-f_{i,N}^*(\frak c_{i,N}))}
\\
&&\hspace*{82pt}{}\times e^{N^{1/3}y(w-\frak
d_{j,N})/ \delta_{j,N}-N(f_{j,N}(w)-f_{j,N}(\frak d_{j,N}))}.
\nonumber
\end{eqnarray}
If $\frak c_{i}>0$, then we replace the contour $\Gamma$ by the contour
$\Upsilon^{(0)}\cup\Upsilon^{(1)}$ (if $\Upsilon^{(0)}$ does not
exist, we just
deform $\Gamma$ into $\Upsilon^{(1)}$) specified by Proposition~\ref
{contoursleft+} with
$\frak c_N=\frak c_{i,N}$, and otherwise deform $\Gamma$ into
$\Upsilon$ as in
Proposition~\ref{contoursleft-}. We moreover deform the contour
$\Theta$ to obtain the contour $\Ttilde$
specified by Proposition~\ref{contoursprop} with $\frak d_N=\frak d_{j,N}$.
The same arguments as those in the proof of (a) show that
\[
\lim_{N\rightarrow\infty}{\bigl\llVert \mathbf 1_{A_i}
\E_i^* \K_N \E _j^{-1}
\mathbf 1_{B_j}\bigr\rrVert }_2^2=0.
\]
Similarly, we have
\[
{\bigl\llVert \mathbf 1_{B_i}\E_i \K_N
\bigl(\E_j^*\bigr)^{-1}\mathbf  1_{A_j}\bigr
\rrVert }_2^2 =\frac{1}{\delta_{i,N}\sigma_{i,N}}\int
_{t_i}^{N^{2/3}\delta
_{i,N}\varepsilon}\int_{s_j}^{N^{2/3}\sigma_{j,N}\varepsilon
}
\widetilde\K_N^{(\frak b_i,\frak a_j)}(x,y)^2 \,\d x\,\d y,
\]
and
\[
{\bigl\llVert \mathbf 1_{A_i}\E_i^* \K_N
\bigl(\E_j^*\bigr)^{-1}\mathbf  1_{A_j}\bigr
\rrVert }_2^2 =\frac{1}{\sigma_{i,N}\sigma_{j,N}}\int
_{s_i}^{N^{2/3}\sigma
_{i,N}\varepsilon}\int_{s_j}^{N^{2/3}\sigma_{j,N}\varepsilon
}
\widetilde\K_N^{(\frak a_i,\frak a_j)}(x,y)^2 \,\d x\,\d y,
\]
where
%
%e145 #&#
\begin{eqnarray}
\label{KNba}
&& \widetilde\K_N^{(\frak b_i,\frak a_j)}(x,y)\nonumber
\\
&&\qquad = \frac{N^{1/3}}{(2i\pi)^2(\frak b_{i,N}-\frak
a_{j,N}+x/(N^{2/3}\delta_{i,N})+y/(N^{2/3}\sigma_{j,N}))}
\nonumber\\[-8pt]\\[-8pt]\nonumber
&&\quad\qquad{}\times \oint _\Gamma\d z\oint_\Theta\d w\,
C_N(z,w)e^{-N^{1/3}x(z-\frak d_{i,N})/ \delta
_{i,N}+N(f_{i,N}(z)-f_{i,N}(\frak d_{i,N}))}
\\
&&\hspace*{82pt}{}\times e^{-N^{1/3}y(w-\frak
c_{j,N})/ \sigma_{j,N}+N(f_{j,N}^*(w)-f_{j,N}^*(\frak c_{j,N}))},
\nonumber
\end{eqnarray}
and
%
%e146 #&#
\begin{eqnarray}
\label{KNaa}
&& \widetilde\K_N^{(\frak a_i,\frak a_j)}(x,y)\nonumber
\\
&&\qquad = \frac{N^{1/3}}{(2i\pi)^2(\frak a_{i,N}-\frak
a_{j,N}-x/(N^{2/3}\sigma_{i,N})+y/(N^{2/3}\sigma_{j,N}))}
\nonumber\\[-8pt]\\[-8pt]\nonumber
&&\quad\qquad{}\times \oint _\Gamma\d z\oint_\Theta\d w\,
C_N(z,w)e^{N^{1/3}x(z-\frak c_{i,N})/ \sigma
_{i,N}-N(f_{i,N}^*(z)-f_{i,N}^*(\frak c_{i,N}))}
\\
&&\hspace*{82pt}{}\times e^{-N^{1/3}y(w-\frak
c_{j,N})/ \sigma_{j,N}+N(f_{j,N}^*(w)-f_{j,N}^*(\frak c_{j,N}))}.
\nonumber
\end{eqnarray}
For kernel (\ref{KNba}), we split the contour $\Gamma$ into
$\Upsilon^{(0)}$ and $\Upsilon^{(1)}$ where these contours are
specified by
Proposition~\ref{contoursprop} for $\frak d_N=\frak d_{i,N}$. (Again,
if $\Upsilon^{(0)}$ does not exist, we just
deform $\Gamma$ into $\Upsilon^{(1)}$.) We also
deform $\Theta$ to obtain the contour $\Ttilde$ as in
Proposition~\ref{contoursleft+} or Proposition~\ref{contoursleft-} with
$\frak c_N=\frak c_{j,N}$, depending on whether or not $\frak c_{j}>0$.
For kernel~(\ref{KNaa}), we similarly split the contour $\Gamma$ into
$\Upsilon^{(0)}$ and $\Upsilon^{(1)}$ and take these contours as in
Proposition~\ref{contoursleft+}
for $\frak c_N=\frak c_{i,N}$ if $\frak c_{i}>0$, and deform $\Gamma$
into $\Upsilon$ as in
Proposition~\ref{contoursleft-} otherwise. Moreover, $\Theta$ is
replaced by $\Ttilde$ as specified in Proposition~\ref{contoursleft+} or
Proposition~\ref{contoursleft-} with $\frak c_N=\frak c_{j,N}$
depending on
whether or not $\frak c_{j}>0$.

The same line of arguments as those in the proof of (a) then shows that
(b) and (c)
hold true, except when $\frak c_{j,N}<0$.
Indeed, in the latter case the contour $\Ttilde$ coming with
Proposition~\ref{contoursleft-} does cross by construction the set
$\mathcal K$ at a point $x_0^{-1}$ where $x_0\in\supp(\nu)$. Thus we
cannot use
bound (\ref{boundCNz,w}) anymore.

To overcome this technical point, having in mind definition (\ref{CNz,w})
of $C_N(z,\break w)$, observe that since by construction $\Upsilon^{(0)}\cup
\Upsilon^{(1)}$ or $\Upsilon$ lies in a
compact subset of $\C\setminus\mathcal K$, the map
$z\mapsto(1-z\lambda_\ell)^{-1}$ is bounded there uniformly in
$1\leq\ell\leq n$ and $N$ large enough. Since moreover by
construction $\Ttilde$ lies in $\C\setminus\{0\}$,
the map $(z,w)\mapsto(zw)^{-1}$ is bounded on the contours
uniformly in $N$ large enough. Observe furthermore that for every
$1\leq\ell\leq n$, we have
\[
\frac{e^{Nf_{j,N}^*(w)}}{1-\lambda_\ell w}=e^{N f_{j,N}^{*[\ell]}(w)},
\]
where
%
%e147 #&#
\begin{equation}
\label{fjrightvertell} f_{j,N}^{*[\ell]}(w)=\frak a_{j,N}(w-
\frak c_{j,N})-\log(w)+\frac
{1}{N}\mathop{\sum
_{k=1}}_{k\neq\ell}^n\log(1-
\lambda_k w).
\end{equation}
Namely,\vspace*{1pt} the pole at $w=\lambda_\ell$ introduced by $ C_N(z,w)$ is
actually canceled by $e^{Nf_{j,N}^*(w)}$. Thus items (b) and (c) of the
proposition follow provided that the previous estimates continue to
hold, uniformly in $1\leq\ell\leq n$, after the replacement of
$e^{Nf_{j,N}^*}$ by $e^{N f_{j,N}^{*[\ell]}}$. However, this is not
hard to obtain because, as a consequence of definitions (\ref{fj*})
and (\ref{fjrightvertell}), for every $k\in\N$ and compact subset
$B\subset\C\setminus\mathcal K$, there exists $C_{B,k}>0$
independent of $N$ such that
\[
\sup_{w\in B}\max_{1\leq\ell\leq n}\bigl\llvert
\bigl(f_{j,N}^{*[\ell
]} \bigr)^{(k)}(w)-
\bigl(f_{j,N}^{*} \bigr)^{(k)}(w)\bigr\rrvert \leq
\frac
{C_{B,k}}{N}.
\]
The proof of Proposition~\ref{AIkey} is therefore complete.
\end{pf*}

%s6 #&#
\section{Proof of Theorem \texorpdfstring{\protect\ref{thBessel}}{5}: Fluctuations at the hard edge}
\label{Besselsection}

In this section, we provide a proof for Theorem~\ref{thBessel}.

Let us fix $s>0$ and $\alpha\in\mathbb Z$. We set $n=N+\alpha$ and
define $\sigma_N$ as in (\ref{varBessel}). The representation for the
gap probabilities of determinantal point processes as Fredholm
determinants yields
\[
\p \bigl(N^2\sigma_N x_{\min}\geq s \bigr)=\det
(I-\K_N )_{L^2(0,s/(N^2\sigma_N))},
\]
where
\[
x_{\min}= \cases{ x_1=\tilde x_{\alpha+1}, &\quad if
$\alpha\geq0$,
\cr
x_{1-\alpha}=\tilde x_{1}, &\quad if $
\alpha<0$.}
\]
If we introduce the integral operator $\widetilde\K_N$ acting on
$L^2(0,s)$ with kernel
%
%e148 #&#
\begin{equation}
\label{kernelhardedge} \widetilde\K_N(x,y)=\frac{1}{N^2\sigma_N}
\K_N \biggl(\frac
{x}{N^2\sigma_N},\frac{y}{N^2\sigma_N} \biggr),
\end{equation}
then it follows from a change of variables that
%
%e149 #&#
\begin{equation}
\label{gapforBessel} \p \bigl(N^2\sigma_N x_{\min}
\geq s \bigr)=\det (I-\widetilde\K _N )_{L^2(0,s)}.
\end{equation}
We\vspace*{1pt} recall that $\K_{\mathrm{Be},\alpha}(x,y)$ has been introduced in
(\ref{Besselkernel}) and also define the operator $\E$ and $\E^{-1}$
acting\vspace*{1pt} on $L^2(0,s)$ by $\E h(x)=x^{\alpha/2}h(x)$ and $\E^{-1}
h(x)=x^{-\alpha/2}h(x)$. Notice that when $\alpha\geq0$ (resp.,
$\alpha<0$), the operator $\E$ (resp., $\E^{-1}$) is well defined on
$L^2(0,s)$, but $\mathrm{E}^{-1}$ (resp., $\mathrm{E}$) is not defined on the whole space.
Nevertheless, in the following these operators will always arise
pre-multiplied or post-multiplied by an appropriate operator so that
the product is well defined on $L^2(0,s)$; see below.

The aim of this section is to prove the following.

%
%pr6.1 #&#
\begin{proposition}
\label{keyBessel}
\[
\lim_{N\rightarrow\infty} \sup_{(x,y)\in(0,s]\times(0,s]}\bigl\llvert
\widetilde\K_N(x,y)-\E\K_{\mathrm{Be},\alpha}\E^{-1}(x,y)\bigr
\rrvert =0.
\]
\end{proposition}

Let us first show how Theorem~\ref{thBessel} follows from this proposition.

\begin{pf*}{Proof of Theorem~\ref{thBessel}}
The relation $x J'_\alpha(x)=\alpha J_\alpha(x)-xJ_{\alpha+1}( x)$
(see \cite{book-erdelyi-1953}, Section~7.2.8, equation (54)) provides
%
%e150 #&#
\begin{equation}
\label{Besselvariation} \K_{\mathrm{Be},\alpha}(x,y)=\frac{\sqrt x J_{\alpha+1}(\sqrt x
)J_\alpha(\sqrt y ) - \sqrt y J_{\alpha+1}(\sqrt y )J_\alpha(\sqrt
x )}{2(x-y)}.
\end{equation}
It then follows from \cite{book-erdelyi-1953}, Section~7.14.1, equation (9), that
\[
\K_{\mathrm{Be},\alpha}(x,y)=\frac{1}{4}\int_0^1
J_\alpha(\sqrt{x u} )J_\alpha(\sqrt{y u} )\,\d u,
\]
and, after the change of variables $u\mapsto u/s$, this yields the
factorization $\K_{\mathrm{Be},\alpha}=\B_s^2$ as operators of
$L^2(0,s)$ where $\B_s$ has for kernel $\B_s(x,y)=J_\alpha(\sqrt
{xy/s})/(2\sqrt s)$. The asymptotic behavior as $x\to0$
\begin{eqnarray*}
J_\alpha(\sqrt x ) &=&\frac{1}{\alpha!} \biggl(\frac{\sqrt x}{2}
\biggr)^\alpha \bigl(1+O\bigl(x^2\bigr) \bigr), \qquad
\mbox{if } \alpha\geq0,
\\
J_\alpha(\sqrt x ) &=&\frac{(-1)^\alpha}{\llvert  \alpha \rrvert  !} \biggl(\frac
{\sqrt x}{2}
\biggr)^{\llvert  \alpha \rrvert  } \bigl(1+O\bigl(x^2\bigr) \bigr), \qquad\mbox
{if } \alpha< 0,
\end{eqnarray*}
which\vspace*{1pt} is provided by the series representation (\ref{seriesrepBessel})
of $J_\alpha$, then shows that $\B_s$, $\B_s\E^{-1}$ and $\K_{\mathrm{Be},\alpha}\E^{-1}$ when $\alpha\geq0$, $\E\B_s$ and $\E\K
_{\mathrm{Be},\alpha}$ when $\alpha<0$, and $\E\K_{\mathrm{Be},\alpha}\E
^{-1}$ are well defined and Hilbert--Schmidt operators. Moreover, $\E$
and $\K_{\mathrm{Be},\alpha}\E^{-1}$ when $\alpha\geq0$, $\E^{-1}$
and $\E\K_{\mathrm{Be},\alpha}$ when $\alpha<0$, and $\E\K_{\mathrm{Be},\alpha} \E^{-1}$ are trace class being products of two
Hilbert--Schmidt operators.

Since $[0,s]$ is compact, it follows from Proposition~\ref{keyBessel} that
\[
\lim_{N\rightarrow\infty}{\bigl\llVert \mathbf 1_{(0,s)} \bigl(
\widetilde\K _N -\E\K_{\mathrm{Be},\alpha}\E^{-1} \bigr)
\mathbf 1_{(0,s)} \bigr\rrVert }_2=0
\]
and
\[
\lim_{N\rightarrow\infty}\Tr (\mathbf 1_{(0,s)} \widetilde
\K_N \mathbf 1_{(0,s)} )=\Tr \bigl(\mathbf 1_{(0,s)}
\E\K_{\mathrm{Be},\alpha}\E ^{-1}\mathbf 1_{(0,s)} \bigr).
\]
We then obtain from Proposition~\ref{propdet2} that
\[
\lim_{N\rightarrow\infty}\det (I-\widetilde\K_N
)_{L^2(0,s)}=\det \bigl(I-\E\K_{\mathrm{Be},\alpha}\E^{-1}
\bigr)_{L^2(0,s)},
\]
which shows together with (\ref{gapforBessel}) and (\ref
{switchidentity}) that
\[
\lim_{N\rightarrow\infty}\p \bigl(N^2\sigma_N
x_{\min}\geq s \bigr)=\det (I-\K_{\mathrm{Be},\alpha} )_{L^2(0,s)}.
\]
Finally,\vspace*{1pt} $\det (I-\K_{\mathrm{Be},0} )_{L^2(0,s)}=e^{-s}$ has
been observed in \cite{forrester-1993-bessel}, and the proof of
Theorem~\ref{thBessel} is complete.
\end{pf*}

We now focus on the proof of Proposition~\ref{keyBessel}.

%s6.1 #&#
\subsection{The Bessel kernel}

We first provide a double complex integral formula for the Bessel kernel.

%
%le6.2 #&#
\begin{lemma}
\label{doubleintBessel}
With $\K_{\mathrm{Be},\alpha}(x,y)$ defined in (\ref{Besselkernel}),
for every $0<r<R$ and \mbox{$x,y>0$}, we have
%
%e151 #&#
\begin{eqnarray}
&& \K_{\mathrm{Be},\alpha}(x,y)\nonumber
\\
&&\qquad =\frac{1}{(2i\pi)^2} \biggl(\frac
{y}{x} \biggr)^{\alpha/2}
\\
&&\quad\qquad{}\times \oint_{\llvert  z\rrvert  =  r} \frac{\d z}{z}\oint _{\llvert  w\rrvert  =  R}\frac{\d w}{w}
\frac{1}{z-w} \biggl(\frac{z}{w} \biggr)^\alpha
e^{-\sklfrac{x}z+\sklfrac{z}4+\sklfrac{y}w-\sklfrac{w}4}.
\nonumber
\end{eqnarray}
\end{lemma}

We recall that by convention, all contours of integrations are oriented
counterclockwise, and thus the notation $\oint_{\llvert  z\rrvert  =r}$ is unambiguous.

\begin{pf*}{Proof of Lemma \ref{doubleintBessel}}
The Laurent series generating function for the Bessel
functions with integer parameters reads (see \cite{book-erdelyi-1953}, 7.2.4 (25))
\[
e^{\sklfrac{x}{2}(z-\sklfrac{1}z)}=\sum_{\alpha\in\mathbb Z}J_\alpha
(x)z^\alpha,\qquad z\in\C\setminus\{0\}.
\]
This yields for every $x,r>0$ and $\alpha\in\mathbb Z$,
\[
J_{\alpha}(\sqrt x )=\frac{1}{2i\pi}\oint_{\llvert  z\rrvert  = r}z^{-\alpha
}e^{\sklfrac{\sqrt x}{2}(z-\sklfrac{1}z)}
\frac{\d z}{z}.
\]
After the changes of variables $z\mapsto2\sqrt x z$ and $w\mapsto
1/(2\sqrt yw)$, this provides for every $x,y>0$, $0<r<R$ and $\alpha
\in\mathbb Z$,
%
%e152 #&#
%e153 #&#
\begin{eqnarray}
\label{besselx} J_{\alpha}(\sqrt x ) & =&\frac{1}{2i\pi(2\sqrt x )^\alpha}\oint
_{\llvert  z\rrvert  = 1/r}z^{-\alpha}e^{xz-\sklafrac{1}{4z}}\frac{\d z}{z},
\\
\label{bessely} J_{\alpha}(\sqrt y ) & =&\frac{ (2\sqrt y )^\alpha}{2i\pi}\oint
_{\llvert  w\rrvert  =  1/R}w^{\alpha}e^{-yw+\sklafrac{1}{4w}}\frac{\d w}{w}.
\end{eqnarray}
By plugging (\ref{besselx}) and (\ref{bessely}) into (\ref
{Besselvariation}), we obtain
%
%e154 #&#
\begin{eqnarray}
\label{besselint1}
&& (x-y)\K_{\mathrm{Be},\alpha}(x,y)\nonumber
\\
&&\qquad =\frac{1}{(2i\pi)^2} \biggl(\frac{y}{x} \biggr)^{\alpha/2}
\\
&&\quad\qquad{} \times\oint
_{\llvert  z\rrvert  =  1/r}\d z\oint_{\llvert  w\rrvert  =  1/R}\d w \frac{w^\alpha}{z^{\alpha
+1}}e^{xz-1/(4z)-yw+1/(4w)}
\biggl(\frac{1}{4zw}-y \biggr).
\nonumber
\end{eqnarray}
We continue the computation by mean of integrations by parts, as
explained to us by Manuela Girotti while we discussed a similar formula
appearing in her work \cite{girotti-13-preprint}. Indeed, since
$-ye^{-yw}=\frac{\partial}{\partial w}e^{-yw}$, a first integration
by parts provides
%
%e155 #&#
\begin{eqnarray}
&& \oint_{\llvert  z\rrvert  =  1/r} \d z\oint_{\llvert  w\rrvert  =  1/R}\d w \frac{w^\alpha
}{z^{\alpha+1}}e^{xz-\sklafrac{1}{4z}-yw+\sklafrac{1}{4w}}
\biggl(\frac{1}{4zw}-y \biggr)\nonumber
\\
&&\qquad= \oint_{\llvert  z\rrvert  =  1/r} \d z\oint_{\llvert  w\rrvert  =  1/R}\d w
\frac{w^\alpha}{z^{\alpha+1}}e^{xz-\sklafrac{1}{4z}-yw+\sklafrac{1}{4w}}\nonumber
\\
&&\hspace*{117pt}{}\times \biggl(\frac{1}{4zw}+
\frac{1}{4w^2}-\frac{\alpha}{w} \biggr)
\\
&&\qquad= \oint_{\llvert  z\rrvert  =  1/r} \d z\oint_{\llvert  w\rrvert  =  1/R}\d w
\frac
{1}{z-w} \biggl(\frac{w}{z} \biggr)^\alpha
e^{xz-\sklafrac{1}{4z}-yw+\sklafrac{1}{4w}} \nonumber
\\
&&\hspace*{116pt}{}\times \biggl(\frac{1}{4w^2}-\frac{1}{4z^2}+
\frac{\alpha}{z}-\frac
{\alpha}{w} \biggr).\nonumber
\end{eqnarray}
Next, by observing that
\begin{eqnarray*}
&& \biggl(\frac{1}{4w^2}-\frac{1}{4z^2} \biggr)e^{-\sklafrac{1}{4z}+\sklafrac{1}{4w}}
\\
&&\qquad =-
\biggl(\frac{\partial}{\partial z}+\frac{\partial
}{\partial w} \biggr)e^{-\sklafrac{1}{4z}+\sklafrac{1}{4w}},
\end{eqnarray*}
another integration by parts yields
%
%e156 #&#
\begin{eqnarray}
\label{besselint2}
&& \oint_{\llvert  z\rrvert  =  1/r} \d z\oint_{\llvert  w\rrvert  =  1/R}\d w
\frac{1}{z-w} \biggl(\frac{w}{z} \biggr)^\alpha
e^{xz-\sklafrac{1}{4z}-yw+\sklafrac{1}{4w}} \nonumber
\\
&&\hspace*{83pt}
{}\times \biggl(\frac{1}{4w^2}-\frac{1}{4z^2}+
\frac{\alpha}{z}-\frac
{\alpha}{w} \biggr)
\\
&&\qquad= (x-y) \oint_{\llvert  z\rrvert  =  1/r} \d z\oint_{\llvert  w\rrvert  = 1/R}\d w
\frac{1}{z-w} \biggl(\frac{w}{z} \biggr)^\alpha
e^{xz-\sklafrac{1}{4z}-yw+\sklafrac{1}{4w}}.
\nonumber
\end{eqnarray}
By combining (\ref{besselint1})--(\ref{besselint2}), we obtain
%
%e157 #&#
\begin{eqnarray}
&& \K_{\mathrm{Be},\alpha}(x,y)\nonumber
\\
&&\qquad =\frac{1}{(2i\pi)^2} \biggl(\frac
{y}{x} \biggr)^{\alpha/2}
\\
&&\quad\qquad{}\times  \oint_{\llvert  z\rrvert  = 1/r} \d z\oint_{\llvert  w\rrvert  = 1/R}\d w \frac{1}{z-w} \biggl(
\frac{w}{z} \biggr)^\alpha e^{xz-\sklafrac{1}{4z}-yw+\sklafrac{1}{4w}},
\nonumber
\end{eqnarray}
and the lemma follows after the change of variables $z\mapsto-1/z$ and
$w\mapsto-1/w$.
\end{pf*}

%
%co6.3 #&#
\begin{corollary}
\label{Besselkernelmodified}
For every $0<r<R$ and $x,y>0$, we have
\begin{eqnarray*}
&& \E\K_{\mathrm{Be},\alpha}\E^{-1}(x,y)
\nonumber\\[-8pt]\\[-8pt]\nonumber
&&\qquad =\frac{1}{(2i\pi)^2}\oint
_{\llvert  z\rrvert  =  r} \frac{\d z}{z}\oint_{\llvert  w\rrvert  =  R}\frac{\d w}{w}
\frac
{1}{z-w} \biggl(\frac{z}{w} \biggr)^\alpha
e^{-\sklfrac{x}z+\sklfrac{z}4 +\sklfrac{y}w-\sklfrac{w}4}.
\end{eqnarray*}
\end{corollary}

Equipped with Corollary~\ref{Besselkernelmodified}, we are now in
position to establish Proposition~\ref{keyBessel}.

%s6.2 #&#
\subsection{Asymptotic analysis}
We now perform an asymptotic analysis for the kernel $\widetilde\K
_N(x,y)$ as in Section~\ref{SectionTW}. The main idea is that when the
leftmost edge is a hard edge, the associated critical point $\frak c$
should be at infinity. This leads us to study the integrand of the
double integral representation of $\widetilde{\K}_N(x,y)$ in a
neighborhood of $z=0$ and $w=0$ after the changes of variables
$z\mapsto1/z$ and $w\mapsto1/w$.

\begin{pf*}{Proof of Proposition~\ref{keyBessel}}
By choosing $q=0$ in (\ref{kernel}), which is possible according to
Remark~\ref{freedomq}, we obtain with (\ref{kernelhardedge})
\begin{eqnarray}
\label{kernelhard1}
\widetilde\K_N(x,y) &=& \frac{1}{(2i\pi)^2N\sigma_N}\oint_{\Gamma} \d z \oint_{\Theta} \d w
\frac{1}{w-z} \biggl(\frac{z}{w} \biggr)^{N}e^{- \sklafrac{zx}{N\sigma_N}+\sklvafrac{wy}{N\sigma_N}}\hspace*{-20pt}
\nonumber\\[-8pt]\\[-8pt]\nonumber
&&\hspace*{93pt}{}\times \prod_{j=1}^n\frac{w-\lambda_j^{-1}}{z-\lambda_j^{-1}},
\nonumber
\end{eqnarray}
where we recall that the contour $\Gamma$ encloses the $\lambda
_j^{-1}$'s whereas the contour $\Theta$ encloses the origin and is
disjoint from $\Gamma$. We deform $\Gamma$ so that it encloses
$\Theta$, which is possible since the integrand is analytic at the
origin as a function of $z$ and the residue picked at $z=w$ vanishes.
Moreover, since the $\lambda_j^{-1}$'s are zeros of the integrand as a
function of $w$, we can deform $\Theta$ such that it encloses all the
$\lambda_j^{-1}$'s. More precisely, we specify the contours to be
$\Gamma=\{z\in\C\dvtx    \llvert  z\rrvert  =N\sigma_N/r\}$ and $\Theta=\{z\in\C\dvtx
\llvert  z\rrvert  =N\sigma_N/R\}$ with $0<r<R<\liminf_N\lambda_1/2$. Notice that
for $N$ large enough, $\Gamma= \Gamma(N)$ and $\Theta=\Theta(N)$
enclose the $\lambda_j$'s.

Next, we perform the changes of variables $z\mapsto N\sigma_N/z$ and
$w\mapsto N\sigma_N/w$ in~(\ref{kernelhard1}) in order to get
\begin{eqnarray*}
\widetilde\K_N(x,y)
&=& \frac{1}{(2i\pi)^2}\oint_{\llvert  z\rrvert  = r} \frac{\d z}{z} \oint
_{\llvert  w\rrvert  = R} \frac{\d w}{w} \frac{1}{z-w} \biggl(
\frac{z}{w} \biggr)^{\alpha}e^{- \sklfrac{x}z+\sklfrac{y}w}
\\
&&\hspace*{104pt}{}\times \prod
_{j=1}^n\frac{\sklafrac{w}{N\sigma
_N}-\lambda_j}{\sklafrac{z}{N\sigma_N}-\lambda_j}
\\
&=& \frac{1}{(2i\pi)^2}\oint_{\llvert  z\rrvert  = r} \frac{\d z}{z} \oint
_{\llvert  w\rrvert  = R} \frac{\d w}{w} \frac{1}{z-w} \biggl(
\frac{z}{w} \biggr)^{\alpha}e^{- \sklfrac{x}z+\sklfrac{y}w}
\\
&&\hspace*{71pt}\qquad\quad{} \times e^{-N(F_N(z)-F_N(0))+N(F_N(w)-F_N(0))},
\end{eqnarray*}
where use the fact that $n=N+\alpha$, and we introduce the map
\[
F_N(z)=\frac{1}{N}\sum_{j=1}^n
\log \biggl(\frac{z}{N\sigma_N}-\lambda _j \biggr).
\]
Note that for every $N$ large enough and $z\in B(0,R+1)$, we have
$\llvert  z\rrvert  /N\sigma_N\leq\liminf_N\lambda_1/2-\delta$ for some $\delta
>0$. Thus we can choose a branch of the logarithm such that $F_N$ is
well defined and holomorphic on $B(0,R+1)$ for all $N$ sufficiently
large. Moreover, recalling that
\[
\sigma_N=\frac{4}{N}\sum_{j={1}}^n
\frac{1}{\lambda_j}
\]
and observing the identity $F_N'(0)=-1/(4N)$, a Taylor expansion of
$F_N$ around zero yields for every $z\in B(0,R+1)$ and for all $N$
large enough,
\begin{eqnarray*}
\biggl\llvert F_N(z) - F_N(0)+ \frac{z}{4N}\biggr
\rrvert &\leq& \frac
{1}{2}\frac{\llvert  z\rrvert  ^2}{N^2\sigma_N^2}\sup_{w\in B(0,R+1)}
\Biggl\llvert \frac
{1}{N}\sum_{j=1}^n
\frac{1}{(\sklafrac{w}{N\sigma_N}-\lambda
_j)^2}\Biggr\rrvert
\\
&\leq& \frac{n}{2N^3\sigma_N^2\delta^2}(R+1)^2 \leq \frac
{\Delta}{N^2},
\end{eqnarray*}
for some $\Delta>0$ independent of $N$.

Finally, by using Corollary~\ref{Besselkernelmodified} and the
inequality (\ref{ineqdiffkernels}) with
\[
u=-N \bigl(F_N(z)-F_N(0) \bigr)+N \bigl(F_N(w)-F_N(0)
\bigr)\quad\mbox {and}\quad v=\frac{z-w}4, %
\]
we obtain for every $0<x,y\leq s$,
\begin{eqnarray*}
&& \bigl\llvert \widetilde\K_N(x,y)-\E\K_{\mathrm{Be},\alpha}\E
^{-1}(x,y)\bigr\rrvert
\\
&&\qquad\leq\frac{\Delta r^{\alpha-1}}{2\pi^2 {R}^{\alpha
+1}(R-r)N}
\\
&&\quad\qquad{} \times\oint_{\llvert  z\rrvert  = r} e^{- x\re(1/z)+\re(z)/4+\Delta/N}\llvert \d z\rrvert
\oint_{\llvert  w\rrvert  =  R} e^{y\re(1/w)-\re
(w)/4+\Delta/N}\llvert \d w\rrvert
\\
&&\qquad \leq \frac{C(s)}{N}
\end{eqnarray*}
for some $C(s)>0$ independent of $N$ and $0<x,y\leq s$, and Proposition
\ref{keyBessel} follows.
\end{pf*}

The proof of Theorem~\ref{thBessel} is therefore complete.

%%%%%%%%%%%%%%%%%%%%%%%%%%%%%%%%%%%%%%%%%%%%%%%%%%%%%%%%%%%%%%%
% LES ANNEXES
%%%%%%%%%%%%%%%%%%%%%%%%%%%%%%%%%%%%%%%%%%%%%%%%%%%%%%%%%%%%%%%

%sA #&#
\begin{appendix}
%sB #&#
\section{Proof of Proposition \texorpdfstring{\protect\ref{1111}}{2.4}}\label{anx-SB-gamma}

The proof of Proposition~\ref{1111} makes use of
\cite{silverstein-choi-1995}, Theorems~4.3 and 4.4. In words,
\cite{silverstein-choi-1995}, Theorem~4.3, says that on any connected
component of
$D$, there is at most one interval on which the function $g$ is decreasing,
while \cite{silverstein-choi-1995}, Theorem~4.4, says that on any two disjoint
open intervals of $D$ where $g$ is decreasing, the images of the
closures of
these intervals by $g$ are disjoint.

\begin{pf*}{Proof of Proposition~\ref{1111}}
Let us prove (a). Assume $\gamma> 1$. Since $m(z)$ is the Cauchy--Stieltjes
transform of a probability measure supported by $[0,+\infty)$, the
function $m(x)$
decreases from zero as $x$ increases from $-\infty$ to the origin.
Hence its inverse $g(x)$
decreases to $-\infty$ as $x$ increases to zero.
Since
\[
x g(x) = 1 + \gamma\int\frac{x\lambda}{1-x\lambda} \nu(\d\lambda ),
\]
the dominated convergence theorem implies that
$x g(x) \to1 - \gamma< 0$ as $x \to-\infty$. It results that
$g(x) \to0^+$ as $x \to-\infty$, and $g(x)$ reaches a positive
maximum on
$(-\infty, 0)$.
By~\cite{silverstein-choi-1995}, Theorems~4.3 and 4.4, we obtain that
the function $g(x)$ exhibits the
behavior described in the statement, and its maximum coincides with
$\frak a$.

To prove (b), recalling the expression of $x g(x)$ and observing that
\[
x^2 g'(x) = - 1 + \gamma\int \biggl(
\frac{x\lambda}{1-x\lambda} \biggr)^2 \nu(\d\lambda),
\]
we deduce that when $\gamma\leq1$, the function $g$ is negative and
decreasing on $(-\infty, 0)$.

We now show (c). For $x > 2/\eta$ and $\lambda\in\supp(\nu)$, we
have $\llvert   1 - x \lambda\rrvert   \geq x\eta- 1 > 1$. Therefore, $g(x) \to0$
and $x^2 g'(x) \to\gamma- 1 < 0$ as $x\to+\infty$ by the dominated
convergence theorem. This shows that $g(x)$ has a positive supremum on
$(1/\eta, \infty)$, and it decreases to zero as $x\to+\infty$.
By~\cite{silverstein-choi-1995}, Theorems~4.3 and 4.4, we obtain that
the function $g(x)$ exhibits the
behavior described in the statement, and its supremum coincides with
$\frak a$.

Turning to (d), assume that $[\frak d,\infty) \subset D$. Then by
Proposition~\ref{endpoints}(a) there exists $\varepsilon> 0$ such that
$g'(x) < 0$ on $(\frak d - \varepsilon, \frak d)$ and
$g'(\frak d) = 0$. It is furthermore clear that $g(x) \to0$ as $x\to
\infty$.
Since $\frak b = g(\frak d) > 0$, we get that there exists an interval
in $(\frak c, \infty)$ over which $g$ is decreasing. However, this contradicts
\cite{silverstein-choi-1995}, Theorem~4.3.

To show (e) we observe that $m(x)$, being the Cauchy--Stieltjes
transform of a probability measure, decreases from
$\frak d = \lim_{x\downarrow\frak b} m(x)$ to $0$ as $x$ increases
over the
interval $(\frak b, \infty)$. Proposition~\ref{SC-support} shows then that
$g$ decreases from $+\infty$ to $\frak b$ as $x$ increases from zero to
$\frak d$, and that $(0,\frak d) \subset(0, 1/\xi)$. Theorem~4.3 of
\cite{silverstein-choi-1995} shows that $g$ decreases nowhere on
$(\frak d, 1/\xi)$.
\end{pf*}

%sC #&#
\section{Deformed Tracy--Widom fluctuations}\label{appBBP}

In this section, we consider a particular case of a \textit{nonregular
positive edge} where our previous analysis still applies. In this case,
the fluctuations of the
associated extremal eigenvalue will be described by the deformed
Tracy--Widom law, as introduced in Baik et al. \cite{BBP-2005}, equation~(17).
Consider the integral operator $\K_{\Aii}^{(k)}$ with kernel
%
%eC.1 #&#
\begin{equation}
\label{defdeformed-TW} \K_{\Aii}^{(k)}(x,y)=\frac{1}{(2i\pi)^2}
\oint_{\Xi}\d z\oint_{\Xi
'}\d w \frac{1}{w-z} \biggl(
\frac{w}{z} \biggr)^ke^{-xz+\sklfrac{z^3}3+yw-\sklfrac{w^3}3},\hspace*{-30pt}
\end{equation}
where the contours $\Xi$ and $\Xi'$ are the same as in the proof of
Lemma~\ref{Airydoubleintegral}, and the associated
distribution\footnote{Notice that definition (\ref{defdeformed-TW})
is consistent with that given in \cite{BBP-2005}, as the
product of the operators associated with
\cite{BBP-2005}, equations~(120)~and~(122), has kernel $\K_\Aii
^{(k)}(x,y)$. }
\[
F_k(s)=\det \bigl(I-\K_\Aii^{(k)}
\bigr)_{L^2(s,\infty)}.
\]
If $k=0$, we recover the usual Airy kernel~(\ref{Airykerneldoublecontours}).\vadjust{\goodbreak}

Given a right edge $\frak b$ associated to the limiting spectral
distribution $\mu(\gamma,\nu)$, we assume the following structure
for $\nu_N$, which readily implies that $\frak b$ is a nonregular
edge for $k\geq1$:

%as3 #&#
\begin{assumption}[(Population eigenvalues at critical $\frak d$)]\label{assdeformed}
Let $k$ be a fixed integer such that there exist
eigenvalues $\zeta_1,\ldots,\zeta_k\in\{\lambda_1,\ldots,\lambda
_n\}$
satisfying ${\zeta_j}^{-1}\to\frak d$ as $N\to\infty$ for every
$1\leq j\leq k$.
\end{assumption}

The following statement may deserve a more
formal status, but since we only sketch its proof and do not provide
the full details, we simply call it a statement.

\begin{statement*} Let Assumptions~\ref{assgauss}~and~\ref{assnu}
hold true, let $\frak b$ be a right edge and
$\frak b=g(\frak d)$ with
$\frak d\in D$ and assume moreover that Assumption~\ref{assdeformed}
holds true.\footnote{In case of a positive left edge, one will
consider instead the straightforward counterpart of Assumption~\ref{assdeformed}.} Denote by
\begin{eqnarray*}
\check{\nu}_N & =& \frac{n}{n-k} \Biggl(\nu_N-
\frac{1}{n}\sum_{j=1}^k
\delta_{\zeta_j} \Biggr),\qquad\gamma_N = \frac
{n-k}N,
\\
g_N(z) & =& \frac{1}z +{ \gamma}_N \int
\frac{\lambda}{1 - z\lambda} \check{\nu}_N(\d\lambda),
\end{eqnarray*}
and let ${\frak d}_N$ and ${\tilde x}_{\phi(N)}$ be the sequences
associated to $\check{\nu}_N$ and
${g}_N$, as provided in Proposition~\ref{propproperties-regular-left}.
Assume moreover that
%
%eC.2 #&#
\begin{equation}
\label{speedcondition} \lim_{N\rightarrow\infty}N^{1/3}\max
_{j=1}^k\bigl\llvert {\zeta_j}^{-1}-{
\frak d}_N\bigr\rrvert =0
\end{equation}
and that the following weak regularity condition holds true:
%
%eC.3 #&#
\begin{equation}
\label{weakregularity} \liminf_{N\to\infty}\min_{j=1,\ldots,n,  \lambda_j\neq\zeta
_1,\ldots,\zeta_k}
\bigl\llvert \frak d-\lambda_j^{-1}\bigr\rrvert >0.
\end{equation}
Then, for every $s\in\R$,
%
%eC.4 #&#
\begin{equation}
\label{toprovekTW} \lim_{N\to\infty}\p \bigl( N^{2/3}{
\delta}_N ({\tilde x}_{\phi(N)}-{\frak b}_N )\leq s
\bigr)= \det \bigl(I-\K_\Aii^{(k)} \bigr)_{L^2(s,\infty)},
\end{equation}
where ${\frak b}_N= {g}_N({\frak d}_N)$ and ${\delta}_N=(2/
{g}_N''({\frak d}_N))^{1/3}$.
\end{statement*}

\begin{pf*}{Outline of proof for the statement}
Introducing the map
\[
{f}_N(z)=-{\frak b}_N(z-{\frak d}_N)+
\log(z)-\frac{n-k}{N}\int\log (1-xz)\check{\nu}_N(\d x),
\]
which is the counterpart of $f_N$ from Section~\ref{SectionTW}.
From Proposition~\ref{asymptoticFredholm} and a change of variables,
we have as $N\to\infty$,
\begin{eqnarray*}
&& \p \bigl( N^{2/3}\delta_N (\tilde x_{\phi(N)}-\frak
b_N )\leq s \bigr)
\\
&&\qquad = \det (I-\mathbf 1_{(s, \varepsilon N^{2/3}\delta
_N)}\widetilde{\K }_N
\mathbf 1_{(s, \varepsilon N^{2/3}\delta_N)} )_{L^2(s,\infty)}+o(1),
\end{eqnarray*}
where the integral operator $\widetilde{\K}_{N}$ is associated with
the kernel
%
%eC.5 #&#
\begin{eqnarray}
&& \widetilde\K_{N}(x,y)\nonumber
\\
&&\qquad = \frac{N^{1/3}}{(2i\pi)^2\delta_N}\nonumber
\\
&&\quad\qquad{}\times \oint _{\Gamma} \d z
\oint_{\Theta} \d w \frac{1}{w-z}
\\
&& \hspace*{82pt}{}\times \prod_{j=1}^k
\biggl(\frac{w-{\zeta_j}^{-1}}{z-{\zeta_j}^{-1}} \biggr)\nonumber
\\
&&\hspace*{108pt}{}\times e^{- N^{1/3}x\sklfrac{(z-\frak d_N)}{\delta_N}+N^{1/3}y\sklfrac
{(w-\frak d_N)}{\delta_N}+N f_N(z)-N f_N(w)}.
\nonumber
\end{eqnarray}
By following the proof of Lemma~\ref{lemmaproperties-f-fN}, we can see
that $\re f_N$ similarly converges locally uniformly toward (\ref
{defRef}) on an appropriate subset of the complex plane containing
$\frak d$, and this yields the existence of appropriate contours as in
Proposition~\ref{contoursprop} by using the same exact proof. Since by
assumption the ${\zeta_j}^{-1}$'s stay in an arbitrary small
neighborhood of $\frak d$ for every $N$ large enough, the product over
the $\zeta_j$'s in the integrand $\widetilde\K_{N}(x,y)$ is bounded
away from that neighborhood. As a consequence, we can show, as in
Section~\ref{secasmptotic-analysis-right} and in step~2 of the proof
of Proposition~\ref{keyTW}, that with $\Upsilon_*$ and $\Ttilde_*$,
respectively, defined in (\ref{defupsilon-star-modified}) and (\ref
{deftheta-modified}),
%
%eC.6 #&#
\begin{eqnarray}
\label{KNapproxkTW}
\hspace*{-3pt}&& \widetilde\K_{N}(x,y)\nonumber
\\
\hspace*{-3pt}&&\qquad =\frac{N^{1/3}}{(2i\pi)^2\delta_N}\nonumber
\\
\hspace*{-3pt}&&\qquad\quad{}\times  \oint
_{\Upsilon_*} \d z \oint_{\Ttilde_*} \d w \frac{1}{w-z}
\\
\hspace*{-3pt}&&\hspace*{91pt}{}\times \prod
_{j=1}^k \biggl(\frac{w-{\zeta_j}^{-1}}{z-{\zeta_j}^{-1}}
\biggr)\nonumber
\\
\hspace*{-3pt}&&\hspace*{116pt}{}\times  e^{- N^{1/3}x\sklfrac{(z-\frak d_N)}{\delta_N}+N^{1/3}y\sklfrac
{(w-\frak d_N)}{\delta_N}+N f_N(z)-N f_N(w)}
\nonumber
\end{eqnarray}
up to negligible terms, in the sense that the remaining terms do not
contribute in the large $N$ limit. Moreover, by proceeding similarly as
in Lemma~\ref{Airydoubleintegral} and step~2\vadjust{\goodbreak} of the proof of
Proposition~\ref{keyTW}, we have that
\begin{eqnarray}
\label{KAiapproxkTW}
\K_{\Aii}^{(k)}(x,y) &=& \frac{N^{1/3}}{(2i\pi)^2\delta_N}\nonumber
\\
&&{}\times \oint
_{\Upsilon_*}\d z\oint_{\Ttilde_*}\d w \frac{1}{w-z} \biggl(
\frac{w-\frak d_N}{z-\frak d_N} \biggr)^k\nonumber
\nonumber\\[-8pt]\\[-8pt]\nonumber
&&\hspace*{58pt}{}\times \exp \biggl\{-N^{1/3}x\frac{(z-\frak d_N)}{\delta_N} + Ng_N''(\frak
d_N)\frac{(z-\frak d_N)^3}6
\\
&&\hspace*{91pt}{}  +N^{1/3}y\frac{(w-\frak d_N)}{\delta_N}
-Ng_N''(\frak d_N)\frac{(w-\frak d_N)^3}6\biggr\}
\nonumber
\end{eqnarray}
up to negligible terms. Finally, to conclude we need to estimate the
difference between the right-hand sides of (\ref{KNapproxkTW}) and
(\ref{KAiapproxkTW}), which is the counterpart of step~1 in the proof
of Proposition~\ref{keyTW}; we claim that similar estimates can be
performed with minor modifications, provided that (\ref
{speedcondition}) holds true.
\end{pf*}
\end{appendix}

% zodis "Acknowledgments" paliekamas pagal autoriu
\section*{Acknowledgments}
Adrien Hardy and Jamal Najim are pleased to thank the organizers of the 2011 France--China
summer school in Changchun ``\textit{Random Matrix Theory and
High-dimensional Statistics}'' where this project began. The authors
are indebted to Steven Delvaux for providing an important argument in
the asymptotic analysis; see Section~\ref{contourssection}. Moreover,
Adrien Hardy would like to thanks Sandrine P\'ech\'e for interesting discussions,
and Manuela Girotti for generously sharing her computations on the
double integral representation for the Bessel kernel; see the proof of
Lemma~\ref{doubleintBessel}.

%\begin{supplement}[id=suppA]
%\sname{Supplement A}
%\stitle{}
%\slink[doi]{10.1214/00-AOPXXXXSUPP} %[doi,text={...}] - jei reikia
%suskaldyti doi
%\sdatatype{.pdf}
%\sfilename{aopXXXX\_supp.pdf}
%\sdescription{}
%\end{supplement}

% imsref loaded by linak, 2015-05-21 11:34:46
% imsref loaded by linak, 2015-05-21 13:09:15
%
% imsref loaded by linak, 2015-05-26 09:38:42
% imsref loaded by linak, 2016-03-16 10:14:51
% imsref loaded by linak, 2016-03-16 10:18:33
% imsref loaded by linak, 2016-03-16 10:22:22
% imsref loaded by linak, 2016-03-16 10:24:25
% imsref loaded by linak, 2016-03-16 10:27:25
% imsref loaded by linak, 2016-03-16 10:43:56
% imsref loaded by linak, 2016-03-16 10:53:02
% imsref loaded by linak, 2016-03-16 10:53:30

\printaddresses
\end{document}